\documentclass[a4paper,12pt]{article}

\usepackage {amssymb,latexsym}
\usepackage {amsmath}
\usepackage[french]{babel}

\newtheorem{theorem}{Th\'{e}or\`{e}me}[section]
\newtheorem{prop}{Proposition}[section]
\newtheorem{lemma}{Lemme}[section]
\newtheorem{cor}{Corollaire}[section]

\newtheorem{defi}{D\'{e}finition}[section]
\newtheorem{exer}{Exercice}[subsection]
\newtheorem{exem}{Exemple}[section]
\newtheorem{rem}{Remarque}[section]
\newtheorem{prob}{Probl\`{e}me}[section]
\newtheorem{conj}{Conjecture}
\newenvironment{prooof}{
                        \noindent{\bf\small D\'{e}monstration: }\small}
                                       {\hfill {$\mathbf \Box$}\medskip}
\newcommand{\bdefi}{\begin{defi}}
\newcommand{\edefi}{\end{defi}}
\newcommand{\bexer}{\begin{exer}\small\rm}
\newcommand{\eexer}{\end{exer}}
\newcommand{\bexem}{\begin{exem}\rm}
\newcommand{\eexem}{\end{exem}}
\newcommand{\bsat}{\begin{theorem}}
\newcommand{\esat}{\end{theorem}}
\newcommand{\bprop}{\begin{prop}}
\newcommand{\eprop}{\end{prop}}
\newcommand{\bcor}{\begin{cor}}
\newcommand{\ecor}{\end{cor}}
\newcommand{\blem}{\begin{lemma}}
\newcommand{\elem}{\end{lemma}}
\newcommand{\brem}{\begin{rem}}
\newcommand{\erem}{\end{rem}}
\newcommand{\bbew}{\begin{prooof}}
\newcommand{\ebew}{\end{prooof}}
\newcommand{\bprob}{\begin{prob}}
\newcommand{\eprob}{\end{prob}}
\newcommand{\bconj}{\begin{conj}}
\newcommand{\econj}{\end{conj}}

\newcommand{\beq}{\begin{equation}}
\newcommand{\bea}{\begin{eqnarray}}
\newcommand{\eea}{\end{eqnarray}}
\newcommand{\beas}{\begin{eqnarray*}}
\newcommand{\eeas}{\end{eqnarray*}}

\newcommand{\real}{\mathbb R}

\newcommand{\complex}{\mathbb C}
\newcommand{\nat}{\mathbb N}

\newcommand{\lbl}{\label}
\newcommand{\ben}{\begin{enumerate}}
\newcommand{\een}{\end{enumerate}}
\newcommand{\ra}{\rightarrow}

\newcommand{\Cinf}{\mathcal{C}^\infty}

\newcommand{\Ginf}{{\Gamma}^\infty}

\newcommand {\id} {{\mathrm{id}}}

\newcommand{\Dop}{\mathbf{D}}

\makeatletter
\@addtoreset{equation}{section}

\makeatother

\begin{document}

\pagestyle{empty}
~~~~~~~~~~~~
\vspace{0cm}
\begin{center}
  {\LARGE \bf Sur l'existence d'une prescription d'ordre
    naturelle projectivement invariante} \\[4mm]

\vspace{1cm}

       {\bf \large Martin Bordemann} \\
        Laboratoire des Math\'{e}matiques et Applications\\
        Facult\'{e} des Sciences et Techniques\\
        Universit\'{e} de Haute Alsace, Mulhouse \\
        4, rue des Fr\`{e}res Lumi\`{e}re\\
        68093 Mulhouse, France \\
        e--mail: M.Bordemann@univ-mulhouse.fr \\[4mm]
         Ao\^{u}t 2002

\end{center}

\vspace{1cm}

\begin{center}
 \begin{minipage}{12cm}
  {\small
   \begin{center}
   {\bf English Abstract}
   \end{center}
   For the space of all differential operators mapping the smooth sections
   of a given vector bundle on a manifold to another one over the same
   manifold, an ordering prescription (or a symbol calculus)
   is a real linear bijection from
   the space of all principal symbols to the space of differential
   operators. Ordering prescriptions are neither unique nor natural
   with respect to local diffeomorphisms. P.Lecomte has proposed to
   take into account the covariant derivatives used to build ordering
   prescriptions for the naturality of transformation properties and
   has conjectured that there exists an natural ordering prescription for
   differential operators between density bundles which in addition is
   invariant under projective changes of the covariant derivatives.
   We prove this conjecture by constructing a projectively invariant lift
   of a torsion-free connexion to a torsion-free connexion on (the positive
   part of) the total space of the bundle of all $a$-densities for nonzero
   $a$, by lifting the symbols in a projectively invariant way by
   showing them to be in bijection to the space of all $\real^+$-equivariant
   and
   divergence-free symmetric tensor fields on the total space, and by
   using the standard ordering procedure (`all the covariant derivatives to
   the right') on the total space. For Ricci-flat manifolds we show that this
   ordering prescription coincides --with the appropiate replacements--
   with an explicit formula in $\real^m$ obtained by Duval, Lecomte and
   Ovsienko.
   }
 \end{minipage}
\end{center}

\newpage
\tableofcontents
\newpage

\pagestyle{plain}

\section*{Introduction}
  \addcontentsline{toc}{section}{Introduction}

Soit $M$ une vari\'{e}t\'{e} diff\'{e}rentiable, $E$ et $E'$ deux fibr\'{e}s vectoriels
sur $M$ et $\Dop\big(\Ginf(E),\Ginf(E')\big)$ l'espace de tous les
op\'{e}rateurs diff\'{e}rentiables de $E$ dans $E'$. Une `prescription d'ordre'
--je m'excuse d'avance pour cette
traduction provisoire de l'expression anglaise `ordering prescription'--
(ou un calcul symbolique) est une bijection $\real$-lin\'{e}aire $\rho$ de
l'espace de tous les symboles principaux, $\Ginf\big(STM\otimes
Hom(E,E'\big)$ sur $\Dop\big(\Ginf(E),\Ginf(E')\big)$ qui respecte le symbole
dans le sens que $\rho\big(\sigma(D)\big)-D$ est un op\'{e}rateur d'ordre
$k-1$ si $D$ \'{e}tait un op\'{e}rateur diff\'{e}rentiel d'ordre $k$ (et non pas d'ordre
$k-1$) et $\sigma(D)$ est son symbole principal. Si on se donne une
connexion dans le fibr\'{e} tangent, une connexion dans $E$ et une connexion
dans $E'$, et si l'on suit la
r\`{e}gle `toutes les d\'{e}riv\'{e}es covariantes \`{a} droite', on obtient une
prescription d'ordre $\rho_s[\nabla]$, appel\'{e}e `standard' (voir par exemple
\cite{BNW98} ou \cite{BNPW98}): ici, il n'y a pas de
d\'{e}riv\'{e}es des symboles; l'application $\rho_s$ est
$\Cinf(M,\real)$-lin\'{e}aire.

Une prescription d'ordre est tr\`{e}s importante pour la m\'{e}canique quantique:
une quantification associe une `observable classique'
(un symbol) \`{a} une `observable quantique' (\`{a} un op\'{e}rateur diff\'{e}rentiel).
On y utilise d'autres prescriptions d'ordre, par exemple celle de
Weyl-Moyal $\rho_W[\nabla]$ qui consiste en une sym\'{e}trisation des d\'{e}riv\'{e}es
covariantes et des symboles et qui est de la forme
$\rho_W[\nabla](A)=\rho_s[\nabla](e^{\mathsf{Div}/2}A)$, o\`{u}
$\mathsf{Div}$ d\'{e}signe la divergence covariante par rapport \`{a} $\nabla$,
voir par exemple \cite{BNW98} pour une discussion.

A l'aide d'une prescription d'ordre, on peut retirer la multiplication
(en g\'{e}n\'{e}ral) noncommutative associative des op\'{e}rateurs diff\'{e}rentiels
\`{a} l'espace des symboles: pour certaines prescriptions d'ordre et pour
le cas o\`{u} $E=E'=M\times \real$ on obtient une multiplication associative
formelle bidiff\'{e}rentielle sur l'espace des fonctions de classe $\Cinf$
\`{a} valeurs r\'{e}elles sur la vari\'{e}t\'{e} symplectique $T^*M$. La g\'{e}n\'{e}ralisation
de ces multiplications est la th\'{e}orie des star-produits, voir
\cite{BFFLS78} pour sa d\'{e}finition, voir \cite{DL83} et \cite{Fed94}
pour la d\'{e}monstration
de l'existence au cas symplectique et voir finalement \cite{Kon97b} pour
l'existence et la classification \`{a} \'{e}quivalence pr\`{e}s dans le cas plus
g\'{e}n\'{e}ral d'une vari\'{e}t\'{e} de Poisson.

Pour restreindre la multitude des diff\'{e}rentes prescriptions d'ordre, on
peut exiger l'\'{e}quivariance de l'application $\rho$ par l'action d'un groupe
ou d'une alg\`{e}bre de Lie: de cette fa\c{c}on, C.Duval, P.Lecomte et
V.Y.Ovsienko ont montr\'{e} l'existence et l'unicit\'{e} (!) d'une prescription d'ordre
sur $M=\real^m$ pour des fibr\'{e}s de densit\'{e}s qui soit \'{e}quivariante par
l'action de l'alg\`{e}bre de Lie $\mathfrak{sl}(m+1,\real)$, voir \cite{DO97},
\cite{LO99},\cite{Lec99} et \cite{Lec00}:
l'action canonique du groupe $GL(m+1,\real)$ se projette sur une action
par homographies $\Phi_g:x\mapsto gx/|gx|$ ($\forall g\in GL(m+1,\real)$,
$\forall x\in S^m$ de la sph\`{e}re $S^m$ de dimension $m$
qui contient $\real^m$ comme domaine d'une carte, par cons\'{e}quent les
g\'{e}nerateurs infinit\'{e}simaux de cette action d\'{e}finissent l'action de
$\mathfrak{sl}(m+1,\real)$ sur $\real^m$
par certains champs de vecteurs quadratiques (dans des coordonn\'{e}es appropri\'{e}es).
Le m\^{e}me r\'{e}sultat de
l'existence et de l'unicit\'{e} a \'{e}t\'{e} obtenu par cette \'{e}quipe pour l'alg\`{e}bre de
Lie $\mathfrak{o}(p+1,q+1)$ (avec $p,q\in\nat$ et $p+q=m$) agissant sur
$\real^m$ par des transformations conformes infinit\'{e}simales, voir
\cite{DLO99}. Ensuite, ces r\'{e}sultats ont \'{e}t\'{e} g\'{e}n\'{e}ralis\'{e}s \`{a} d'autres espaces
homog\`{e}nes, voir \cite{BM01}.

Pour rendre ces r\'{e}sultats globaux, Pierre Lecomte a d'abord propos\'{e} dans
plusieurs expos\'{e}s de tenir en compte la structure d'une d\'{e}riv\'{e}e covariante
$\nabla$ sur la vari\'{e}t\'{e} $M$ --bel exemple du principe de jauge ou couplage
minimal des physiciens-- pour d\'{e}finir une meilleure fa\c{c}on de
transformation
par toutes les diff\'{e}omorphismes locaux: par exemple, pour les fibr\'{e}s
triviaux $E=M\times\real=E'$, on peut demander une esp\`{e}ce de `naturalit\'{e}',
c.-\`{a}-d. la r\`{e}gle de transformation
\[
     ~~~~~~~~~~~~~~ ~~~~~~\Phi^*\big(\rho[\nabla](A)(\varphi)\big)
              = \rho[\Phi^*\nabla](\Phi^*A)\big(\Phi^*\varphi\big),
              ~~~~~~~~~~~~~~~~~~~~~~(*)
\]
quels que soient le diff\'{e}omorphisme local $\Phi$, le symbole $A$ et la
fonction de classe $\Cinf$ \`{a} valeurs r\'{e}elles sur $M$. Si on enlevait
$\nabla$, il n'y aurait aucune prescription d'ordre naturelle: \`{a} partir des
op\'{e}rateurs diff\'{e}rentiels d'ordre $2$ c'est impossible, comme un calcul
\'{e}l\'{e}mentaire en $\real^m$ le montre. La prescription
d'ordre standard $\rho_s[\nabla]$ est certainement naturelle dans le sens
$(*)$.
Mais il y en a d'autres, diff\'{e}rentes de $\rho_s$, par exemple la
prescription d'ordre de Weyl pour les demi-densit\'{e}s, comme dans
\cite{BNPW98}, p.21, eqn (6.9).\\
 Pour restreindre davantage ces
prescriptions d'ordre naturelles, Pierre Le\-comte a propos\'{e} --sans doute
inspir\'{e} par l'exemple des prescriptions d'ordre \'{e}qui\-vari\-antes par
$\mathfrak{sl}(m+1,\real)$-- de demander que $\rho$ soit invariante par
un changement de connexion $\nabla\mapsto \nabla'$, o\`{u} $\nabla'$ est
projectivement \'{e}quivalente \`{a} $\nabla$, c.-\`{a}-d. que $\nabla$ et $\nabla'$
ont les m\^{e}mes g\'{e}od\'{e}siques \`{a} r\'{e}param\'{e}trage pr\`{e}s. Il a conjectur\'{e}
l'existence d'une telle prescription d'ordre. Pour les op\'{e}rateurs
diff\'{e}rentiels d'ordre $2$ et $3$ c'\'{e}tait d\'{e}j\`{a} montr\'{e}, voir
\cite{Bou00} et \cite{Bou01}.

Le but principal de cet article est de d\'{e}montrer cette conjecture de
Lecomte de fa\c{c}on g\'{e}om\'{e}trique, c.-\`{a}-d. de construire explicitement
une prescription d'ordre $\rho_L[\nabla]$ naturelle et projectivement
invariante pour toutes les vari\'{e}t\'{e}s de dimension $m\geq 2$.

L'id\'{e}e de la construction est simple: \\
Dans l'exemple de
$\tilde{M}:=\real^{m+1}\setminus \{0\}$ vu comme fibr\'{e} trivial \`{a} fibre
type $\real^+$ sur la sph\`{e}re $M:=S^m$, on voit que toutes les bijections
lin\'{e}aires $g$ de $\real^{m+1}$ pr\'{e}servent la connexion canonique
$\tilde{\nabla}$ plate sans
torsion dans $\tilde{M}$, tandis que les homographies
induites $\Phi_g$ en g\'{e}n\'{e}ral ne pr\'{e}servent pas la connexion de Levi Civita
$\nabla$ canonique sur $S^m$. Quand j'ai
entendu parler de l'\'{e}quivariance de la prescription d'ordre de Duval,
Lecomte et Ovsienko, j'ai toute suite voulu `relever' la situation `difficile
en bas sur $S^m$' \`{a} la situation `plus simple en haut sur $\real^{m+1}$':
plus pr\'{e}cis\'{e}ment, je voulais trouver un rel\`{e}vement $GL(m+1,\real)$-\'{e}quivariant
$\varphi\mapsto\tilde{\varphi}$ des densit\'{e}s sur $S^m$ et un rel\`{e}vement
$GL(m+1,\real)$-\'{e}quivariant $A\mapsto \tilde{A}$ des symboles sur $S^m$.
Ainsi la formule
\[
    ~~~~~~~~~~~~~~~~~~~~\big(\rho_L[\nabla](A)(\varphi)\big)^{\sim}
          :=\rho_s[\tilde{\nabla}](\tilde{A})(\tilde{\varphi})
   ~~~~~~~~~~~~~~~~~~~~~~~~~(**)
\]
donnerait une prescription d'ordre $GL(m+1,\real)$-\'{e}quivariante, car la
prescription d'ordre standard sur $\real^{m+1}$ l'est \'{e}videmment vue
l'invariance de la connexion $\tilde{\nabla}$.\\
Puisque l'espace des sections, sur lequel les op\'{e}rateurs diff\'{e}rentiels
agissent, est un espace de densit\'{e}s, j'ai g\'{e}n\'{e}ralis\'{e} l'exemple, en mettant
$\tilde{M}^a$ \'{e}gal \`{a} la partie strictement positive de l'espace total
du fibr\'{e} des $a$-densit\'{e}s, $|\Lambda^mT^*M|^a$. Dans l'exemple
$a=1/(m+1)$. A l'aide de la multiplication des nombres r\'{e}els
strictement positifs --qui existe pour tout fibr\'{e} vectoriel--
$\tilde{M}^a$ est muni de la structure d'un fibr\'{e} principal \`{a} groupe
struturel $\real^+$ (qui est associ\'{e} au fibr\'{e} des rep\`{e}res lin\'{e}aires).
Par cons\'{e}quent, toute $b$-densit\'{e} $\varphi$ se rel\`{e}ve en tant que fonction
$\tilde{\varphi}$ \`{a} valeurs
r\'{e}elles \'{e}quivariante par rapport \`{a} l'action de $\real^+$. De plus,
toute connexion sans torsion $\nabla$ sur $M$ d\'{e}finit
un rel\`{e}vement horizontal des champs de vecteurs sur $M$ \`{a} des champs de
vecteurs sur $\tilde{M}^a$. En jouant avec les champs relev\'{e}s et le
champ d'Euler et en me laissant inspirer par la formule reliant la d\'{e}riv\'{e}e
covariante canonique dans $\real^{m+1}\setminus \{0\}$ et la connexion de Levi
Civita sur $S^m$, j'ai \'{e}t\'{e} surpris d'avoir construit une connexion sans
torsion $\tilde{\nabla}$ sur $\tilde{M}^a$, qui ne d\'{e}pendait que de la classe
projective de $\nabla$, qui \'{e}tait invariant par la multiplication avec
des nombres r\'{e}els strictement positifs et qui se transformait bien par
l'action des diff\'{e}omorphismes locaux. Ensuite, il fallait relever les
symboles principaux de mani\`{e}re projectivement invariante; apr\`{e}s avoir fait
de gros calculs en vain,
la consid\'{e}ration suivante m'a donn\'{e} la bonne piste:
supposons qu'on ait un tel rel\`{e}vement $A\mapsto \tilde{A}$,
alors la formule $(**)$ nous donnerait une prescription d'ordre naturelle
projectivement invariante. Mais on aurait pu prendre la prescription
d'ordre de Weyl-Moyal $\rho_W[\tilde{\nabla}]$ dans la formule $(**)$
pour avoir une autre prescription d'ordre selon Lecomte.
Alors si non seulement l'\'{e}nonc\'{e}
de l'existence dans la conjecture de Lecomte, mais encore
l'\'{e}nonc\'{e} de l'unicit\'{e}
est vrai--ce que je ne sais pas \`{a} pr\'{e}sent--, il ne faut pas que
l'op\'{e}rateur de Neumaier $N=e^{\widetilde{\mathsf{Div}}/2}$ (qui fait la
transition $\rho_s\leadsto \rho_W$) modifie r\'{e}ellement les symboles relev\'{e}s.
En particulier, on arrive ainsi \`{a} la condition
\[
~~~~~~~~~~~~~~~~~~ ~~~~~~~~ ~~~~\widetilde{\mathsf{Div}}\tilde{A}=0.
~~~~~~~~~~~~~~~~~~~~~~~~~~~~~~~~~~~~~~~~~~~~(***)
\]
Je me suis aper\c{c}u que cette condition, qui \'{e}videmment ne d\'{e}pend que de
$\tilde{\nabla}$, c.-\`{a}-d. de la
classe projective de $\nabla$, --plus une certaine
$\real^+$-\'{e}quivariance-- suffit pour d\'{e}finir le rel\`{e}vement des symboles
souhait\'{e} et de donner des formules explicites pour le cas o\`{u} une vari\'{e}t\'{e}
est munie d'une
connexion sans torsion pour laquelle la partie sym\'{e}trique du tenseur de
Ricci s'annule.

Le premier paragraphe rappelle les d\'{e}finitions des prescriptions d'ordre
et leur \'{e}quivariance par rapport aux actions de groupes et d'alg\`{e}bres de
Lie. Le deuxi\`{e}me paragraphe est consacr\'{e} \`{a} la d\'{e}finition des prescriptions
d'ordre
naturelles. Pour pouvoir traiter d'autres fibr\'{e}s que ceux de $b$-densit\'{e}s,
j'ai choisi comme cadre la th\'{e}orie des fibr\'{e}s et op\'{e}rateurs
naturels, d\'{e}velopp\'{e}e par Palais, Terng, Epstein, Thurston,
Kol\'{a}\v{r}, Michor et Slov\'{a}k: le livre \cite{KMS93} \'{e}crit par ces trois derniers
auteurs est surtout recommandable et j'ai essay\'{e} d'utiliser
leur notation. Dans l'appendice \ref{SecApp2}, j'ai \'{e}num\'{e}r\'{e} quelques
constructions pour fixer la notation. Les fibr\'{e}s naturels sont grosso modo
des fibr\'{e}s qui ont un caract\`{e}re fonctoriel, c.-\`{a}-d. qui permettent de relever
des diff\'{e}omorphismes locaux, comme
$M\mapsto TM$ ou $M\mapsto S^kTM$ ou $Q_\tau P^1M$ (le fibr\'{e} affine dont les
sections sont des connexions d'ordre $q$ sur $M$ qui se projettent sur des
connexions sans torsion dans le fibr\'{e} tangent de $M$). Un op\'{e}rateur naturel
est une collection d'applications locales (param\'{e}tr\'{e}e par toutes les vari\'{e}t\'{e}s
de dimension $m$) entre l'espace de sections du premier
fibr\'{e} naturel dans l'espace de sections du deuxi\`{e}me fibr\'{e} naturel, qui
`commutent avec {\em pull back}'.
Une prescription d'ordre
naturelle pour des op\'{e}rateurs diff\'{e}rentiels entre deux fibr\'{e}s naturels
$F$ et $F'$ est
interpr\'{e}t\'{e}e comme une collection d'op\'{e}rateurs naturels entre les fibr\'{e}s
naturels
$Q_\tau P^q\times \big(S^kT\otimes Hom(F,F')\big)$ et $Hom(J^k\circ
F,F')$, o\`{u} $J^k$ est le fibr\'{e} naturel des jets (de sections) d'ordre $k$.
Alors, je suppose la localit\'{e}  a priori, contrairement \`{a} Duval, Lecomte et
Ovsienko, qui la d\'{e}duisent de l'\'{e}quivariance.
Ensuite, je donne la construction du
fibr\'{e} principal $\tilde{M}^a$ en paragraphe 3: ici il me faut quelques
notions de la th\'{e}orie des fibr\'{e}s principaux, pour laquelle les sources
principales sont \cite{KN63}, \cite{KN69} et \cite{KMS93}. J'ai donn\'{e} un
petit aper\c{c}u de ces notions dans l'appendice
\ref{SecApp1}. Le quatri\`{e}me paragraphe rappelle la notion des connexions
projectivement \'{e}quivalentes et formule une g\'{e}n\'{e}ralisation de la conjecture
de Lecomte pour d'autres fibr\'{e}s naturels $F$ et $F'$ d'ordre 1.
Le cinqui\`{e}me paragraphe
est consacr\'{e} \`{a} la construction du rel\`{e}vement naturel projectivement
invariant $\nabla\mapsto \tilde{\nabla}$. Je d\'{e}montre que ce rel\`{e}vement
est unique: la d\'{e}monstration tr\`{e}s longue de l'unicit\'{e} dans laquelle j'utilise
la th\'{e}orie des op\'{e}rateurs naturels est envoy\'{e}e dans l'appendice
\ref{SecApp3}. Dans la derni\`{e}re partie du cinqui\`{e}me paragraphe, je discute
l'exemple clef de la sph\`{e}re $S^m$.
Dans le sixi\`{e}me paragraphe, le rel\`{e}vement naturel projectivement invariant des
symboles est trait\'{e}, en analysant l'\'{e}quation $(***)$ de divergence nulle
sur $\tilde{M}^a$:
ces symboles sont des champs de tenseurs sym\'{e}triques tensoris\'{e}s avec
des $c$-densit\'{e}s, et il s'av\`{e}re que le rel\`{e}vement par cette m\'{e}thode
n'est pas possible pour tous les
nombres r\'{e}els $c$: on doit exclure les valeurs `r\'{e}sonnantes'
\[
   c\in \left\{ {\textstyle \frac{j+m}{m+1}}~|~j\in\nat~\right\}
\]
qui figurent aussi dans \cite{Lec99} et \cite{Lec00}. Le septi\`{e}me
paragraphe
donne une r\'{e}ponse positive \`{a} la conjecture de Lecomte, en construisant
explicitement la prescription d'ordre $\rho_L$ \`{a} l'aide de l'\'{e}quation $(**)$
(th\'{e}or\`{e}me
\ref{TExPreNatProjInv}). Ensuite on d\'{e}rive une formule explicite pour
le cas o\`{u} la partie sym\'{e}trique du tenseur de Ricci s'annule. Paragraphe
8 traite les prescriptions d'ordre \'{e}quivariantes: si l'action d'un groupe
relev\'{e}e de $M$ \`{a} $\tilde{M}^a$ pr\'{e}serve la connexion relev\'{e}e $\tilde{\nabla}$,
alors $\rho_L[\nabla]$ est \'{e}quivariante (Corollaire \ref {CorEquivar}).
L'application de ce r\'{e}sultat g\'{e}n\'{e}ral \`{a} l'exemple de la sph\`{e}re $S^m$
donne ainsi une explication g\'{e}om\'{e}trique pour l'existence d'une
prescription d'ordre $GL(m+1,\real)/\real^+$-\'{e}quivariante sur $S^m$,
et donc l'existence d'une prescription d'ordre
$\mathfrak{sl}(m+1,\real)$-\'{e}quivariante sur $\real^m$. Dans la derni\`{e}re
section j'\'{e}num\`{e}re quelques probl\`{e}mes ouverts qui me semblent int\'{e}ressants.

\subsubsection*{Notations}

Les notations sont ou bien directement expliqu\'{e}es ou bien contenues
dans un des appendices. $P^rM$ d\'{e}signe le fibr\'{e} principal
de tous les jets en $0$ d'ordre $r$ inversibles des applications
$\real^m\ra M$. $P^1M$ est le fibr\'{e} principal de tous les rep\`{e}res
lin\'{e}aires. $Q_\tau P^qM$ d\'{e}signe le fibr\'{e} affine de toutes les connexions
principales de $P^qM$ qui se projettent sur le fibr\'{e} $Q_\tau P^1M$ de
toutes les connexions sans torsion dans le fibr\'{e} tangent. Pour un entier
$k$, le fibr\'{e} de tous les jets d'ordre $k$ des sections d'un fibr\'{e}
(vectoriel) $E$ sur $M$ est not\'{e} par $J^kE$.

\subsubsection*{Remerciment}

Je tiens \`{a} remercier Christian Duval, Pierre Lecomte et Valentin Ovsienko
pour m'avoir introduit au sujet par des expos\'{e}s bien pr\'{e}sent\'{e}s et des
dicussions tr\`{e}s fructueuses. Je voudrais \'{e}galement remercier Sophie
Lef\`{e}vre pour avoir jet\'{e} un coup d'oeil critique sur le fran\c{c}ais de
ce manuscrit.

\section{Prescriptions d'ordre}

Soit $M$ une vari\'{e}t\'{e} diff\'{e}rentiable de dimension $m$. Pour un entier positif
$k$ soit $S^kTM$ (resp. $S^kT^*M$) {\em le fibr\'{e} des $k$-vecteurs
sym\'{e}triques} (resp. {\em le fibr\'{e} des $k$-formes sym\'{e}triques}): pour
$k=0$ le fibr\'{e} $S^0TM:=M\times\real =: S^0T^*M$ et $S^kTM$ (resp. $S^kT^*M$)
est le
sous-fibr\'{e} de la $k^{\mathrm{\grave{e}me}}$ puissance tensorielle de $TM$
(resp. $T^*M$) qui soit invariant par l'action naturelle du groupe sym\'{e}trique
d'ordre $k$, $S_k$. On peut regarder les sections de classe $\Cinf$ de $S^kTM$
(resp. de $S^kT^*M$), c.-\`{a}-d.
{\em les champs de $k$-vecteurs sym\'{e}triques} (resp. {\em les champs de $k$-formes
sym\'{e}triques}), comme des application $k$-lin\'{e}aires sym\'{e}triques sur
l'espace de toutes les $1$-formes (resp. sur l'espace de tous les champs
de vecteurs). On rappelle la multiplication sym\'{e}trique $\vee$ de
$A_1\in \Ginf(S^kTM)$ et $A_2\in \Ginf(S^lTM)$:
\bea
   \lefteqn{(A_1\vee A_2)\big(\alpha_1,\ldots,\alpha_{k+l}\big)
    := } \nonumber \\
    & & \frac{1}{k!l!}\sum_{\sigma\in S_{k+l}}
              A_1\big(\alpha_{\sigma(1)},\ldots,\alpha_{\sigma(k)}\big)
             A_2\big(\alpha_{\sigma(k+1)},\ldots,\alpha_{\sigma(k+l)}\big)
\eea
o\`{u} $\alpha_1,\ldots,\alpha_{k+l}\in\Ginf(T^*M)$ et on a utilis\'{e} la convention
du livre de Greub, \cite{Gre78}. Une formule analogue duale
est vraie pour le produit sym\'{e}trique d'une $k$-forme sym\'{e}trique et d'une
$l$-forme sym\'{e}trique. Il est bien connu que $\vee$ est une multiplication
associative et commutative. Alors
$\big(\Ginf(STM):=\oplus_{k=0}^\infty \Ginf(S^kTM),\vee\big)$ (resp.
$\big(\Ginf(ST^*M):=\oplus_{k=0}^\infty \Ginf(S^kT^*M)$, $\vee\big)$)
est une alg\`{e}bre
commutative associative gradu\'{e}e. Soit $X$ un champ de vecteurs et $\gamma$ une
$l$-forme sym\'{e}trique. On rappelle le produit int\'{e}rieur
\[
    \big(i(X)\gamma\big)\big(X_1,\ldots,X_{k-1}\big)
         := \gamma\big(X,X_1,\ldots,X_{l-1}\big)
\]
(o\`{u} $X_1,\ldots,X_{l-1}\in\Ginf(TM)$) pour $l\geq 1$ et $i(X)\gamma:=0$
lorsque $l=0$. Le produit int\'{e}rieur est une d\'{e}rivation de degr\'{e} $-1$
de $\big(\Ginf(ST^*M),\vee\big)$. Pour $A:= X_1\vee\cdots\vee X_k$ on
d\'{e}finit
\beq
    i(A)\gamma = i(X_1\vee\cdots\vee X_k)\gamma
               := i(X_1)\cdots i(X_k)\gamma
\end{equation}
qui se prolonge de fa\c{c}on naturelle en un produit int\'{e}rieur de
$i:\Ginf(STM)\times \Ginf(ST^*M)\ra \Ginf(ST^*M)$. Pour $A_1,A_2\in\Ginf(STM)$
on a \'{e}videmment la propri\'{e}t\'{e} de module suivante:
\beq
    i(A_1\vee A_2) = i(A_1)i(A_2).
\end{equation}
En \'{e}changeant les r\^{o}les de $\Ginf(STM)$ et $\Ginf(ST^*M)$ on obtient par
dualit\'{e} un produit int\'{e}rieur de $i:\Ginf(ST^*M)\times \Ginf(STM)\ra
\Ginf(STM)$.

Soient $\tau:E\ra M$ et $\tau':E'\ra M$ deux fibr\'{e}s vectoriels sur $M$.
Alors l'espace des sections
$\Ginf(ST^*M\otimes E):= \oplus_{k=0}^\infty \Ginf(S^kT^*M\otimes E)$ est
un $\Ginf(ST^*M)$-module de fa\c{c}on naturelle, et on prolonge les produits
int\'{e}rieurs de $\Ginf(STM)$ \`{a} $\Ginf(ST^*M\otimes E)$. En outre, pour
les sections dans $\Ginf\big(STM\otimes ST^*M\otimes Hom(E,E')\big)$ il y a un
produit int\'{e}rieur d\'{e}fini par
\beq
  i(\gamma\otimes A\otimes \phi)\big(\gamma'\otimes A'\otimes \psi\big):=
     i(A)\big(\gamma'\big)\otimes i(\gamma)\big(A'\big)\otimes \phi(\psi)
\end{equation}
quels que soient $A,A'\in\Ginf(STM)$, $\phi\in Hom(E,E')$, $\gamma,\gamma'\in
\Ginf(ST^*M)$ et $\psi\in\Ginf(E)$.

Soit maintenant $\nabla^{TM}$ une connexion sans torsion dans le fibr\'{e} tangent,
et on utilise le m\^{e}me symbole $\nabla^{TM}$ pour sa prolongation naturelle
aux espaces de sections $\Ginf(STM)$ et $\Ginf(ST^*M)$. Soit $\nabla^E$
une connexion dans $E$, et soit $\nabla$ la connexion dans l'espace des
sections $\Ginf(E\otimes ST^*M)$ form\'{e}e de $\nabla^{TM}$ et $\nabla^E$. On
rappelle la {\em diff\'{e}rentielle sym\'{e}trique $\mathsf{D}$ par rapport \`{a}
$\nabla^{TM}$ (et $\nabla^E$)}: soient $\gamma\in\Ginf(E\otimes S^lT^*M)$ et
$X_1,\ldots,X_{l+1}\in\Ginf(TM)$, alors
\beq
   (\mathsf{D}\gamma)\big(X_1,\ldots,X_{l+1}\big)
     :=\frac{1}{l!}\sum_{\sigma\in S_l}
       (\nabla_{X_{\sigma(1)}}\gamma)
           \big(X_{\sigma(2)},\ldots,X_{\sigma(l+1)}\big).
\end{equation}
Pour le cas $E=M\times\real$ la diff\'{e}rentielle sym\'{e}trique $\mathsf{D}$ est une
d\'{e}rivation de degr\'{e} $+1$ de l'alg\`{e}bre $\big(\Ginf(ST^*M),\vee\big)$.
En g\'{e}n\'{e}ral, $\mathsf{D}$ est une d\'{e}rivation du $\Ginf(ST^*M)$-module
$\Ginf(E\otimes ST^*M)$, c.-\`{a}-d.:
\beq
   \mathsf{D}(\beta\vee\gamma)=(\mathsf{D}\beta)\vee\gamma + \beta\vee
                     (\mathsf{D}\gamma)
\end{equation}
quels que soient les \'{e}lements $\beta\in\Ginf(ST^*M)$ et
$\gamma\in \Ginf(E\otimes STM)$. On d\'{e}finit la {\em divergence covariante
par rapport \`{a} $\nabla$}, \'{e}crite $\mathsf{Div}$, par l'application
lin\'{e}aire $\Ginf(STM\otimes E)\ra \Ginf(STM\otimes E)$ de degr\'{e} $-1$
donn\'{e}e par
\beq \lbl{EqDefDiv}
   \mathsf{Div}A:=i(\mathbf{1})\mathsf{D}A
\end{equation}
o\`{u} $A\in \Ginf(STM\otimes E)$, et l'application identique $\mathbf{1}$
appartenant \`{a} l'espace $Hom(TM,TM)\otimes Hom(E,E)$ est consid\'{e}r\'{e}e comme \'{e}l\'{e}ment de l'espace des
sections
$\Ginf(S^1T^*M\otimes S^1TM \otimes Hom(E,E))$ de fa\c{c}on naturelle.
A l'aide d'une base locale $\partial_1,\ldots,\partial_m$ du fibr\'{e} tangent
et sa base duale $dx^1,\ldots,dx^m$ on peut \'{e}crire la divergence de
la forme
\[
   \mathsf{Div}A=\sum_{j=1}^m i(dx^{j})\big(\nabla_{\partial_j}A\big),
\]
et la diff\'{e}rentielle sym\'{e}trique par
\[
   \mathsf{D}A = \sum_{j=1}^m dx^{j}\vee \big(\nabla_{\partial_j}A\big).
\]

On rappelle qu'un {\em op\'{e}rateur diff\'{e}rentiel d'ordre $k$ de $E$ dans $E'$}
 est une
application lin\'{e}aire
$D:\Ginf(E)\ra\Ginf(E')$ qui satisfait la condition suivante:
dans toute carte $\big(U,(x^1,\ldots,x^m)\big)$ telle que la restriction des
fibr\'{e}s
\`{a} $U$, $E|_U$ et $E'|_U$, soient trivialisables on a pour tout $\psi\in
\Ginf(E)$ (o\`{u} $e_1,\ldots,e_{K}\in\Ginf(U,E|_U)$ est une base locale des
sections de $E$, $\psi|_U=\sum_{j=1}^K\psi^je_j$ et $e'_1,\ldots,e'_{K'}$
est une base locale des sections de $E'$)
\beq
    \big(D(\psi)\big)|_U=\sum_{a=0}^k\sum_{i_1,\ldots,i_a=1}^m
    \sum_{j=1}^K\sum_{j'=1}^{K'}D^{a;i_1\ldots i_a;j'}_j
    \frac{\partial^a \psi^j|_U}{\partial x^{i_1}\cdots \partial
    x^{i_a}}e'_{j'}
\end{equation}
o\`{u} les $D^{a;i_1\ldots i_a;j'}_j$ sont des fonctions de classe $\Cinf$ sur
$U$ \`{a} valeurs r\'{e}elles ou complexes (d\'{e}pendant de la nature des fibr\'{e}s
vectoriels $E$ et $E'$).
L'espace de tous les op\'{e}rateurs diff\'{e}rentiels de
$E$ dans $E'$ est not\'{e} par $\Dop\big(\Ginf(E),\Ginf(E')\big)$.

\bdefi \lbl{DPreOrdStandard}
  Soient $M,E,E',\nabla,\nabla^E$ les structures d\'{e}finies ci-dessus. La
  {\em pr\'{e}\-scrip\-tion d'ordre standard (par rapport \`{a} $\nabla$ et $\nabla^E$)},
  \'{e}crite $\rho_s$, est l'application
  $\real$-lin\'{e}aire $\rho_s[\nabla]=\rho_s: \Ginf\big(STM\otimes Hom(E,E')\big)\ra
  \Dop\big(\Ginf(E),\Ginf(E')\big)$ donn\'{e}e par
  \beq \lbl{EqPreOrdStandard}
    \big(\rho_s(A)\big)(\psi):=
            i(A)\big(\mathsf{D}^k\psi\big)
  \end{equation}
  o\`{u} $A\in \Ginf\big(S^kTM\otimes Hom(E,E)\big)$ et $\psi\in\Ginf(E)$.
\edefi
Si l'on regarde les diff\'{e}rentielles sym\'{e}triques it\'{e}r\'{e}es en coordonn\'{e}es on
voit ais\'{e}ment que $\rho_s$ est bien d\'{e}finie et bijective.
De plus, il est bien connu que pour chaque op\'{e}rateur diff\'{e}rentiel $D$
non nul
la composante non nulle du plus haut degr\'{e} de
$\rho_s^{-1}(D)$ ne d\'{e}pend pas de la connexion choisie: ceci est appel\'{e}
le {\em symbole principal de $D$}.

\bdefi \lbl{DPreOr}
 Soient $E$ et $E'$ deux fibr\'{e}s vectoriels sur une vari\'{e}t\'{e} dif\-f\'{e}\-ren\-tiable $M$
 de dimension $m$. Une {\em pr\'{e}scription d'ordre} est une bijection
 $\real$-lin\'{e}aire
 \[
    \rho:\Ginf\big(STM\otimes Hom(E,E')\big)\ra
             \Dop\big(\Ginf(E),\Ginf(E')\big)
 \]
 telle que pour chaque op\'{e}rateur diff\'{e}rentiel $D$ non nul la composante
 non nulle du plus haut degr\'{e} de $\rho^{-1}(D)$ co\"{\i}ncide avec le
 symbole principal de $D$.
\edefi
A part la prescription d'ordre standard il y a d'autres pr\'{e}scriptions d'ordre
motiv\'{e}es par la m\'{e}canique quantique, voir par exemple \cite{BNW98} et
\cite{BNPW98}: \`{a} l'aide de la
divergence covariante (\ref{EqDefDiv}) on peut former {\em l'op\'{e}rateur de
Neumaier}
\beq
     N[\nabla]:=e^{\frac{1}{2}\mathsf{Div}}
\end{equation}
(qui est bien d\'{e}finie sur les \'{e}l\'{e}ments de $\Ginf\big(ST^M\otimes
Hom(E,E')\big)$) et d\'{e}finir la {\em pr\'{e}scription d'ordre de type Weyl}
\beq
     \rho_w(A):=\rho_s[\nabla]\big(N[\nabla](A)\big).
\end{equation}

Soit $G$ un groupe de Lie (o\`{u} $\mathfrak{g}$ d\'{e}signe son alg\`{e}bre de Lie)
et soient $\Phi:G\times M\ra M$,
$\Phi^E:G\times E\ra E$ et $\Phi^{E'}:G\times E'\ra E'$ des actions \`{a}
gauche de $G$ telles que $G$ agit sur $E$ et $E'$ par des morphsimes
de fibr\'{e}s vectoriels tels que
\[
  \tau\big((\Phi^E(g,e)\big)=\Phi\big(g,\tau(e)\big)
   \mathrm{~~~et~~~}
  \tau'\big((\Phi^{E'}(g,e')\big)=\Phi\big(g,\tau'(e')\big)
\]
quels que soient $g\in G$, $e\in E$ et $e'\in E'$.
En \'{e}crivant $\Phi_g,\Phi^E_g$ et $\Phi^{E'}_g$ pour les applications
$\Phi(g,~),\Phi^E(g,~)$ et $\Phi^{E'}(g,~)$ pour $g\in G$ on peut
retirer les sections $\varphi\in\Ginf(E)$ et
$A\in \Ginf\big(M,S^kTM\otimes Hom(E,E')\big)$ par
\[
   (\Phi_g^*\varphi)_x:=\Phi^E_{g^{-1}}\varphi_{\Phi_g(x)}
   ~~~\mathrm{et}~~~
   (\Phi_g^*A)_x(e) :=\big(\Phi^{E'}_{g^{-1}}\otimes S^kT\Phi_{g^{-1}}\big)
      A_{\Phi_g(x)}\big(\Phi^E_g(e)\big).
\]
\bdefi
 Soient $G,\Phi,\Phi^E$ et $\Phi^{E'}$ comme d\'{e}finies ci-dessus.
 Une prescription d'ordre $\rho$ est dite {\em $G$-\'{e}quivariante} lorsque
 \beq \lbl{EqDefGEquivariPreOr}
    \Phi_g^*\big(\rho(A)\varphi\big)=\rho(\Phi^*_gA)\Phi_g^*\varphi
 \end{equation}
 quels que soient $g\in G$, $\varphi\in \Ginf(E)$ et
 $A\in \Ginf\big(M,S^kTM\otimes Hom(E,E')\big)$.
\edefi
Bien s\^{u}r, il y a une version infinit\'{e}simale de cette $G$-\'{e}quivariance
que l'on obtient en mettant $g=\exp(t\xi)$ avec $t\in\real$ et
$\xi\in\mathfrak{g}$  et en diff\'{e}rentiant l'\'{e}quation
(\ref{EqDefGEquivariPreOr}) par rapport \`{a} $t$ en $t=0$:

Une {\em action \`{a} gauche $\phi$ de l'alg\`{e}bre de Lie $\mathfrak{g}$ sur $M$} est
 un
antihomomorphisme d'alg\`{e}bres de Lie $\mathfrak{g}\ra\Ginf(TM)$, \'{e}crit
$\xi\mapsto \xi_M$, o\`{u} $[\xi_M,\eta_M]=-[\xi,\eta]_M$ quels que soient
$\xi,\eta\in\mathfrak{g}$. Une {\em action \`{a} gauche de $\mathfrak{g}$ sur
le fibr\'{e} vectoriel $E$
(compatible avec l'action \`{a} gauche $\phi$ de $\mathfrak{g}$ sur $M$)} est un
antihomomorphisme $\phi^E$ d'alg\`{e}bres de Lie
$\mathfrak{g}\ra\Dop\big(\Ginf(E),\Ginf(E')\big)$,
\'{e}crit $\xi\mapsto \phi^E_\xi$ (c.-\`{a}-d.
$\phi^E_{[\xi,\eta]}=-\phi^E_\xi\phi^E_\eta+\phi^E_\eta\phi^E_\xi$), tel que
de plus $\phi^E_\xi(\tau^*f\varphi)=\tau^*\big(\xi_M(f)\big)\phi^E_\xi(\varphi)$
quels que soient $\xi\in\mathfrak{g}$, $f\in\Cinf(M,\real)$,
$\varphi\in\Ginf(E)$. Lorsque $\phi^E$ et $\phi^{E'}$ sont deux actions \`{a}
gauche de $\mathfrak{g}$ sur les fibr\'{e}s vectoriels $E$ et $E'$, respectivement,
compatibles avec l'action \`{a} gauche $\phi$, on a une action \`{a} gauche $\phi^{(k)}$
de $\mathfrak{g}$ sur le fibr\'{e} vectoriel $S^kTM\otimes Hom(E,E')$ d\'{e}finie
par
\beas
  \lefteqn{\big(\phi^{(k)}_\xi(X_1\vee\cdots\vee X_k\otimes
  B)\big)(\varphi)} \hspace{3cm}\\
   & := &\sum_{l=1}^k X_1\vee\cdots\vee[\xi_M,X_l]\vee\cdots\vee X_k\otimes
   B(\varphi)\\
   &    &
     +X_1\vee\cdots\vee X_k\otimes \left(\phi^{E'}_\xi\big(B(\varphi)\big)
     - B\big(\phi_\xi^{E}(\varphi)\big)\right)
\eeas
quels que soient $\xi\in\mathfrak{g}$, $X_1,\ldots, X_k\in\Ginf(TM)$,
$B\in Hom(E,E')$ et $\varphi\in\Ginf(E)$.
\bdefi
 Soient $\mathfrak{g},\phi,\phi^E$,$\phi^{E'}$ et $\phi^{(k)}$
 comme d\'{e}finies ci-dessus.
 Une prescription d'ordre $\rho$ est dite {\em $\mathfrak{g}$-\'{e}quivariante}
 lorsque
 \beq \lbl{EqDefgEquivariPreOr}
    \phi^{E'}_\xi\big(\rho(A)\varphi\big)=
    \rho\big(\phi^{(k)}_\xi (A)\big)\varphi
    +\rho(A)\big(\phi_\xi^{E}\varphi\big)
 \end{equation}
 quels que soient $\xi\in\mathfrak{g}$, $\varphi\in \Ginf(E)$ et
 $A\in \Ginf\big(S^kTM\otimes Hom(E,E')\big)$.
\edefi

\section{Prescriptions d'ordre naturelles}

Dans l'appendice on a donn\'{e} un rappel de la th\'{e}orie des fibr\'{e}s principaux
et des fibr\'{e}s et op\'{e}rateurs naturels d'apr\`{e}s les ouvrages
\cite{KN63},\cite{KN69} et surtout \cite{KMS93} pour fixer la notation
et quelques \'{e}quations importantes.

On peut maintenant donner la d\'{e}finition d'une prescription d'ordre
naturelle qui sera une l\'{e}g\`{e}re g\'{e}n\'{e}ralisation d'une d\'{e}finition donn\'{e}e
par P.B.A. Lecomte:
\bdefi[P.B.A. Lecomte]\lbl{DefLecNat}
 Soient $F$ et $F'$ deux foncteurs de fibr\'{e}s vectoriels d'ordre $r$ et
 $r'$, respectivement, et soit $q:=\max(r,r')$. Pour une varit\'{e}t\'{e} diff\'{e}rentiable
 $M$ de dimension $M$  et une $1$-forme de
 connexion dans $\Ginf(Q_\tau P^qM)$ on a toujours la connexion induite
 $\nabla^{TM}$
 sans torsion dans le fibr\'{e} tangent et les connexions induites $\nabla^{FM}$
 et $\nabla^{F'M}$ dans les fibr\'{e}s vectoriels associ\'{e}s $FM$ et $F'M$.
 Pour tout entier positif $k$ on consid\`{e}re les foncteurs de fibr\'{e}s d'ordre
 $q+1$ et $\max(k+r,r')$
 \beas
   L^k_C  & := & Q_\tau P^q\times \big(S^kT\otimes
                                    Hom(F,F')\big)~~\mathrm{et} \\
   L^k_Q  & := & Hom(J^k\circ F,F').
 \eeas
 Une {\em prescription d'ordre naturelle} est une collection d'op\'{e}rateurs
 naturels d'ordre fini $\rho=\big(\rho^k:L^k_C\leadsto L^k_Q\big)_{k\in\nat}$ telle
 que pour
 toute vari\'{e}t\'{e} $M$ de dimension $m$ et pour toute $1$-forme de connexion
 $\omega$ dans
 $P^qM$  l'application $\rho_M[\nabla]$
 de $\oplus_{k=0}^\infty\Ginf\big(S^kTM\otimes Hom(FM,F'M)\big)$ dans
 $\Dop(FM,F'M)$ suivante soit une pr\'{e}scription d'ordre: soient
 $(f_k)_{0\leq k\leq r}$ une famille de section o\`{u}
 $f_k\in \Ginf\big(S^kTM\otimes Hom(FM,F'M)\big)$ et $f:=f_0+\cdots+ f_r$,
 alors
 \beq
       \rho_M[\nabla](f):=\rho^0_M(\omega,f_0)+\cdots
       +\rho^k_M(\omega,f_k).
 \end{equation}
 o\`{u} $\nabla$ d\'{e}signe la connexion induite par $\nabla^{TM},\nabla^{FM}$ et
 $\nabla^{F'M}$ dans les fibr\'{e}s $S^kTM\otimes Hom(FM,F'M)$.
\edefi
Grosso modo, une prescription d'ordre naturelle est une prescription
d'ordre d\'{e}pendant d'une connexion qui se transforme `bien' quand on
applique des diff\'{e}omorphismes locaux. On va expliciter ce changement: soit
$\Phi:M\ra N$ une immersion entre deux vari\'{e}t\'{e}s diff\'{e}rentiables de
dimension $m$, soit $\nabla'$ la connexion dans les fibr\'{e}s
$S^kTN\otimes Hom(FN,F'N)$ r\'{e}sultant d'une $1$-forme de connexion
$\omega'$ dans le fibr\'{e} $P^qN$ comme ci-dessus, soit $\varphi'$ une section dans
$\Ginf(FN)$ et soit
$f'$ une section dans $\oplus_{k=0}^\infty\Ginf\big(S^kTN\otimes
Hom(FN,GN)\big)$. Alors la `naturalit\'{e}' de $(\rho^k)_{k\in\nat}$ veut dire que
 \beq \lbl{EqDefNaturalitePreOrd}
   \Phi^*\left(\big(\rho_N[\nabla'](f')\big)(\varphi')\right)
       =\big(\rho_M[\Phi^*\nabla'](\Phi^*f')\big)(\Phi^*\varphi').
 \end{equation}
o\`{u} $\Phi^*\nabla'$ d\'{e}signe la connexion dans les fibr\'{e}s $S^kTM\otimes Hom(FM,F'M)$
induite par la $1$-forme de connexion $(P^q\Phi)^*\omega'$ retir\'{e}e.

\noindent L'exemple le plus important pour une telle prescription d'ordre est le
suivant:
\bcor
   Soit $M$ une vari\'{e}t\'{e} diff\'{e}rentiable de dimension $m$. Soient les fibr\'{e}s vectoriels
   $E$ et $E'$ sur $M$ les objets $FM$ et $F'M$ des foncteurs de fibr\'{e}s
  vectoriels d'ordre $r$ et $r'$, $F$ et $F'$, respectivement, soit
  $\omega$ une connexion dans le fibr\'{e} $P^qM$ (o\`{u} $q:=\max(r,r')$)
  induisant les connexions
  $\nabla^{TM}$ sans torsion dans $TM$, $\nabla^{E}$ dans $E$ et
  $\nabla^{E'}$ dans $E'$ et soit $\nabla$ la connexion induite dans les fibr\'{e}s
  $S^kTM\otimes Hom(E,E')$.\\
  Alors la collection de toutes les prescriptions d'ordre standard
  $\rho_s[\nabla]$ (voir la d\'{e}finition \ref{DPreOrdStandard}) d\'{e}termine une
  prescription d'ordre naturelle.
\ecor

\section{Fibr\'{e} (principal) de densit\'{e}s}

Pour la notation de la th\'{e}orie des fibr\'{e}s principaux on peut consulter
l'appendice.

Soit $a\in\real$. On rappelle la d\'{e}finition
du {\em fibr\'{e} des $a$-densit\'{e}s sur $M$}, not\'{e} $\tau^a:|\Lambda^m T^*M|^a\ra M$:
ce fibr\'{e}
est le fibr\'{e} vectoriel associ\'{e} \`{a} $P^1M$ \`{a} fibre type $\real$ o\`{u} l'action
du groupe structurel $GL(m,\real)$ de $P^1M$ est donn\'{e}e par
\beq \lbl{EqActGl}
     GL(m,\real)\times \real \ra \real: (g,\lambda)\mapsto |\det g|^{-a}\lambda.
\end{equation}
Alors
\beq
     |\Lambda^m T^*M|^a := P^1M\times_{GL(m,\real)} \real.
\end{equation}
On va noter $q^a:P^1M\times \real \ra |\Lambda^m T^*M|^a:
(p,\lambda)\mapsto q^a(p,\lambda)$
la projection canonique. Les sections de classe $\Cinf$ de $|\Lambda^m
T^*M|^a$ sont dites les {\em $a$-densit\'{e}s de $M$}.
Dans le cas $a=0$ l'application canonique
\[
   |\Lambda^m T^*M|^0\ra M\times \real : q^0(p,\lambda)\mapsto
   (\pi^1_0(p),\lambda)
\]
est un isomorphisme de fibr\'{e}s vectoriels bien d\'{e}fini, alors les
$0$-densit\'{e}s s'identifient de fa\c{c}on naturelle aux fonctions de classe
$\Cinf$ \`{a} valeurs r\'{e}elles.
Puisque toute vari\'{e}t\'{e} peut \^{e}tre munie d'une m\'{e}trique riemannienne $\mathsf{g}$
(\`{a}
laquelle correspond une fonction $GL(m,\real)$-\'{e}quivariante
$f_{\mathsf{g}}:P^1M\ra \real^{n*}\times
\real^{n*}$) il s'ensuit que tout fibr\'{e} $|\Lambda^m T^*M|^a$ admet une
section $|\mathsf{g}|^{a/2}$ non nulle dont sa fonction \'{e}quivariante est d\'{e}finie par
\beq \lbl{EqDefDetmathsfgademi}
   f_{|\mathsf{g}|^{a/2}}(e_1,\ldots,e_m)
      :=|\det \mathsf{g}(e_i,e_j)|^{\frac{a}{2}}
\end{equation}
o\`{u} $(e_1,\ldots,e_m)=:p\in P^1 M$. Par cons\'{e}quent, tous les
fibr\'{e}s $|\Lambda^m T^*M|^a$ sont trivialisables: soit $\sigma$ une $a$-densit\'{e}
non nulle, par exemple $|\mathsf{g}|^{a/2}$. On consid\`{e}re
l'application
\beq \lbl{EqDefKsigma}
   K_\sigma:M\times \real \ra |\Lambda^m T^*M|^a:
       (x,\lambda)\mapsto \sigma(x)\lambda
\end{equation}
qui donne \'{e}videmment un isomorphisme de fibr\'{e}s vectoriels.
Pour $a\neq 0$ cette trivialisation d\'{e}pend en g\'{e}n\'{e}ral de la section
$\sigma$.

En tant que fibr\'{e} associ\'{e} du fibr\'{e} des rep\`{e}res lin\'{e}aires, l'association
$M\mapsto |\Lambda^m T^*M|^a$ et $|\Lambda^m T^*\Phi|^a\big(q^a(p,\lambda)\big)
:=q^a\big(P^1\Phi(p),\lambda\big)$ (pour une immersion $\Phi$ entre deux
vari\'{e}t\'{e}s $M$ et $N$ de dimension $m$) est un foncteur de fibr\'{e}s
vectoriels.

La multiplication des nombres r\'{e}els induit une multiplication fibre-par-fibre
sur les fibr\'{e}s de densit\'{e}s: $q^a(p,\lambda)q^b(p,\lambda')
:=q^{a+b}(p,\lambda\lambda')$ quels que
soient $a,b,\lambda,\lambda'\in real$, $p\in P^1M$, alors
\beq
  |\Lambda^m T^*M|^a \times |\Lambda^m T^*M|^b \ra |\Lambda^m
  T^*M|^{a+b}~~\forall a,b\in\real.
\end{equation}
Par cons\'{e}quent, la multiplication avec un \'{e}l\'{e}ment
$y\in|\Lambda^m T^*M|^a$ avec $\tau^a(y)=x\in M$ d\'{e}finit un homomorphisme
de la fibre sur $x$ de $|\Lambda^m T^*M|^b$ dans la fibre sur $x$ de
$|\Lambda^m T^*M|^{a+b}$. D'un autre cot\'{e}, puisque le fibr\'{e}
de tous les homomorphismes,
$Hom(|\Lambda^m T^*M|^a,|\Lambda^m T^*M|^b)$ est un fibr\'{e} associ\'{e} de $P^1M$
\`{a} fibre type $Hom(\real,\real)\cong\real$ o\`{u} la repr\'{e}sentation de $GL(m,\real)$
est donn\'{e}e par $(g,\lambda)\mapsto |\det g|^{a}|\det g|^{-b}\lambda$, alors
\beq
  Hom(|\Lambda^m T^*M|^a,|\Lambda^m T^*M|^b)\cong |\Lambda^m T^*M|^{b-a}~~
         \forall a,b\in\real.
\end{equation}

On d\'{e}finit la partie ouverte $\tilde{M}^a$ de $|\Lambda^mT^*M|^a$ suivante:
\beq \lbl{EqDefDensFibr}
   \tilde{M}^a := P^1M\times_{GL(m,\real)}\real^+
\end{equation}
o\`{u} $G:=\real^+$ est l'ensemble de tous les nombres r\'{e}els strictement
positifs et l'action de $GL(m,\real)$ est comme dans eqn(\ref{EqActGl}).
On note la restriction de $q^a$ \`{a} $P^1M\times \real^+$ et la restriction
de $\tau^a$ \`{a} $\tilde{M}^a$ par les m\^{e}mes symboles $q^a$ et $\tau^a$,
respectivement.
\bprop \lbl{PPropriTildeMa}
  Pour tout r\'{e}el $a$ il vient:
  \ben
   \item $\tilde{M}^a$ munie de l'action droite du groupe multiplicatif
    $\real^+$ donn\'{e}e par
    \[
        \tilde{M}^a \times \real^+ \ra \tilde{M}^a:
            \big(q^a(p,\lambda),s\big)\mapsto \big(q^a(p,\lambda)\big)s
             := q^a(p,\lambda s)
    \]
    et de la projection $\tau^a:\tilde{M}^a\ra M$ est un fibr\'{e} principal sur
    $M$ \`{a} groupe structurel $\real^+$. De plus, le champ fondamental $1^*$ est
    \'{e}gal \`{a} la restriction du champ d'Euler $\mathsf{E}$ de
    $|\Lambda^mT^*M|^a$ (dont le flot $G_t$ est donn\'{e} par $y\mapsto
    ye^t$).
   \item L'application
    \[
       \Xi^a:P^1M\ra \tilde{M}^a: p \mapsto q^a(p,1)
    \]
    est un morphisme de fibr\'{e}s principaux sur $M$ induisant
    l'homomorphisme $\xi^a:GL(m,\real)\ra \real^+:g\mapsto |\det g|^{-a}$ et
    l'application identique sur $M$.
    $\Xi^a$ est une submersion surjective si et seulement si $a\neq 0$.
   \item Soit $b\in\real$ et $a\neq 0$. Alors le fibr\'{e} $|\Lambda^mT^*M|^b$
     s'obtient
     comme fibr\'{e} associ\'{e} au fibr\'{e} principal $\tilde{M}^a$ \`{a} fibre type $\real$
     o\`{u} l'action du groupe $\real^+\times \real\ra \real$ est donn\'{e}e par
     \[
          (\lambda,s)\mapsto \lambda^{\frac{b}{a}}s.
     \]
     On va noter $\tilde{q}^a_b:\tilde{M}^a\times \real \ra
     |\Lambda^mT^*M|^b$ pour la projection.
     En particulier, \`{a} chaque section
     $\varphi \in \Ginf\big(|\Lambda^mT^*M|^b\big)$
     on associe une fonction $\tilde{\varphi}:\tilde{M}^a\ra\real$ qui est
     $(-\frac{b}{a})$-\'{e}quivariante dans le sens suivant:
     \[
        \tilde{\varphi}(y\lambda)=\lambda^{-\frac{b}{a}}\tilde{\varphi}(y)
          ~~\forall y\in
             \tilde{M}^a,~~\forall \lambda\in\real^+
     \]
     En notant l'inverse de l'application $t\mapsto \tilde{q}^a_b(y,t)=s$ de
     $\real$
     dans la fibre $|\Lambda^mT^*M|^b_x$ par $t=y^{-1}s$ on a
     $\tilde{\varphi}(y)=y^{-1}\varphi\big(\tau^a(y)\big)$.
  \een
\eprop
\bbew
 1. $\tilde{M}^a$ est un ouvert du fibr\'{e} vectoriel $|\Lambda^m T^*M|^a$,
 et en utilisant la d\'{e}finition de $\tilde{M}^a$ on voit qu'une trivialisation
 locale de $|\Lambda^m T^*M|^a$ induit une trivialisation locale de
 $\tilde{M}^a$ dont la fibre type est \'{e}gale $\real^+$. En outre, la
 multiplication des fibres par les nombres r\'{e}els strictement positifs dans le fibr\'{e}
 vectoriel $|\Lambda^m T^*M|^a$ (qui peut toujours \^{e}tre exprim\'{e}e en termes de la
 projection
 $q^a$) induit une action droite du groupe
 $\real^+$ sur $\tilde{M}^a$ qui est libre parce que $\tilde{M}^a$ ne
 rencontre pas la section nulle de $|\Lambda^m T^*M|^a$. L'\'{e}quation
 $\mathsf{E}=1^*$ est \'{e}vidente. \\
 2. Puisque $\tau^a\big(q^a(p,r)\big)$ est d\'{e}finie par $\pi^1_0(p)$ il
 est clair que $\Xi^a$ pr\'{e}serve les fibres induisant l'application
 identique sur $M$. De plus, pour tout $g\in GL(m,\real)$:
 \[
     \Xi^a(pg) = q^a(pg,1) = q^a(p,|\det g|^{-a}1)= q^a(p,1)|\det g|^{-a}
                            = \Xi^a(p)|\det g|^{-a}.
 \]
 ce qui prouve toutes les propri\'{e}t\'{e}s d'homomorphismes. Soit $\lambda >0$.
 Si $a\neq 0$ il
 existe un $s>0$ tel que $\lambda=s^{-ma}=|\det(s\mathbf{1})|^{-a}$ o\`{u}
 $\mathbf{1}$
 d\'{e}signe la matrice identique dans $GL(m,\real)$, d'o\`{u}
 $q^a(p,\lambda)=q^a(pg,1)=\Phi^a(pg)$ avec $g:=s\mathbf{1}$.
 Alors $\Xi^a$ est
 surjective. Puisque $\Phi^a$ pr\'{e}serve les fibres, il suffit de prouver la
 surjectivit\'{e} de son application tangente sur les vecteurs verticaux: soit
 $b\in\mathfrak{gl}(m,\real)$ et $b^*$ le champ fondamental sur $P^1M$,
 $b^*_p:=\frac{d}{dt}\big(p\exp(tb)\big)|_{t=0}$. Alors
 \beas
   T_p\Xi^a~b^*_p & = &
   \frac{d}{dt}\big(\Xi^a(p)|\det(e^{tb})|^{-a}\big)|_{t=0}
    = \frac{d}{dt}\big(\Xi^a(p)e^{-ta\mathrm{tr}(b)}\big)|_{t=0}
    = -a\mathrm{tr}(b)\mathsf{E}_{\Xi^a(p)},
 \eeas
 et le membre droit ne s'annule pas si $b$ est un multiple non nul de la
 matrice identique. Donc $\Xi^a$ est une submersion si $a\neq 0$.
 Par contre, si $a=0$, alors
 $\Xi^0(pg)=\Xi^0(p)$ quel que soit $g\in GL(m,\real)$ et $\Xi^0$ n'est
 pas surjective.\\
 3. Il existe une unique fonction $f_\varphi:P^1M\ra\real$ de classe $\Cinf$
 telle que $f_\varphi(pg)=|\det g|^b f_\varphi(p)$. Soit $y:=q^a(p,s)\in \tilde{M}^a$.
 Gr\^{a}ce \`{a} l'\'{e}quivariance de $f_\varphi$ il vient que la fonction
 $\tilde{\varphi}:\tilde{M}^a\ra\real$ suivante est bien d\'{e}finie:
 \beq
     \tilde{\varphi}(y)=\tilde{\varphi}\big(q^a(p,s)\big):=
                    f_\varphi(p)s^{-\frac{b}{a}}
 \end{equation}
 et a l'\'{e}quivariance \'{e}nonc\'{e}e. Le fait que $|\Lambda^m T^*M|^b$ s'obtient
 comme fibr\'{e} associ\'{e} de $\tilde{M}^a$ se d\'{e}montre de mani\`{e}re analogue.
\ebew
\bdefi
 Le fibr\'{e} principal $(\tilde{M}^a,\tau^a,M,\real^+)$ est dit le
 {\em fibr\'{e} principal des $a$-densit\'{e}s sur $M$}.
\edefi
Si $\sigma$ est une section non nulle du fibr\'{e} $|\Lambda^mT^*M|^a$, alors
ou bien $\sigma$ ou bien $-\sigma$ est une section de $\tilde{M}^a$.
Dans le premier cas (comme par exemple pour $|\mathsf{g}|^{a/2}$)
la restriction \`{a} $\tilde{M}^a$ de l'application $K_\sigma$ (voir
(\ref{EqDefKsigma})) donne un isomorphisme du fibr\'{e} principal
trivial $M\times \real^+$ avec $\tilde{M}^a$.

Soit $\Phi:M\ra M'$ une immersion o\`{u} $M,M'$ sont des vari\'{e}t\'{e}s
diff\'{e}rentiables de dimension $m$. Pour tout $a\in\real$ on d\'{e}finit une application
$\tilde{\Phi}^a:\tilde{M}^a\ra \widetilde{M'}^a$ de classe $\Cinf$ par
\beq \lbl{EqDefFuncDens}
  \tilde{\Phi}^a\big(q^a(p,\lambda)\big):= {q'}^a\big(P^1\Phi(p),\lambda\big)
\end{equation}
quels que soient $p\in P^1M,\lambda\in\real^+$.
On a le th\'{e}or\`{e}me suivant:
\bsat
 L'association $\tilde{F}^a:\mathcal{M}f_m\ra \mathcal{FM}$ d\'{e}finie par
 $M\mapsto \tilde{M}^a$ (voir eqn \ref{EqDefDensFibr}) et
 $\Phi\mapsto\tilde{\Phi}^a$ (voir eqn \ref{EqDefFuncDens}) est
 un foncteur de fibr\'{e}s. De plus, tout $\tilde{\Phi}^a$ est une immersion,
 alors $\tilde{F}^a$ d\'{e}finit un foncteur de $\mathcal{M}f_m$ dans
 $\mathcal{M}f_{m+1}$.
\esat
\bbew
 Puisque $P^1\Phi:P^1M\ra P^1M'$ est un morphisme de fibr\'{e}s principaux il
 vient que $P^1\Phi(pg)=P^1\Phi(p)g$ quel que soit $g\in GL(m,\real)$, alors
 l'application $\tilde{\Phi}^a$ est bien d\'{e}finie et un morphisme de fibr\'{e}s
 principaux. Le dernier \'{e}nonc\'{e} est vraie pour tous les foncteurs de
 fibr\'{e}s.
\ebew

\noindent Pour calculer $\tilde{\Phi}^a$ on peut utiliser une section non nulle
$\sigma$ de $\tilde{M}^a$:
\blem
 Avec les notations mentionn\'{e}es ci-dessus, soit $f_\Phi$ la fonction de
 classe $\Cinf$ \`{a} valeurs r\'{e}elles strictement positives sur $M$ d\'{e}finie
 par $f_\Phi\sigma:=\Phi^*\sigma$. Alors on a la formule suivante:
 \beq
    \tilde{\Phi}^a\big(K_\sigma(x,\lambda)\big)
      = K_\sigma\big(\Phi(x),{\textstyle \frac{\lambda}{f_\Phi} }\big)
 \end{equation}
\elem
\bbew
 Puisque $\tilde{\Phi}^a\big(\sigma(x)\big)$ est un \'{e}l\'{e}ment de la fibre
 de $\tilde{M}^a$ sur $\Phi(x)$, alors il existe une fonction $g_\Phi$ de classe
 $\Cinf$ \`{a} valeurs r\'{e}elles strictement positives sur $M$ telle que
 $\tilde{\Phi}^a\big(\sigma(x)\big)=\tilde{\Phi}_x^a\big(\sigma(x)\big)
 =g_\Phi \sigma\big(\Phi(x)\big)$.
 La d\'{e}fintion de $\Phi^*\sigma$ montre que $g_\Phi=f_\Phi^{-1}$.
\ebew

Soit $\omega$ une $1$-forme de connexion sur $P^1M$ et $\nabla$ la connexion
dans le fibr\'{e} tangent associ\'{e}e. Puisque
l'homomorphisme $\Xi^a:P^1M\ra \tilde{M}^a$ induit l'application identique
sur la base $M$, il existe une unique $1$-forme de connexion $\omega^a$
\`{a} valeurs r\'{e}elles sur $\tilde{M^a}$ d\'{e}finie par
\beq \lbl{EqConnPTilde}
   \omega^a_{\Xi(p)}\big(T_p\Xi\hspace{1mm}v\big)
     := -a\mathrm{tr}\big(\omega_p(v)\big)~~\forall p\in P^1M, v\in T_pP
\end{equation}
car $T_e\xi=T_e|\det|^{-a}=-a\mathrm{tr}$.
La vari\'{e}t\'{e}
 $\tilde{M}^a$ est un ouvert de l'espace total
 $|\Lambda^m T^*M|^a$, et le sous-fibr\'{e} horizontal $H\tilde{M}^a$ de
 $T\tilde{M}^a$ d\'{e}finie par la $1$-forme de connexion $\omega^a$ co\"{\i}ncide
 avec la restriction du sous-fibr\'{e} horizontal $H|\Lambda^m T^*M|^a$ de
 $T|\Lambda^m T^*M|^a$ induit
 par celui qui est d\'{e}fini par la $1$-forme de connexion $\omega$ sur $P^1M$.
 Donc les restrictions \`{a} $\tilde{M}^a$ des rel\`{e}vements horizontaux
 $X^{\mathbf{h}}$ et $Y^{\mathbf{h}}$ des champs de vecteurs $X$ et $Y$ sur $M$
 \`{a} $|\Lambda^m T^*M|^a$ co\"{\i}ncident avec les rel\`{e}vements horizontaux
 $X^h$ et $Y^h$ \`{a} $\tilde{M}^a$
Soit $\mathrm{tr}R$ la trace du tenseur de courbure de $\nabla$, c.-\`{a}-d.
$(\mathrm{tr}R)_x(v,w):=\mathrm{trace}\big(z\mapsto R_x(v,w)z\big)$ quels
que soient $x\in M$, $v,w,z\in T_xM$.
\blem
 Soient $X$, $Y$ deux champs de vecteurs sur $M$ et $X^h$, $Y^h$ leurs
 rel\`{e}vements horizontaux \`{a} $\tilde{M}^a$. Alors
 \beq \lbl{EqCrLieHorHorMtilde}
   [X^h,Y^h]=[X,Y]^h +a (\mathrm{tr}R)(X,Y)\mathsf{E}.
 \end{equation}
\elem
\bbew
 D'apr\`{e}s l'\'{e}quation (\ref{EqCrLieHorPrincipal}) et le fait que le champ d'Euler
 est le champ fondamental $1^*$ il faut montrer que
 $\Omega^a(X^h,Y^h)=-a(\mathrm{tr}R)(X,Y)$. Soient $X^H$ et $Y^H$ les
 rel\`{e}vements horizontaux de $X,Y$ \`{a} $P^1M$. Il vient que $X^H$ et $X^h$
 sont $\Xi^a$-li\'{e}s, c.-\`{a}-d. $T_p\Xi^a\hspace{1mm}X^H_p=X^h_{\Xi(p)}$
 $\forall p\in P^1M$, et l'on a
 \beas
    \Omega^a(X^h,Y^h)_{\Xi^a(p)} &= & d\omega^a_{\Xi^a(p)}(X^h,Y^h)
    +[\omega^a_{\Xi^a(p)}(X^h),\omega^a_{\Xi^a(p)}(Y^h)]\circ\Xi^a \\
                      & = & d({\Xi^a}^*\omega^a)_p(X^H,Y^H) + 0 \\
                      & = &
                      -a\mathrm{tr}\big((d\omega)_p\big(X^H,Y^H\big)
                                +[\omega_p(X^H),\omega_p(Y^H)]\big) \\
                      & = & -a\mathrm{tr}\big(\Omega_p(X^H,Y^H)\big) \\
                      & = & -a\mathrm{trace}\big(
                             Z\mapsto
                              \Omega_p(X^H,Y^H)\big({\theta_M}_p(Z^H)\big)\big)\\
                      & = & -a(\mathrm{tr}R)(X,Y)_{\Xi^a(p)}.
 \eeas
\ebew

\section{Connexions projectivement \'{e}quivalentes}

\noindent Soit $M$ une vari\'{e}t\'{e} diff\'{e}rentiable de dimension $m$.

\bdefi
  Deux connexions $\nabla$ et $\nabla'$ dans le fibr\'{e} tangent de $M$ sont
  dites {\em projectivement \'{e}quivalentes} (que l'on note $\nabla'\sim\nabla$)
  lorsqu'il existe une $1$-forme
  $\alpha\in\Ginf(T^*M)$ telle que
  \beq \lbl{EqDefProjEq}
    \nabla'_X Y ~=~\nabla_X Y ~+~ \alpha(X)\hspace{1mm}Y
                                ~+~ \alpha(Y)\hspace{1mm}X.
  \end{equation}
  quels que soient les champs de vecteurs $X,Y\in\Ginf(TM)$.
\edefi
Il est \'{e}vident que $\sim$ est une relation d'\'{e}quivalence, et on appelle la
classe d'\'{e}quivalence de $\nabla$ la {\em classe projective de $\nabla$}.
 Une classe
d'\'{e}quivalence de connexions sans torsion est \'{e}galement dite une {\em structure
projective} sur $M$, voir \cite{Kob72}, p.147, Prop. 7.2.
\bprop
 Soient $\nabla$ et $\nabla'$ deux connexions sans torsion sur $M$.
 Alors les deux \'{e}nonc\'{e}s suivants sont \'{e}quivalents:
 \ben
  \item $\nabla$ et $\nabla'$ sont projectivement \'{e}quivalentes.
  \item $\nabla$ et $\nabla'$ {\em ont les m\^{e}mes g\'{e}od\'{e}siques \`{a} r\'{e}param\'{e}trage
    pr\`{e}s} dans le sens
      suivant: soit $x\in M$, $U_x$ un voisinage ouvert
      de l'origine de l'espace tangent $T_xM$, qui est dans
      l'intersection des domaines des applications exponentielles de
      $\nabla$ et $\nabla'$, et $v\in U_x$; alors pour tout triplet
      $(x,U_x,v)$ il existe $\epsilon,\epsilon'>0$,
      $0<c,c'\in\real$ et un diff\'{e}omorphisme $g:I:=]-\epsilon,c[\ra
      I':=]-\epsilon',c'[$ avec $g(0)=0$ et $(dg/dt)(0)=1$ tel que la
      g\'{e}od\'{e}sique $\gamma:I\ra M$ par rapport \`{a} $\nabla$ \'{e}manant de $x$
      \`{a} vitesse initiale $v$ est donn\'{e}e par $\gamma'\circ g$ o\`{u}
      $\gamma'$ est la g\'{e}od\'{e}sique $I'\ra M$ par rapport \`{a} $\nabla'$ \'{e}manant de $x$
      \`{a} vitesse initiale $v$.
 \een
\eprop
\bbew
 Soit $S\in \Ginf(TM\otimes T^*M\otimes T^*M)$ le champ de tenseur d\'{e}fini
 par $S(X,Y):=\nabla'_XY-\nabla_XY$. Puisque $\nabla$ et $\nabla'$ sont
 de torsion nulle, le champ $S$ est sym\'{e}trique, $S(X,Y)=S(Y,X)$.
 Soit $\gamma':I'\ra M$ une courbe de classe $\Cinf$ et $g:I\ra I'$ de
 classe $\Cinf$. Pour une courbe donn\'{e}e $c:I\ra M$, on rappelle
 son acc\'{e}l\'{e}ration $D^2c/dt^2:I\ra TM$ par rapport \`{a} $\nabla$, c.-\`{a}-d.
 \[
     \frac{D^2 c}{dt^2}
         :=(c^*\nabla)_{\frac{\partial}{\partial t}}\frac{dc}{dt}.
 \]
 o\`{u} $c^*\nabla$ est la connection $\nabla$ retir\'{e}e \`{a} $c^*TM$ (la
 restriction du fibr\'{e} tangent \`{a} $c$) et la vitesse $\frac{dc}{dt}$ est
 consid\'{e}r\'{e}e comme une section de $c^*TM$. Bien s\^{u}r, $c$ est une g\'{e}od\'{e}sique
 si et seulement si son acc\'{e}l\'{e}ration s'annule. Alors
 \bea
    \frac{D^2(\gamma'\circ g)}{dt^2} & = &
          \left(\frac{dg}{dt}\right)^2
          \left(\frac{{D'}^2\gamma'}{d\tau^2}\circ g\right)\nonumber \\
          & &  + \hspace{1mm}\frac{d^2g}{dt^2}
            \left(\frac{d\gamma'}{d\tau}\circ g\right)\hspace{1mm}
            -\hspace{1mm}\left(\frac{dg}{dt}\right)^2
                 S_{\gamma'\circ g}\left(\frac{d\gamma'}{d\tau}\circ g,
                                      \frac{d\gamma'}{d\tau}\circ
                                      g\right)
                \lbl{EqAccGeo}
 \eea
 Si l'\'{e}nonc\'{e} 1. est satisfait, et si $\gamma'$ est une g\'{e}od\'{e}sique de
 $\nabla'$, alors le membre droit de l'\'{e}quation ci-dessus est
 \'{e}gal \`{a} $(d\gamma'/d\tau)\circ g$ multipli\'{e} par
 \[
     \frac{d^2g}{dt^2}
      - 2 \left(\frac{dg}{dt}\right)^2
        \alpha_{\gamma'\circ g}\left(\frac{d\gamma'}{d\tau}\circ g\right)
        =: \ddot{g} -2F(g)\dot{g}^2
 \]
 o\`{u} $F(g):=\alpha_{\gamma'\circ g}\big((d\gamma'/d\tau)\circ g\big)$.
 On voit ais\'{e}ment que l'application
 \[
      \int_0^{g(t)}e^{-2\int_0^x F(y)dy}dx = t
 \]
 est inversible pour donner une fonction $g$ de classe $\Cinf$ \`{a} valeurs
 r\'{e}elles d\'{e}finie sur un intervalle ouvert
 contenant $0$ qui satisfait l'\'{e}quation diff\'{e}rentielle
 $0=\ddot{g} -2F(g)\dot{g}^2$ et les conditions initiales
 $g(0)=0$ et $\dot{g}(0)=1$. Alors, autour de $0$ la fonction $g$ est un
 diff\'{e}omorphisme, et la courbe r\'{e}param\'{e}tr\'{e}e
 $\gamma=\gamma'\circ g$ est une g\'{e}od\'{e}sique de $\nabla$.\\
 Si l'\'{e}nonc\'{e} 2. est vrai, alors les courbes $\gamma$ et $\gamma'$ sont des
 g\'{e}od\'{e}siques, et pour $t=0$ l'equation (\ref{EqAccGeo}) entra\^{\i}ne que
 $S_x(v,v)$ est proportionnel \`{a} $v$ pour tout $v\in T_xM$. On obtient l'\'{e}quation alg\'{e}brique
 $S_x(v,v)\wedge v=0~\forall v\in T_xM$ qui donne apr\`{e}s une polarisation
 \[
     0=S_x(v_1,v_2)\wedge v_3 + S_x(v_3,v_1)\wedge v_2 +
     S_x(v_2,v_3)\wedge v_1~~~\forall v_1,v_2,v_3\in T_xM
 \]
 La trace par rapport \`{a} $v_3$ et le deuxi\`{e}me facteur du produit tensoriel
 donne l'\'{e}quation
 \[
    S_x(v_1,v_2)=
    \frac{1}{m+1}\big((\mathrm{tr}S_x)(v_1)~v_2~+~(\mathrm{tr}S_x)(v_2)~v_1\big)
               ~~~\forall v_1,v_2\in T_xM
 \]
 o\`{u} $(\mathrm{tr}S_x)(v_1):=\mathrm{trace}\big(v\mapsto S_x(v_1,v)\big)$. Ceci
 montre que $S$ est de la forme (\ref{EqDefProjEq}) o\`{u} la $1$-forme
 $\alpha$ est donn\'{e}e par $\frac{1}{m+1}\mathrm{tr}S$. Par cons\'{e}quent,
 $\nabla'\sim\nabla$.
\ebew

Pour une connexion sans torsion $\nabla$ dans le fibr\'{e} tangent on rappelle
les d\'{e}finitions pour le tenseur de courbure
$R\in\Ginf(TM\otimes T^*M\otimes \Lambda^2 T^*M)$, la trace de la courbure
$\mathrm{tr}R\in\Ginf(\Lambda^2 T^*M)$ et le tenseur de Ricci
$Ric\in \Ginf(T^*M\otimes T^*M)$ o\`{u}
$X,Y,Z\in\Ginf(TM)$:
\bea
  R(X,Y)Z  & := & \nabla_X\nabla_YZ -\nabla_Y\nabla_XZ - \nabla_{[X,Y]}Z \\
  (\mathrm{tr}R)(X,Y) & := & \mathrm{trace}\big(Z\mapsto R(X,Y)Z\big)      \\
  Ric(Z,Y)   & := & \mathrm{trace}\big(X\mapsto R(X,Y)Z\big)
\eea
De la premi\`{e}re identit\'{e} de Bianchi ($0=R(X,Y)Z+R(Y,Z)X+R(Z,X)Y$) on d\'{e}duit
\beq
   (\mathrm{tr}R)(X,Y) = Ric(X,Y)-Ric(Y,X).
\end{equation}
Dans la proposition suivante on calcule ces champs de tenseurs pour une
connexion $\nabla'$ qui soit projectivement \'{e}quivalente \`{a} $\nabla$:
\bprop
  Soit $\nabla$ une connexion sans torsion dans le fibr\'{e} tangent de $M$,
  $\alpha$ une $1$-forme, et $\nabla'$ une connexion sans torsion dans
  le fibr\'{e} tangent de $M$ qui soit projectivement \'{e}quivalente \`{a} $\nabla$
  via $\alpha$, voir eqn (\ref{EqDefProjEq}).\\
  Alors on les formules suivantes qui relient les tenseurs de courbure
  $R$ de $\nabla$ et $R'$ de $\nabla'$, leurs traces
  $\mathrm{tr}R$ et $\mathrm{tr}R'$ et les tenseurs de Ricci $Ric$ et
  $Ric'$ o\`{u} $X,Y,Z$ sont des champs de vecteurs:
  \bea
      R'(X,Y)Z & = & R(X,Y)Z + \alpha(Y)\alpha(Z)X - \alpha(X)\alpha(Z)Y
                                                           \nonumber \\
                   & & ~~\big((\nabla_X\alpha)(Y)-(\nabla_Y\alpha)(X)\big)Z
                          + (\nabla_X\alpha)(Z)~Y -(\nabla_Y\alpha)(Z)~X
                          \nonumber \\
                   & &                                     \\
      \mathrm{tr}R'(X,Y) & = & \mathrm{tr}R(X,Y) +
      (m+1)\big((\nabla_X\alpha)(Y)-(\nabla_Y\alpha)(X)\big) \\
      Ric'(X,Y)     & = & Ric(X,Y) + (m-1)\alpha(X)\alpha(Y)
                           + (\nabla_X\alpha)(Y)-m(\nabla_Y\alpha)(X)
                            \nonumber \\
                   & &
  \eea
\eprop
\bbew
 Les trois formules s'obtiennent par un long calcul direct.
\ebew
\bprop
 Soient $\Phi:M\ra N$ une immersion entre deux vari\'{e}t\'{e}s diff\'{e}rentiables de
 dimension $m$. Soient $\nabla$ et $\nabla'$ deux connexions sans torsion
 dans le fibr\'{e} tangent de $N$.\\
 Si $\nabla$ et $\nabla'$ sont projectivement \'{e}quivalentes via une
 $1$-forme $\alpha\in\Ginf(T^*N)$, alors les connexions retir\'{e}es
 $\Phi^*\nabla$ et $\Phi^*\nabla'$ sont \'{e}quivalentes via la $1$-forme
 retir\'{e}e $\Phi^*\alpha$.
\eprop
\bbew
Soient $\omega$ et $\omega'$ les formes de connexion sur le fibr\'{e} des
rep\`{e}res lin\'{e}aires $P^1N$ qui correspondent aux connexions $\nabla$ et
$\nabla'$ sans torsion.
On rappelle
la $1$-forme canonique $\theta_N$ \`{a} valeurs dans $\real^m$ sur $P^1N$.
Alors pour $p\in P^1N$, $X\in T_pP^1N$ et $v\in\real^m$ on a
\beq \lbl{EqConnPProj}
   \omega'_p(X)(v) = \omega_p(X)(v) + \big((\pi^1_0)^*\alpha\big)_p(X)~v
                         + f_\alpha(p)\big(v\big)\theta_{Np}(X)
\end{equation}
voir \cite{Kob72}, p.145, paragraphe 7.
On peut consid\'{e}rer le pull-back de
l'\'{e}quation ci-dessus par $P^1\Phi$: puisque $(P^1\Phi)^*\theta_N=\theta_{M}$
et $(P^1\Phi)^*f_\alpha=f_{\Phi^*\alpha}$ ceci montre que les connexions
retir\'{e}es $\Phi^*\nabla$ et $\Phi^*\nabla'$ sont projectivement
\'{e}quivalentes par rapport \`{a} la $1$-forme $\Phi^*\alpha$.
\ebew

Le probl\`{e}me suivant me semble int\'{e}ressant
dans le cadre des prescriptions d'ordre naturelles:

\bprob
 Soient $F$ et $F'$ deux foncteurs de fibr\'{e}s vectoriels d'ordre $1$.
 Que peut-on dire sur l'existence et l'unicit\'{e} d'une prescription d'ordre
 naturelle $\rho$ qui soit {\em projectivement invariante}, c.-\`{a}-d.
 \beq
       \rho[\nabla](A)=\rho[\nabla'](A)
\end{equation}
 quelles que soient les connexions $\nabla$ et $\nabla'$ sans torsion
 projectivement \'{e}qui\-va\-lentes et la section
 $A\in\Ginf\big(STM\otimes Hom(FM,F'M)\big)$?
\eprob

L'inspiration pour ce probl\`{e}me est venue de la conjecture suivante de
Pierre Lecomte:

\bconj[P.B.A.~Lecomte] Soit $b\in\real$ et soient les deux foncteurs de fibr\'{e}s
vectoriels $F$ et $F'$ d'ordre $1$ \'{e}gaux \`{a} $|\Lambda^m T^*|^b$.\\
Alors il existe une pr\'{e}scription d'ordre naturelle qui soit
projectivement invariante.
\econj

\section{Rel\`{e}vement projectivement invariant d'une connexion sans torsion
         au fibr\'{e} principal de densit\'{e}s}

Soit $M$ une vari\'{e}t\'{e} de dimension $m$ et $\nabla$ une connexion sans
torsion dans le fibr\'{e} tangent de $M$. Soit $\omega$ la $1$-forme de
connexion correspondante dans le fibr\'{e} des rep\`{e}res lin\'{e}aires et
$\omega^a$ la $1$-forme de connexion induite dans le fibr\'{e} principal de
$a$-densit\'{e}s $\tilde{M}^a$ pour $a\in \real,a\neq 0$. On rappelle la
notation $X^h$ pour le rel\`{e}vement horizontal \`{a} $\tilde{M}^a$ d'un
champ de vecteurs $X$ sur $M$, et $\mathsf{E}$ pour le champ d'Euler
sur $\tilde{M}^a$.
\bdefi
 Soit $\bar{\nabla}$ une connexion sans torsion dans le fibr\'{e} tangent de
 $\tilde{M}^a$. Elle est dite {\em invariante} lorsque
 \beq \lbl{EqDefLENabla}
     \big(L_{\mathsf{E}}\bar{\nabla}\big)(Z,Z'):=
      [\mathsf{E},\bar{\nabla}_ZZ']-\bar{\nabla}_{[\mathsf{E},Z]}Z'
                            -\bar{\nabla}_Z[\mathsf{E},Z'] = 0
  \end{equation}
  quels que soient les champs de vecteurs $Z,Z'$ sur $\tilde{M}^a$.
\edefi
On rappelle que $L_{\mathsf{E}}\bar{\nabla}$ est un champ de tenseurs dans
$\Ginf(T\tilde{M}^a\otimes T^*\tilde{M}^a\otimes T^*\tilde{M}^a)$.
\bdefi
 On appelle un {\em rel\`{e}vement naturel de connexions sans torsion pour le
 foncteur $\tilde{F}^a$} une collection d'applications
 $A_M:Q_\tau P^1M\ra Q_\tau P^1 \tilde{M}^a$ (pour toute vari\'{e}t\'{e} $M$ de
 dimension $m$) telle que
 \ben
  \item Pour toute immersion $\Phi:M\ra N$ avec $\dim N=m$ on a
     \beq \lbl{EqDefConnNat}
         A_M(\Phi^*\nabla')=
             {\tilde{\Phi}^{a*}}A_N(\nabla')
     \end{equation}
     quelle que soit la connexion $\nabla'$ dans le fibr\'{e} tangent de $N$.
  \item Toutes les connexions $A_M(\nabla)$ sont invariantes.
 \een
 On appelle $A_M(\nabla)$ {\em le rel\`{e}vement naturel de la connexion $\nabla$
 \`{a} $\tilde{M}^a$ par rapport \`{a} $A$}.
\edefi
\bprop \lbl{PRelNatConn}
 Soient $X$ et $Y$ deux champs de vecteurs sur $M$ et
 $\mu,\nu,\rho\in\real$.
 \ben
  \item Alors les formules
 suivantes d\'{e}terminent un rel\`{e}vement naturel $\bar{\nabla}$ d'une connexion sans
 torsion $\nabla$ dans le fibr\'{e} tangent de $\tilde{M}^a$:
  \bea
    \bar{\nabla}_{X^h}Y^h & := &
                   (\nabla_XY)^h + \frac{a}{2}
                   {\tau^a}^*\big((\mathrm{tr}R)\big(X,Y\big)\big)\hspace{1mm}
                              \mathsf{E}
                        \nonumber \\
                            &    & ~~~
                            +\mu{\tau^a}^*\big(Ric(X,Y)+Ric(Y,X)\big)\hspace{1mm}
                            \mathsf{E}
                      \\
    \bar{\nabla}_{X^h}\mathsf{E} & := & \nu X^h, \\
    \bar{\nabla}_{\mathsf{E}}X^h & := & \nu X^h, \\
    \bar{\nabla}_{\mathsf{E}}\mathsf{E} & := & \rho\mathsf{E}.
  \eea
  \item Tout rel\`{e}vement naturel de connexions sans torsion pour le
  foncteur $\tilde{F}^a$ est de la forme pr\'{e}c\'{e}dente pour un choix arbitraire
  de $\mu,\nu,\rho\in\real$.
 \een
\eprop
La d\'{e}monstration de cette proposition est tr\`{e}s longue (surtout l'\'{e}nonc\'{e} de
l'unicit\'{e})
et sans importance pour la suite. On l'a donc mise dans l'appendice
\ref{SecApp3}.

Nous allons montrer maintenant l'existence et l'unicit\'{e} d'un rel\`{e}vement
naturel de connexions sans torsion $\nabla\mapsto\tilde{\nabla}$ tel que
la connexion $\tilde{\nabla}$ est projectivement invariante:

\bsat\lbl{TExUnConnProjEq}
 Soit $a$ un nombre r\'{e}el non nul et $m$ un entier positif avec $m\geq 2$.
 Il existe un unique rel\`{e}vement naturel de
 connexions sans torsion $\nabla\mapsto\tilde{\nabla}$ pour le
 foncteur $\tilde{F}^a$ tel que $\tilde{\nabla}=\tilde{\nabla'}$ si
 $\nabla$ et $\nabla'$ sont projectivement \'{e}quivalentes. $\tilde{\nabla}$
 est donn\'{e}e par la formule suivante (dans la notation de la proposition
 \ref{PRelNatConn}):
 \bea
    \tilde{\nabla}_{X^h}Y^h & := &
                   (\nabla_XY)^h + {\textstyle \frac{a}{2}}
                   {\tau^a}^*\big((\mathrm{tr}R)\big(X,Y\big)\big)\hspace{1mm}
                              \mathsf{E}
                        \nonumber \\
                            &    & ~~~
      -{\textstyle \frac{a}{2}\frac{m+1}{m-1}}
        {\tau^a}^*\big(Ric(X,Y)+Ric(Y,X)\big)\hspace{1mm}
                            \mathsf{E}
                      \\
    \tilde{\nabla}_{X^h}\mathsf{E} & := & {\textstyle \frac{1}{a(m+1)}}X^h, \\
    \tilde{\nabla}_{\mathsf{E}}X^h & := & {\textstyle \frac{1}{a(m+1)}}X^h, \\
    \tilde{\nabla}_{\mathsf{E}}\mathsf{E} & := & {\textstyle \frac{1}{a(m+1)}}
                             \mathsf{E}.
  \eea
\esat
\bbew
 D'apr\`{e}s la proposition \ref{PRelNatConn} on n'a qu'\`{a} \'{e}tudier la famille
 de connexions y mentionn\'{e}e. Soit $M$ une vari\'{e}t\'{e} de dimension $m$,
 $\nabla$ une connexion et $\alpha$ une $1$-forme. Soit $\nabla'$ la
 connexion $(X,Y)\mapsto \nabla_X Y +\alpha(X)Y+\alpha(Y)X$, alors
 $\nabla$ et $\nabla'$ sont projectivement \'{e}quivalentes. Soit $X^{h'}$ le
 rel\`{e}vement horizontal du champ de vecteurs $X$ sur $M$ \`{a} $\tilde{M}^a$
 par rapport \`{a} $\nabla'$. Soit $\omega^a$ la $1$-forme de connexion sur
 $\tilde{M}^a$ correspondant \`{a} $\nabla$ et ${\omega'}^a$ la $1$-forme de
 connexion sur $\tilde{M}^a$ correspondant \`{a} $\nabla'$. Les \'{e}quations
 (\ref{EqConnPTilde}) et (\ref{EqConnPProj}) nous donnent la formule
 suivante pour tout $p\in P^1M, v\in T_pP$:
 \beas
    {\omega'}^a_{\Xi^a(p)}\big(T_p\Xi^a\hspace{1mm}v\big)
     & = &  -a\mathrm{tr}\big({\omega'}_p(v)\big) \\
     & = &  -a\mathrm{tr}\big(\omega_p(v)\big)
       -a \big({\Xi^a}^*(\pi^1_0)^*\alpha\big)_{p}(v)
                \mathrm{tr}(\mathbf{1})
           -a\langle f_\alpha(p),\theta_{Mp}(v)\rangle \\
     & = &  \omega^a_{\Xi^a(p)}\big(T_p\Xi^a\hspace{1mm}v\big)
           -a(m+1)({\tau^a}^*\alpha)_{\Xi^a(p)}\big(T_p\Xi^a\hspace{1mm}v\big),
 \eeas
 alors
 \[
    {\omega'}^a = {\omega}^a -a(m+1){\tau^a}^*\alpha.
 \]
 Puisque la diff\'{e}rence $X^{h'}-X^h$ est un champ vertical on en d\'{e}duit la
 formule
 \beq \lbl{EqRelHorPrim}
   X^{h'} = X^{h} + a(m+1){\tau^a}^*\big(\alpha(X)\big)\mathsf{E}.
 \end{equation}
 Fixons maintenant $\mu,\rho,\nu$ et essayons de d\'{e}terminer leurs valeurs
 pour que $\nabla\mapsto\tilde{\nabla}$ soit projectivement invariant:
 Il est clair que
 \beq \lbl{EqConnProjVV}
    \widetilde{\nabla'}_{\mathsf{E}}\mathsf{E}
        = \rho \mathsf{E} = \widetilde{\nabla}_{\mathsf{E}}\mathsf{E}.
 \end{equation}
 Ensuite, on a --gr\^{a}ce \`{a} l'invariance de ${\tau^a}^*\big(\alpha(X)\big)$:
 \bea
    \widetilde{\nabla'}_{\mathsf{E}}X^h & = &
         \widetilde{\nabla'}_{\mathsf{E}}\left(X^{h'}
               -a(m+1){\tau^a}^*\big(\alpha(X)\big)\mathsf{E}\right)
                          \nonumber \\
               & = & \nu X^{h'}
                    - \rho a (m+1){\tau^a}^*\big(\alpha(X)\big)\mathsf{E}
                    \nonumber \\
               & = & \nu X^h
                    - (\rho-\nu)a (m+1){\tau^a}^*\big(\alpha(X)\big)\mathsf{E}
                    \nonumber \\
               & = & \widetilde{\nabla}_{\mathsf{E}}X^h
                    - (\rho-\nu)a
                    (m+1){\tau^a}^*\big(\alpha(X)\big)\mathsf{E}.
                        \lbl{EqConnProjHV}
 \eea
 Ensuite, on calcule les champs horizontaux:
 \beas
  \widetilde{\nabla'}_{X^h}Y^h & = &
         \widetilde{\nabla'}_{X^{h'}-a(m+1){\tau^a}^*(\alpha(X))\mathsf{E}}
          \big(Y^{h'}-a(m+1){\tau^a}^*(\alpha(Y))\mathsf{E}\big) \\
          & = &
          \big(\nabla'_{X}Y\big)^{h'}
             + \frac{a}{2}{\tau^a}^*\big((\mathrm{tr}R')(X,Y)\big)\mathsf{E}
             +\mu{\tau^a}^*\big(Ric'(X,Y)+Ric'(Y,X)\big)\mathsf{E} \\
          &   &
           -a(m+1){\tau^a}^*\big(X(\alpha(Y))\big)\mathsf{E}
           -\nu a(m+1){\tau^a}^*\big(\alpha(Y)\big) X^{h'}  \\
          &   &
           -\nu a(m+1){\tau^a}^*\big(\alpha(X)\big) Y^{h'}
           +\rho a^2(m+1)^2
           {\tau^a}^*\big(\alpha(X)\alpha(Y)\big)\mathsf{E}\\
          & = &
          \big(\nabla_{X}Y\big)^{h}
             + \frac{a}{2}{\tau^a}^*\big((\mathrm{tr}R)(X,Y)\big)\mathsf{E}
             +\mu{\tau^a}^*\big(Ric(X,Y)+Ric(Y,X)\big)\mathsf{E} \\
          &   &
             +a(m+1){\tau^a}^*\big(\alpha(\nabla_XY)\big)\mathsf{E} \\
          &   &
             + {\tau^a}^*\big(\alpha(X)\big) Y^{h}
             + a(m+1){\tau^a}^*\big(\alpha(X)\alpha(Y)\big)\mathsf{E} \\
          &   &
             + {\tau^a}^*\big(\alpha(Y)\big) X^{h}
             + a(m+1){\tau^a}^*\big(\alpha(Y)\alpha(X)\big)\mathsf{E} \\
          &   &
             + \frac{a}{2}(m+1){\tau^a}^*\big((\nabla_X\alpha)(Y)-
                                  (\nabla_Y\alpha)(X)\big)\mathsf{E} \\
          &   &
              +\mu(m-1)
              {\tau^a}^*\big(2\alpha(X)\alpha(Y)-(\nabla_X\alpha)(Y)
                                    -(\nabla_Y\alpha)(X)\big)
                                     \mathsf{E}\\
          &   &
              -a(m+1){\tau^a}^*\big((\nabla_X\alpha)(Y)
                                    +\alpha(\nabla_XY)\big)\mathsf{E} \\
          &   &
              -\nu a(m+1)\big({\tau^a}^*\big(\alpha(Y)\big) X^{h}
                               + {\tau^a}^*\big(\alpha(X)\big) Y^{h}\big)
                                 \\
          &   & a^2(m+1)^2\big(\rho-2\nu\big)
             {\tau^a}^*\big(\alpha(X)\alpha(Y)\big)\mathsf{E},
 \eeas
alors
\bea
 \widetilde{\nabla'}_{X^h}Y^h & = &
                           \widetilde{\nabla}_{X^h}Y^h \nonumber \\
                           &    &
                            +\big(1-\nu a (m+1)\big)
                            \big({\tau^a}^*\big(\alpha(Y)\big) X^{h}
                               + {\tau^a}^*\big(\alpha(X)\big) Y^{h}\big)
                               \nonumber \\
                           &    &
                 -\big(\frac{a}{2}(m+1)+\mu(m-1)\big)
                     {\tau^a}^*\big((\nabla_X\alpha)(Y)
                       +(\nabla_X\alpha)(Y)\big)\mathsf{E} \nonumber \\
                           &    &
        +\left(2a(m+1)+2\mu(m-1)+(\rho-2\nu)a^2(m+1)^2\right)
               {\tau^a}^*\big(\alpha(X)\alpha(Y)\big)\mathsf{E}.
                                         \nonumber \\
                           &    &      \lbl{EqConnProjHH}
\eea
Les \'{e}quations (\ref{EqConnProjHV}) et (\ref{EqConnProjHH}) montrent que
$\nabla\mapsto\tilde{\nabla}$ est projectivement invariant si et seulement
si
\[
  \begin{array}{ccc}
     \mu=-\frac{a}{2}\frac{m+1}{m-1}, &
        \nu= \frac{1}{a(m+1)}, & \rho = \frac{1}{a(m+1)},
  \end{array}
\]
ce qui montre le th\'{e}or\`{e}me.
\ebew

Ma motivation principale pour trouver la connexion naturelle $\tilde{\nabla}$
dans le th\'{e}or\`{e}me pr\'{e}c\'{e}dent a \'{e}t\'{e}
l'\'{e}tude de l'exemple g\'{e}om\'{e}trique suivant:

\bexem \lbl{ExemSmI}
 Soit $m$ un entier positif avec $m\geq 2$,
soit $\tilde{M}=\real^{m+1}\setminus \{0\}$, soit $\langle~,~\rangle$
le produit scalaire canonique de $\real^{m+1}$ et
$|y|=\sqrt{\langle y,y\rangle}$ la norme associ\'{e}e, soit $M$ la sph\`{e}re
$S^m$ et soit
$\tau:\tilde{M}\ra M$ la projection radiale $y\mapsto \frac{y}{|y|}$.
Un champ de vecteurs $X$ sur $M$ est une application de classe $\Cinf$
de $M$ dans $\real^{m+1}$ telle que $\langle X(x),x\rangle=0$ quel que
soit $x\in M$. Soit $\tilde{\nabla}$ la d\'{e}riv\'{e}e covariante canonique dans
$\real^{m+1}$. La connexion de Levi-Civita $\nabla$ induite par la restriction
$\mathsf{g}$ du produit scalaire $\langle~,~\rangle$ \`{a} $TM$ se calcule
par projection
\[
   (\nabla_XY)_x :=
   (\tilde{\nabla}_XY)_x-\langle(\tilde{\nabla}_XY)_x,x\rangle x
   =(\tilde{\nabla}_XY)_x+\mathsf{g}_x(X_x,Y_x)x
\]
o\`{u} $X,Y$ sont deux champs de vecteurs sur $M$ et $\tilde{\nabla}_XY$ est
bien d\'{e}fini sur $M$ parce que $X_x$ est un vecteur tangent de $M$.
On peut d\'{e}finir un rel\`{e}vement horizontal $X^h$ de $X$ \`{a} $\tilde{M}$ en
demandant que $\langle X^h_y,y\rangle = 0$ et $T_y\tau X^h_y
=X_{\tau(y)}$. Il en r\'{e}sulte que
\[
    X^h_y = |y|X_{y/|y|}.
\]
Puisque le champ d'Euler de $\tilde{M}$ est donn\'{e} par $\mathsf{E}_y=y$ on
arrive facilement aux formules suivantes:
\bea
    \tilde{\nabla}_{X^h}Y^h & = & (\nabla_XY)^h
                                  - \tau^*\big(\mathsf{g}_x(X_x,Y_x)\big)
                                      \mathsf{E}, \\
\tilde{\nabla}_{X^h}\mathsf{E}=\tilde{\nabla}_{\mathsf{E}}{X^h}
                            & = & X^h, \\
      \tilde{\nabla}_{\mathsf{E}}\mathsf{E} & = & \mathsf{E}.
\eea
Puisque $\nabla$ pr\'{e}serve la m\'{e}trique riemannienne $\mathsf{g}$ sur $M$,
le champ de tenseur $\mathrm{tr}R$ s'annule, et pour la sph\`{e}re --dont la
courbure sectionnelle est \'{e}gale \`{a} $1$-- il vient que le tenseur de Ricci
est sym\'{e}trique et \'{e}gal \`{a} $(m-1)\mathsf{g}$, voir par exemple
\cite{KN63}, p.204, Theorem 3.1 et p.203, Corollary 2.3. Alors on voit
que les formules pr\'{e}cedentes co\"{\i}ncident avec celles du th\'{e}or\`{e}me
\ref{TExUnConnProjEq} pour le choix $a:=\frac{1}{m+1}$.
\eexem

\section{Rel\`{e}vement naturel projectivement invariant des champs de tenseurs
         sym\'{e}triques
         au fibr\'{e} principal de densit\'{e}s}

Soit $M$ une vari\'{e}t\'{e} diff\'{e}rentiable, $\nabla$ une connexion sans torsion
dans le fibr\'{e} tangent de $M$ et $A$ un champ de tenseurs dans
$\Ginf(|\Lambda^m T^*M|^c\otimes S^kTM)$. On rappelle que
$|\Lambda^m T^*M|^c$ est isomorphe \`{a} $Hom(|\Lambda^m T^*M|^b,|\Lambda^m
T^*M|^{b+c})$. Dans ce paragraphe on veut rel\`{e}ver $A$ \`{a} $\tilde{M}^a$,
$a\neq 0$, d'une fa\c{c}on qui peut d\'{e}pendre de $\nabla$, mais seulement
de sa classe projective.

Pour tout entier positif $k$ et pour tout $y\in \tilde{M}^a$ on d\'{e}finit d'abord
un rel\`{e}vement horizontal $(~)^h_y:S^kT_{\tau(y)}M\otimes
|\Lambda^mT^*M|^c_{\tau(y)}$ dans $S^kT_y\tilde{M}^a$, d\'{e}fini pour
$v_1,\ldots, v_k\in T_{\tau(y)}M$
et $s\in |\Lambda^mT^*M|^c_{\tau(y)}$ par
\beq
  \big(s\otimes (v_1\vee\cdots\vee v_k)\big)^h_y
      :=(y^{-1}s)\big({v_1^h}_y\vee\cdots\vee{v_k^h}_y\big)\in S^kT_y\tilde{M}^a.
\end{equation}
Il en r\'{e}sulte le rel\`{e}vement horizontal $A\mapsto A^h$ des sections $A$
appartenant \`{a} $\Ginf(|\Lambda^mT^*M|^c\otimes S^kTM)$ o\`{u} $A^h$ est un
\'{e}l\'{e}ment de $\Ginf(S^kT\tilde{M}^a)$.
\blem \lbl{LEquivRelHor}
 On a
 \beq
   L_\mathsf{E}A^h= -{\textstyle \frac{c}{a}}A^h.
 \end{equation}
\elem
\bbew
  Soit $g\in\real^+$. On a
  \beas
    \big(s\otimes (v_1\vee\cdots\vee v_k)\big)_{yg}
               & = &
                \big((yg)^{-1}s\big)
                \big({v_1^h}_{yg}\vee\cdots\vee{v_k^h}_{yg}\big) \\
               & = &
                 g^{-\frac{c}{a}}S^kT_yr_g\big(
                      (y^{-1}s)({v_1^h}_{y}\vee\cdots\vee{v_k^h}_{y})
   \eeas
   En particulier, pour $g=e^t$ il s'ensuit (o\`{u} $G_t$ d\'{e}signe le flot du
   champ d'Euler $\mathsf{E}$) que
   \[
        G_t^*A = e^{-\frac{c}{a}t}A,
   \]
   d'o\`{u} le r\'{e}sultat.
\ebew

Le multiple de la partie sym\'{e}trique du tenseur de
Ricci $Ric$ de $\nabla$ suivant sera important: soient $X,Y$ des champs de
vecteurs sur M:
\beq
    \mathsf{r}(X,Y):={\textstyle \frac{1}{2}\frac{1}{m-1}}\big(
                         Ric(X,Y)+Ric(Y,X)\big).
\end{equation}
On calcule la divergence d'un rel\`{e}vement horizontal:
\bprop \lbl{PTildeDiv}
 Soient $j,l$ deux entiers positifs et soit
 $A\in \Ginf(|\Lambda^mT^*M|^c\otimes S^jTM)$. Alors:
 \bea
   \widetilde{\mathsf{Div}}(A^h\vee \mathsf{E}^{\vee l})
             & = & (\mathsf{Div}A)^h \vee \mathsf{E}^{\vee l}
                  -2a(m+1)\big(i(\mathsf{r})A\big)^h\vee \mathsf{E}^{\vee (l+1)}
                       \nonumber \\
             &   & + {\textstyle \frac{l(2j+l+m-(m+1)c)}{a(m+1)}}
                       A^h\vee \mathsf{E}^{\vee (l-1)}
 \eea
\eprop
\bbew
 Il suffit de montrer cette \'{e}quation localement: soit $e_1,\ldots,e_m$ une
 base locale du fibr\'{e} tangent de $M$ et soit $e^1,\ldots,e^m$ sa base duale.
 On rappelle que la $1$-forme de connexion $\omega^a$ sur $\tilde{M}^a$
 est telle que
 $\omega^a(e^h_i)=0$ $\forall 1\leq i\leq m$ et $\omega^a(\mathsf{E})=1$.
 On peut supposer que $A$ est de la forme suivante: soit $\varphi$ une
 $c$-densit\'{e} locale, $X$ un champ de vecteurs local et
 \[
     A=\varphi \otimes X^{\vee j}\mathrm{~~~alors~~}
       A^h=\tilde{\varphi}{X^h}^{\vee j}.
 \]
 Alors
 \beas
   \mathsf{Div}A & = & \sum_{i=1}^m i(e^i)\nabla_{e_i}
                       (\varphi \otimes X^{\vee j})
                  =  \sum_{i=1}^m i(e^i)
                 \big(\nabla_{e_i}\varphi\otimes X^{\vee j}
                       +j\varphi\otimes \nabla_{e_i}X\vee X^{\vee (j-1)}
                             \big)\\
                 & = & j\nabla_X\varphi\otimes X^{\vee (j-1)}
                       + j\varphi\otimes (\mathsf{Div}X)X^{\vee (j-1)}
                                \\
                 &   & ~~~~~~~~~~~
                       + j(j-1)\varphi\otimes (\nabla_XX)\vee X^{\vee
                       (j-2)}.~~~~~~~~~~~~~~~~~~~~~~~~~~~~~~~~~~(*)
 \eeas
 En utilisant la base locale $(e_1^h,\ldots,e_m^h,\mathsf{E})$ de
 $T\tilde{M}^a$ et la base duale
 $(\tau^{a*}e^1,\ldots,\tau^{a*}e^m$, $\omega^a)$ on peut calculer
 la divergence $\widetilde{\mathsf{Div}}$ en notant que
 \[
   \tilde{\nabla}_{e_i^h}\tilde{\varphi}=e_i^h(\tilde{\varphi})=
 \widetilde{\nabla_{e_i}\varphi}~~~\mathrm{et}~~~
 \tilde{\nabla}_{\mathsf{E}}\tilde{\varphi}=\mathsf{E}(\tilde{\varphi})=
 -(c/a)\tilde{\varphi}
 \]
 alors, avec $\bar{a}:=a(m+1)$:
 \beas
    \widetilde{\mathsf{Div}}(A^h\vee \mathsf{E}^{\vee l})
          & = &
          \sum_{i=1}^m i(\tau^{a*}e^i)\left(
                        \widetilde{\nabla_{e_i}\varphi}\hspace{1mm}
                          {X^h}^{\vee j}\vee\mathsf{E}^{\vee l}~
                        +~j\tilde{\varphi}\hspace{1mm}(\nabla_{e_i}X)^h\vee
                          {X^h}^{\vee (j-1)}\vee\mathsf{E}^{\vee l}\right.\\
          & & ~~~~~~~~~~~~~~+j\tilde{\varphi}\hspace{1mm}
                               \tau^{a*}\big(\frac{a}{2}(\mathrm{tr}R)(e_i,X)
                                              -\bar{a}\mathsf{r}(e_i,X)\big)
                                {X^h}^{\vee (j-1)}\vee\mathsf{E}^{\vee (l+1)}
                                                                          \\
          & & ~~~~~~~~~~~~~~\left.+~{\textstyle \frac{l}{a(m+1)}}
                             \tilde{\varphi}\hspace{1mm}
                             {X^h}^{\vee j}\vee e_i^h\vee \mathsf{E}^{\vee (l-1)}
                                                                   \right)\\
          & & +~~~i(\omega^a)\left( {\textstyle \frac{j+l-(m+1)c}{a(m+1)}}
                       \tilde{\varphi}\hspace{1mm}{X^h}^{\vee j}\vee
                               \mathsf{E}^{\vee l}\right)\\
          & = &
            \left( j\widetilde{\nabla_X\varphi}\hspace{1mm}{X^h}^{\vee (j-1)}
                       + j\tilde{\varphi}\hspace{1mm}\tau^{a*}(\mathsf{Div}X)
                         {X^h}^{\vee (j-1)}\right. \\
          & & ~~~~~~~~~+~ \left.j(j-1)\tilde{\varphi}\hspace{1mm}(\nabla_XX)^h
                          \vee {X^h}^{\vee(j-2)}\right)\vee \mathsf{E}^{\vee l}
                                \\
          & & -j(j-1)\bar{a}\tilde{\varphi}\hspace{1mm}
                               \tau^{a*}\big(\mathsf{r}(X,X)\big)
                                {X^h}^{\vee (j-2)}\vee\mathsf{E}^{\vee (l+1)}
                                 \\
          & & +~{\textstyle \frac{l(j+m)}{a(m+1)}}
                             \tilde{\varphi}\hspace{1mm}
                             {X^h}^{\vee j}\vee \mathsf{E}^{\vee (l-1)}
                  ~+~{\textstyle \frac{l(j+l-(m+1)c)}{a(m+1)}}
                       \tilde{\varphi}\hspace{1mm}{X^h}^{\vee j}\vee
                               \mathsf{E}^{\vee (l-1)}.
 \eeas
 Avec l'\'{e}quation ($*$) et $i(\mathsf{r})X^{\vee j}
 =\frac{j(j-1)}{2}\mathsf{r}(X,X)X^{\vee (j-2)}$ on arrive \`{a} l'\'{e}nonc\'{e}
 d\'{e}sir\'{e}.
\ebew

\bdefi
 Soit $M$ une vari\'{e}t\'{e} diff\'{e}rentiable de dimension $m$ et $\nabla$ une
 connexion sans torsion dans le fibr\'{e} tangent de $M$. Pout tout entier
 positif $k$ et pour tous r\'{e}els $a,c$ tels que $a\neq 0$ on d\'{e}finit
 $\Ginf(S^kT\tilde{M}^a)^{-c/a}$ comme le sous-espace vectoriel de
 toutes les sections {\em $(-\frac{c}{a})$-\'{e}quivariantes} $B$ appartenant
 \`{a}
 $\Ginf(S^kT\tilde{M}^a)$, c.-\`{a}-d. qui satisfont \`{a} la condition suivante:
 \beq \lbl{EqEquivBca}
    L_{\mathsf{E}}B
                    =
                 {\textstyle -\frac{c}{a}}B.
 \end{equation}
 De plus, on \'{e}crira
 $\Ginf(S^kT\tilde{M}^a)^{-c/a}_{\tilde{\nabla}}$ pour le sous-espace
 vectoriel de toutes les sections $B$ dans $\Ginf(S^kT\tilde{M}^a)^{-c/a}$
 pour lesquelles la condition suivante soit satisfaite:
 \beq \lbl{EqDivBNul}
     \widetilde{\mathsf{Div}}B =  0.
 \end{equation}
\edefi
\blem \lbl{LDefPsi}
  Pour tout entier positif $k$ et tous r\'{e}els $a,c$, $a\neq 0$,
  l'application $\real$-lin\'{e}aire $\Psi$ de $\Ginf(S^kT\tilde{M}^a)^{-c/a}$
  dans $\Ginf(|\Lambda^mT^*M|^c\otimes S^kTM)$ donn\'{e}e par
  \beq \lbl{EqProjTildeAA}
  (\Psi B)_{\tau^a(y)}(\beta_1,\ldots,\beta_k)
     :=  \tilde{q}^a_c\big(y,B_y
         \big((T_y\tau^{a})^*\beta_1,\ldots,(T_y\tau^{a})^*\beta_k\big)\big),
 \end{equation}
 quels que soient $y\in \tilde{M}^a$ et
 $\beta_1,\ldots,\beta_k\in T^*_{\tau^a(y)}M$,
  est bien d\'{e}finie.
\elem
\bbew
 Soit $B\in\Ginf(S^kT\tilde{M}^a)$ et soit $B$ \'{e}quivariant par rapport
 \`{a} l'action du flot $G_t$ du champ d'Euler $\mathsf{E}$, c.-\`{a}-d.
 $G_t^*B=e^{-(c/a)t}B$ quel que soit $t\in\real$ ou
 \[
   B_{ye^t}= e^{-(c/a)t}(T_yG_t\otimes\cdots\otimes T_yG_t)B_y
 \]
 quel que soit $y\in \tilde{M}^a$, ce qui est \'{e}quivalent \`{a} l'\'{e}quation
 (\ref{EqEquivBca}).
 Alors pour tout $y\in \tilde{M}^a$ et pour tous
 $\beta_1,\ldots,\beta_k\in T_{\tau^a(y)}\tilde{M}^a$ le membre droit
 de l'\'{e}quation (\ref{EqProjTildeAA})
 ne d\'{e}pend que du point $x=\tau^a(y)$ et d\'{e}finit une unique section
 $A_k$ appaartenant \`{a} l'espace $\Ginf(|\Lambda^mT^*M|^c\otimes S^kTM)$.
 On obtient ainsi une application
 lin\'{e}aire $\Psi$ de $\Ginf(S^kT\tilde{M}^a)^{-c/a}$ dans
 $\Ginf(|\Lambda^mT^*M|^c\otimes S^kTM)$.
\ebew
\bsat \lbl{TConstrTildeANabla}
 Soient $k,m$ deux entiers positifs et $a,c\in\real$ tel que $a\neq 0$,
 \beq \lbl{EqCondcmk}
    m\geq 2~~~\mathrm{et}~~~
     c\not\in\left\{{\textstyle \frac{j+k+m}{m+1}}~|~
       j\in\nat~\mathrm{et}~0\leq j \leq k-1\right\}.
 \end{equation}
 Soit $M$ une vari\'{e}t\'{e} diff\'{e}rentiable de dimension $m$.\\
 \ben
 \item Pour toute connexion sans torsion $\nabla$ la restriction de
  l'application $\real$-lin\'{e}aire $\Psi$ d\'{e}finie dans le lemme
 \ref{LDefPsi} \`{a} $\Ginf(S^kT\tilde{M}^a)^{-c/a}_{\tilde{\nabla}}$
 est une bijection sur $\Ginf(|\Lambda^mT^*M|^c\otimes S^kTM)$.
  On note son inverse par $A\mapsto \tilde{A}[\nabla]$,
 pour toute section $A$ appartenant \`{a} $\Ginf(|\Lambda^mT^*M|^c\otimes S^kTM)$.
 \item Soit $\nabla'$ une connexion sans torsion dans le fibr\'{e} tangent de
 $M$ qui soit projectivement \'{e}quivalente \`{a} $\nabla$. Alors on a
 \beq
        \tilde{A}[\nabla]
                    =             \tilde{A}[\nabla']
 \end{equation}
 quel que soit $A\in\Ginf(|\Lambda^mT^*M|^c\otimes S^kTM)$.
 \item Soit $N$ une autre vari\'{e}t\'{e} diff\'{e}rentiable de dimension $m$ et
 $\Phi:N\ra M$ une immersion. Alors on a
 \beq \lbl{EqRelNatProjNaturalite}
   \tilde{\Phi}^{a*}\big(\tilde{A}[\nabla]\big)
              =\widetilde{\Phi^*A}[\Phi^*\nabla]
 \end{equation}
 quelles que soient la section $A\in\Ginf(|\Lambda^mT^*M|^c\otimes S^kTM)$
 et la connexion sans torsion $\nabla$ dans le fibr\'{e} tangent de $M$.
 \een
 On appelle
 $\tilde{A}[\nabla]$ le {\em rel\`{e}vement naturel projectivement
 invariant de $A$ par rapport \`{a} $\nabla$}.
\esat
\bbew
 Fixons une connexion sans
 torsion $\nabla$ dans le fibr\'{e} tangent de $M$.
 Soit $B\in\Ginf(S^kT\tilde{M}^a)$ et soit $A_k:=\Psi B$.
 D'apr\`{e}s le lemme
 \ref{LEquivRelHor} la section relev\'{e}e (par rapport \`{a} $\nabla$) $A_k^h$
 est aussi un \'{e}l\'{e}ment de $\Ginf(S^kT\tilde{M}^a)^{-c/a}$, et gr\^{a}ce \`{a} la
 d\'{e}finition du rel\`{e}vement horizontal d'un champ de vecteurs, il vient
 que $\Psi A_k^h = A_k$. Ceci montre que l'application $\Psi$ est
 toujours surjective.
 Alors $\Psi(B-A_k^h)=0$. La d\'{e}composition de l'espace
 tangent $T_y\tilde{M}^a$ en somme directe du sous-espace horizontal
 $H_y\tilde{M}^a$ et du sous-espace vertical $V_y\tilde{M}^a$ montre que
 le noyau de $T_y\tau^a\otimes\cdots\otimes T_y\tau^a$ restreint
 \`{a} $S^kT_y\tilde{M}^a$ est \'{e}gal \`{a} $V\tilde{M}^a\vee S^{k-1}T_y\tilde{M}^a$.
 Puisque le champ d'Euler est invariante et ne s'annule jamais sur $\tilde{M^a}$
 il s'ensuit qu'il existe une unique section $B_{k-1}$ appartenant
 \`{a} $\Ginf(S^{k-1}T\tilde{M}^a)^{-c/a}$ telle que
 $B=A_k^h+B_{k-1}\vee\mathsf{E}$. On peut r\'{e}p\'{e}ter le proc\'{e}d\'{e} pr\'{e}c\'{e}dent pour
 $B_{k-1}$, et par r\'{e}currence on obtient une suite $A_k,A_{k-1},\cdots,A_0$
 de sections o\`{u} $A_j\in\Ginf(|\Lambda^mT^*M|^c\otimes S^jTM)$, $0\leq j\leq
 k$ telle que
 \[
    B = A_k^h + A_{k-1}^h\vee\mathsf{E}+A^h_{k-2}\vee\mathsf{E}^{\vee 2}
           +\cdots+ A_1^h\vee \mathsf{E}^{\vee (k-1)}
           + \tilde{A_0}\vee\mathsf{E}^{\vee k}.
 \]
 Caculons la divergence de $B$ par rapport \`{a} $\widetilde{\nabla}$ d'apr\`{e}s
 le formule de la proposition \ref{PTildeDiv} o\`{u} $A_j:=0$ lorsque
 $j\leq -1$ ou $j\geq k+1$ et $\bar{a}:=a(m+1)$ et $\bar{m}:=m-(m+1)c$
 \beas
  \widetilde{\mathsf{Div}}B & = &
    \sum_{l=0}^k \widetilde{\mathsf{Div}}\big(A_{k-l}^h\vee
                               \mathsf{E}^{\vee l} \big) \\
      & = &
        \sum_{l=0}^k\big((\mathsf{Div}A_{k-l})^h \vee \mathsf{E}^{\vee l}
                  -2\bar{a}\big(i(\mathsf{r})A_{k-l}\big)^h\vee \mathsf{E}^{\vee (l+1)}
                           \\
             &   & ~~~~~~~~+~~{\textstyle \frac{l(2k-l+\bar{m})}{\bar{a}}}
                       A_{k-l}^h\vee \mathsf{E}^{\vee (l-1)}\big)\\
      & = & (\mathsf{Div}A_{k})^h +
               {\textstyle \frac{2k-1+\bar{m}}{\bar{a}} }A_{k-1}^h \\
      & &
         +\sum_{l=1}^{k-1}\left(\big(\mathsf{Div}A_{k-l}
                  -2\bar{a} i(\mathsf{r})A_{k-(l-1)}
                   +{\textstyle \frac{(l+1)(2k-l-1+\bar{m})}{\bar{a}} }
                       A_{k-(l+1)}\big)^h\vee \mathsf{E}^{\vee l}\right)
 \eeas
Alors il s'ensuit que $\widetilde{\mathsf{Div}}B=0$ si et seulement si
les \'{e}quations suivantes sont satisfaites quel que soit l'entier $l$,
$0\leq l\leq k-1$:
\bea
  A_{k-1} & = & -{\textstyle \frac{\bar{a}}{2k-1+\bar{m}} } \mathsf{Div}A_k
                             \lbl{EqRecAkmu}     \\
     & \vdots &  \nonumber \\
  A_{k-(l+1)} & = & -{\textstyle \frac{\bar{a}}{(l+1)(2k-l-1+\bar{m})} }
             \left(\mathsf{Div}A_{k-l}-2\bar{a}i(\mathsf{r})A_{k-(l-1)}
              \right) \lbl{EqRecAkmlpu}     \\
     & \vdots &  \nonumber \\
  A_0      & = & -{\textstyle \frac{\bar{a}}{k(k+\bar{m})} }
             \left(\mathsf{Div}A_1-2\bar{a}i(\mathsf{r})A_{2}
              \right) \lbl{EqRecAnul}
\eea
Gr\^{a}ce \`{a} la condition (\ref{EqCondcmk}) tous les d\'{e}nominateurs dans les
\'{e}quations pr\'{e}c\'{e}dentes sont non nuls. On voit par r\'{e}currence que
toutes les sections $A_{k-j}$, $1\leq j\leq k$, sont uniquement d\'{e}termin\'{e}es
par la section $A_k$. Alors le
sous-espace de toutes les solutions des \'{e}quations (\ref{EqDivBNul})
et (\ref{EqEquivBca}), $\Ginf(S^kT\tilde{M}^a)^{-c/a}_{\tilde{\nabla}}$,
est en bijection avec l'espace de tous les sections $A_k$, \`{a} savoir
$\Ginf(|\Lambda^mT^*M|^c\otimes S^kTM)$ via l'application $\Psi$. Il
s'ensuit que la restriction de l'application $\Psi$
\`{a} $\Ginf(S^kT\tilde{M}^a)^{-c/a}_{\tilde{\nabla}}$ est une bijection, ce qui
montre l'\'{e}nonc\'{e} (1).\\
(2) L'espace $\Ginf(S^kT\tilde{M}^a)^{-c/a}_{\tilde{\nabla}}$ d\'{e}pend de la
connexion $\tilde{\nabla}$ qui ne d\'{e}pend que de la classe
projective de $\nabla$. Puisque la projection $\Psi$ ne d\'{e}pend d'aucune
connexion, sa restriction $\Psi$ \`{a} $\Ginf(S^kT\tilde{M}^a)^{-c/a}_{\tilde{\nabla}}$
ne d\'{e}pend que de la classe
d'\'{e}quivalence de $\nabla$. Il en est de m\^{e}me pour son application
r\'{e}ciproque $A\mapsto\tilde{A}[\nabla]$, ce qui montre (2).\\
(3) Gr\^{a}ce \`{a} la naturalit\'{e} du rel\`{e}vement $\nabla\mapsto\tilde{\nabla}$,
(\ref{EqDefConnNat}), montr\'{e}e dans le th\'{e}or\`{e}me \ref{TExUnConnProjEq}, on a
pour tout $B\in\Ginf(S^kT\tilde{M}^a)$
\[
     \widetilde{\Phi^*\mathsf{Div}}\big(\Phi^*B\big)
       =\tilde{\Phi}^{a*}\big(\widetilde{\mathsf{Div}}B\big)
\]
ce qui montre que $\tilde{\Phi}^{a*}$ envoie le sous-espace
$\Ginf(S^kT\tilde{M}^a)^{-c/a}_{\tilde{\nabla}}$ dans le sous-espace
$\Ginf(S^kT\tilde{N}^a)^{-c/a}_{\widetilde{\Phi^*\nabla}}$ (o\`{u} l'on utilise
\'{e}videmment le fait que
$\tilde{\Phi}^a$ est un homomorphisme de fibr\'{e}s principaux). En \'{e}crivant
${\tau'}^a:\tilde{N}^a\ra N$ et $\Psi'$ pour l'application
(\ref{EqProjTildeAA}) pour la vari\'{e}t\'{e} $N$ on a pour tout $y'\in \tilde{N}^a$
et $\beta'_1,\ldots,\beta'_k\in T_{{\tau'}^a(y')}^*N$:
\beas
 \lefteqn{ (\Psi'\tilde{\Phi}^{a*}B)_{{\tau'}^a(y')}
    \big(\beta'_1,\ldots,\beta'_k\big)} \\
     &=& \tilde{q}^a_c\left(y',B_{\tilde{\Phi}^a(y')}\big(
              (T_{y'}\tilde{\Phi}^a)^{-1*}(T_{y'}{\tau'})^*
                \beta'_1,\ldots,
                (T_{y'}\tilde{\Phi}^a)^{-1*}(T_{y'}{\tau'})^*
                \beta'_k\big)\right) \\
     &=& \tilde{q}^a_c\left(y',B_{\tilde{\Phi}^a(y')}\big(
              (T_{\tilde{\Phi}^a(y')}{\tau'}^a)^*
              (T_{{\tau'}^a(y')}\Phi)^{-1*}
                \beta'_1,\ldots,
              (T_{\tilde{\Phi}^a(y')}{\tau'}^a)^*
              (T_{{\tau'}^a(y')}\Phi)^{-1*}
                \beta'_k\big)\right) \\
     &=&
      (\Phi^*\Psi B)_{{\tau'}^a(y')} \big(\beta'_1,\ldots,\beta'_k\big),
\eeas
alors, en particulier
\[
  \Psi'(\tilde{\Phi}^{a*}\tilde{A}[\nabla])=\Phi^*(\Psi\tilde{A}[\nabla])
            = \Phi^*A,
\]
et puisque $\tilde{\Phi}^{a*}\tilde{A}[\nabla]$ est un \'{e}l\'{e}ment de
$\Ginf(S^kT\tilde{N}^a)^{-c/a}_{\widetilde{\Phi^*\nabla}}$ on peut
passer \`{a} l'inverse de la restriction de $\Psi'$ \`{a} cet espace,
ce qui nous donne l'\'{e}quation \'{e}nonc\'{e}e (\ref{EqRelNatProjNaturalite}).
\ebew

\noindent On peut d\'{e}duire une formule explicite pour les sections $A_{k-l}$
apparaissant dans l'expression pour $B$ avec
$\widetilde{\mathsf{Div}}B=0$:
\bcor
 Avec les hypoth\`{e}ses du th\'{e}or\`{e}me pr\'{e}c\'{e}dent on d\'{e}finit le morphisme de
 fibr\'{e}s vectoriels
 $\mathsf{r}_k:S^jTM\otimes |\Lambda^mT^*M|^c\ra S^{j-2}TM\otimes
 |\Lambda^mT^*M|^c$ pour tout
 $0\leq j\leq k$:
 \beq
   (\mathsf{r}_kC_j)_x:=(\bar{m}+k+j-1)(k-j+1)2i(\mathsf{r})\big(C_j\big)
 \end{equation}
 quel que soit $C_j\in S^jT_xM\otimes |\Lambda^mT_x^*M|^c$. Alors, pour
 tout $1\leq l\leq k$ o\`{u} $\bar{m}:=m-c(m+1)$:\\
 \beq \lbl{EqExpliRelNatANabla}
  A_{k-l} ={\textstyle \frac{(-\bar{a})^l}{l!(2k+\bar{m}-1)\cdots
                                  (2k+\bar{m}-l)} }
     \mathsf{Deg}_{-l}
     \left(\sum_{j=0}^l(\mathsf{Div}+\mathsf{r}_k)^j\right)A_k,
 \end{equation}
 o\`{u} le symbole $\mathsf{Deg}_{-l}$ d\'{e}signe la projection sur l'op\'{e}rateur
 de degr\'{e} $-l$ dans la somme ($\mathsf{Div}$ est de degr\'{e} $-1$ et
 $\mathsf{r}_k$ est de degr\'{e} $-2$).
\ecor
\bbew
 On utilise r\'{e}currence sur $l$: le cas $l=1$ est \'{e}gal \`{a}
 (\ref{EqRecAkmu}). \\
 D'apr\`{e}s l'\'{e}quation (\ref{EqRecAkmlpu}) on obtient
 \beas
  A_{k-(l+1)} & = &
   {\textstyle \frac{(-\bar{a})^{l+1}}{(l+1)!(2k+\bar{m}-1)\cdots
                                  (2k+\bar{m}-(l+1))} }
       \left(\mathsf{Div}\mathsf{Deg}_{-l}
     \Big(\sum_{j=0}^l(\mathsf{Div}+\mathsf{r}_k)^j\Big)A_k\right. \\
              &   &
       ~~~~~~~~~~~~~~~~~~~~~~~~+\left.\mathsf{r}_k\mathsf{Deg}_{-(l-1)}
     \Big(\sum_{j=0}^{l-1}(\mathsf{Div}+\mathsf{r}_k)^j\Big)A_k\right)
 \eeas
 On a $\mathsf{Div}\mathsf{Deg}_{-l}=\mathsf{Deg}_{-(l+1)}\mathsf{Div}$
 et $\mathsf{r}_k\mathsf{Deg}_{-(l-1)}
 =\mathsf{Deg}_{-(l+1)}\mathsf{r}_k$; de plus
 \[
    \mathsf{Deg}_{-(l+1)}\Big(\mathsf{r}_k(\mathsf{Div}+\mathsf{r}_k)^l\Big)
          =0
 \]
 car $(\mathsf{Div}+\mathsf{r}_k)^l$ est une somme de mon\^{o}mes
 noncommutatifs \`{a} deux lettres $\mathsf{Div}$ et $\mathsf{r}_k$ dont
 la somme du nombre $d$ de lettres $\mathsf{Div}$ et du nombre
 $r$ de lettres $\mathsf{r}_k$ est \'{e}gale \`{a} $l$, donc le degr\'{e} d'un mon\^{o}me
 de $\mathsf{r}_k(\mathsf{Div}+\mathsf{r}_k)^l$ est \'{e}gal
 $-2-d-2r=-2-r-l$ ce qui est toujours strictement inf\'{e}rieur \`{a} $-(l+1)$.
 Un calcul simple montre le reste de la r\'{e}currence.
\ebew

\brem
 Le cas important $c=0$ satisfait toujours \`{a} l'hypoth\`{e}se (\ref{EqCondcmk})
 quel que soit l'entier $m\geq 2$ et $k$. En particulier, pour $k=1$, le
 rel\`{e}vement naturel d'un champ de vecteurs $X$ sur $M$ est \'{e}gal
 \[
     \tilde{X}[\nabla]=X^h-a(\mathsf{Div}A)\mathsf{E}.
 \]
 Ce champ de vecteurs ne d\'{e}pend pas du tout de la connexion $\nabla$:
 soit $F_t$ le flot de $X$, alors le champ de vecteurs correspondant
 au flot $\tilde{F}^a_t$ est exactement donn\'{e} par $\tilde{X}[\nabla]$.
\erem

\section{Existence d'une prescription d'ordre na\-tu\-relle projectivement
     invariante}

Soit $M$ une vari\'{e}t\'{e} de dimension $m$, $m\geq 2$, soient $a,b\in\real$
o\`{u} $a\neq 0$ et soit $k\in\nat$. Pour tout $y\in \tilde{M}^a$ on d\'{e}finit
d'abord un rel\`{e}vement $(~)^h_y:S^kT^*_{\tau^a(y)}M\otimes
|\Lambda^mT^*M|^b_{\tau^a(y)}$ dans $S^kT_y^*M$ par
\beq
  \big(s\otimes (\beta_1\vee\cdots\vee\beta_k)\big)^h_y :=
    (y^{-1}s)\big((T_y\tau^a)^*\beta_1 \vee\cdots\vee(T_y\tau^a)^*\beta_k)\big)
\end{equation}
Il en r\'{e}sulte le rel\`{e}vement $\gamma\mapsto \gamma^h$ des sections
appartenant \`{a} $\Ginf(S^kT^*M\otimes |\Lambda^mT^*M|^b)$ o\`{u}
$\gamma^h$ est un \'{e}l\'{e}ment de $\Ginf(S^kT^*\tilde{M}^a)$. Ce rel\`{e}vement
ne d\'{e}pend pas d'une connexion.

Soit $\Ginf(S^kT^*\tilde{M}^a)^{-b/a}$ le sous-espace des sections $\zeta$
dans $\Ginf(S^kT^*\tilde{M}^a)$ qui soient
{\em $(-\frac{b}{a})$-\'{e}quivariantes} dans le sens suivant:
\[
    G_t^*\zeta = e^{-(b/a)t}\zeta
\]
o\`{u} $G_t$ d\'{e}signe le flot du champ d'Euler $\mathsf{E}$. Soit $\nabla$ une
connexion sans torsion dans le fibr\'{e} tangent de $M$. Alors l'application
$\real$-lin\'{e}aire $\Psi_\nabla$ suivante de $\Ginf(S^kT^*\tilde{M}^a)^{-b/a}$ dans
$\Ginf(S^kT^*M\otimes |\Lambda^mT^*M|^b)$ est bien d\'{e}finie:
\beq
   (\Psi_\nabla\zeta)_{\tau^a(y)}
       \big(v_1,\ldots,v_k\big) :=
       \tilde{q}^a_b\big(y,\zeta_z({v_1^h}_y,\ldots,{v_k^h}_y)\big)
\end{equation}
quels que soient $v_1,\ldots,v_k\in T_{\tau^a(y)}M$, car le membre droit
ne d\'{e}pend pas de la fibre sur $\tau^a(y)$. Soit $\gamma_k:=\Psi_\nabla
\zeta$. Alors, $\Psi_\nabla\big(\zeta-\gamma_k^h\big)=0$, alors \`{a} l'aide
d'un raisonnement enti\`{e}rement analogue \`{a} celui donn\'{e} dans la premi\`{e}re
partie de la d\'{e}monstration du th\'{e}or\`{e}me \ref{TConstrTildeANabla} on montre
qu'il existent des sections $\gamma_j\in\Ginf(S^jT^*M\otimes
|\Lambda^mT^*M|^b)$ avec $0\leq j\leq k$, uniquement d\'{e}termin\'{e}es par
$\zeta$, telles que
\beq
  \zeta = \gamma_k^h + \omega^a\vee \gamma_{k-1}^h
   +\ldots+ {\omega^a}^{\vee (k-1)}\vee \gamma_1 +
            {\omega^a}^{\vee k}\vee \gamma_0.
\end{equation}
\blem\lbl{LNablaTildeDual}
 Soit $\gamma\in \Ginf(S^jT^*M\otimes
|\Lambda^mT^*M|^b)$, soit $X$ un champ de vecteurs sur $M$ et soit
 $\bar{b}:=(m+1)b$. Alors
 \bea
  \tilde{\nabla}_{X^h}\gamma^h & = & (\nabla_X\gamma)^h
  -  {\textstyle \frac{1}{\bar{a}} }
          \omega^a\vee \big(i(X)\gamma\big)^h \lbl{EqNablaTildeDualHH} \\
  \tilde{\nabla}_{\mathsf{E}}\gamma^h & = &
           -{\textstyle \frac{\bar{b} + j}{\bar{a}}}\gamma^h
                       \lbl{EqNablaTildeDualVH} \\
  \tilde{\nabla}_{X^h}\omega^a & = &
        -{\textstyle \frac{a}{2}}(\mathrm{tr}R)\big(X,~\big)^h
        +\bar{a}\mathsf{r}(X,~)^h \lbl{EqNablaTildeDualHV}\\
  \tilde{\nabla}_{\mathsf{E}}\omega^a & = &
       -  {\textstyle \frac{1}{\bar{a}} }\omega^a
                  \lbl{EqNablaTildeDualVV}
 \eea
\elem
\bbew
  Il suffit de v\'{e}rifier ces formules localement: on peut supposer que
  $\gamma$ est de la forme $\varphi\otimes \beta^{\vee k}$ o\`{u}
  $\varphi$ est une $b$-densit\'{e} locale et $\beta$ est une $1$-forme locale
  sur $M$. Soit $Y$ un autre champ de vecteurs. Alors
  \beas
   \big(\tilde{\nabla}_{X^h}(\varphi\otimes\beta)^h\big)(Y^h) & = &
                (\nabla_X\varphi)^h \beta^h(Y^h)
                   +\varphi^h\Big(X^h\big(\beta(Y)^h\big)
                             - \beta^h(\tilde{\nabla}_{X^h}Y^h)\Big) \\
                             & = &
                (\nabla_X\varphi\otimes\beta)^h(Y^h)
                   +\varphi^h\Big(X\big(\beta(Y)\big)
                             -\beta(\nabla_XY)\Big)^h \\
                             & = &
                \Big(\nabla_X\varphi\otimes\beta
                      + \varphi\otimes \nabla_X\beta\Big)^h(Y^h)
  \eeas
  et
  \beas
   \big(\tilde{\nabla}_{X^h}(\varphi\otimes\beta)^h\big)(\mathsf{E})
       & = &
       0-\varphi^h \beta^h(\tilde{\nabla}_{X^h}\mathsf{E}) \\
       & = &
       {\textstyle -\frac{1}{\bar{a}}}\varphi^h\beta(X)^h
  \eeas
  d'o\`{u} l'\'{e}quation (\ref{EqNablaTildeDualHH}) en utilisant la r\`{e}gle
  de Leibniz.
  Les autres \'{e}quations se calculent de mani\`{e}re similaire.
\ebew

\bprop
 Avec les notations pr\'{e}c\'{e}dentes, on a la formule suivante pour la
 diff\'{e}rentielle sym\'{e}trique sur $\tilde{M}^a$ par rapport \`{a} la connexion
 relev\'{e}e $\tilde{\nabla}$:
 \beq \lbl{EqRecurrenceTildeD}
    \tilde{\mathsf{D}}({\omega^a}^{\vee l}\vee \gamma^h)
    = {\omega^a}^{\vee l}\vee (\mathsf{D}\gamma)^h
      -2\bar{a}l{\omega^a}^{\vee (l-1)}\vee (\mathsf{r}\vee\gamma)^h
      -{\textstyle \frac{2j+l+\bar{b}}{\bar{a}} }
           {\omega^a}^{\vee (l+1)}\vee \gamma^h
 \end{equation}
\eprop
\bbew
 Il suffit de montrer cette \'{e}quation localement: soit $e_1,\ldots,e_m$ une
 base locale du fibr\'{e} tangent de $M$ et soit $e^1,\ldots,e^m$ sa base duale.
 On a la formule
 \beq
  \tilde{\mathsf{D}}({\omega^a}^{\vee l}\vee \gamma^h)
    = \omega^a\vee \big(\tilde{\nabla}_{\mathsf{E}}
            ({\omega^a}^{\vee l}\vee \gamma^h)\big)
            +\sum_{i=1}^m e^{i\hspace{1mm}h}\vee\big(
           \tilde{\nabla}_{e_i^h}({\omega^a}^{\vee l}\vee \gamma^h)\big),
 \end{equation}
 et en utilisant le lemme \ref{LNablaTildeDual} pr\'{e}c\'{e}dent, on arrive \`{a}
 la formule \'{e}nonc\'{e}e.
\ebew

\noindent Par une r\'{e}currence \'{e}vidente on arrive au corollaire suivant:
\bcor \lbl{CorTildeDketDk}
 Soit $\varphi\in \Ginf(|\Lambda^mT^*M|^b)$ et $k\in\nat$. Alors
 $\tilde{\mathsf{D}}^k\tilde{\varphi}$ est de la forme
 \beq
  \tilde{\mathsf{D}}^k\tilde{\varphi} =
    (\mathsf{D}^k\varphi)^h+\sum_{l=1}^kB_{k-l}\vee (\mathsf{D}^l\varphi)^h
 \end{equation}
 o\`{u} les sections $B_{k-l}\in\Ginf(S^{k-l}T^*\tilde{M}^a$ sont invariantes
 par le flot du champ d'Euler et sont toutes une somme de termes de la
 forme
 \beq
     \lambda_r\big(-{\textstyle \frac{1}{\bar{a}}}{\omega^a}\big)^{\vee r}
                 \vee\beta_{k-l-r}^h
 \end{equation}
 o\`{u} $r\in\nat$ avec $0\leq r\leq k-l$, le nombre r\'{e}el $\lambda_r$ ne d\'{e}pend
 pas de $\bar{a}$, et la
 section $\beta_{k-l-r}$ appartient \`{a} $\Ginf(S^{k-l-r}T^*M)$ est est un
 polyn\^{o}me \`{a} coefficients rationnels en diff\'{e}rentielles sym\'{e}triques du
 tenseur $\mathsf{r}$.
\ecor

\bsat \lbl{TExPreNatProjInv}
\ben
 \item
 Soit $m$ un entier positif, $m\geq 2$. Soient $b,c\in \real$ tels que
 \beq \lbl{EqCondcmsansk}
   c\not\in\left\{ {\textstyle \frac{j+m}{m+1}}~|~j\in\nat~\right\}.
 \end{equation}
 Alors il existe une prescription d'ordre naturelle $(\rho_L^k)_{k\in\nat}$
 projectivement
 invariante pour les foncteurs de fibr\'{e}s vectoriels
 $F=|\Lambda^mT^*|^b$ et $F'=|\Lambda^mT^*|^{b+c}$.\\
 On a la formule explicite suivante pour $\rho_L$: soit $M$ une vari\'{e}t\'{e}
 dif\-f\'{e}\-ren\-tiable de dimension $m$, soit $k\in \nat$,
 $a\in\real\setminus\{0\}$,
 $A\in\Ginf(S^kTM\otimes |\Lambda^mT^*M|^c)$, $\nabla$ une connexion
 sans torsion dans le fibr\'{e} tangent de $M$ et $\varphi$ une $b$-densit\'{e}
 sur $M$, alors pour tout $y$ dans le fibr\'{e} principal des $a$-densit\'{e}s sur
 $M$:
 \beq \lbl{EqDefRhoL}
   \big(\rho_L[\nabla](A)(\varphi)\big)_{\tau^a(y)}
     :=\tilde{q}^a_{b+c}\left(y,
       \big(\rho_s[\tilde{\nabla}](\tilde{A}[\nabla])(\tilde{\varphi})\big)_y
        \right)
 \end{equation}
 o\`{u} $\rho_s$ d\'{e}signe la prescription d'ordre standard (voir la d\'{e}finition
 \ref{DPreOrdStandard}). La prescription d'ordre $\rho_L$ ne d\'{e}pend pas de $a$.
 \item La conjecture de Lecomte est vraie.
 \een
\esat
\bbew
 (1) Soit $G_t$ le flot du champ d'Euler $\mathsf{E}$ sur $\tilde{M}^a$.
 La fonction $\tilde{\varphi}$ est $(-\frac{b}{a})$-\'{e}quivariante et
 la connexion relev\'{e}e $\tilde{\nabla}$ est invariante par rapport \`{a} $G_t$
 d'apr\`{e}s sa construction (th\'{e}or\`{e}me \ref{TExUnConnProjEq}). Il s'ensuit
 que les diff\'{e}rentielles sym\'{e}triques it\'{e}r\'{e}es
 $\tilde{\mathsf{D}}^k\tilde{\varphi}$ de $\tilde{\varphi}$
 par rapport \`{a} $\tilde{\nabla}$ forment un champ de tenseurs dans
 $\Ginf(\tilde{M}^a,S^kT^*\tilde{M}^a)$ qui est
 $(-\frac{b}{a})$-\'{e}quivariant.
 D'autre part, le rel\`{e}vement $\tilde{A}[\nabla]$ de $A$
 est un champ de tenseurs $(-\frac{c}{a})$-\'{e}quivariant
 d'apr\`{e}s sa construction dans le th\'{e}or\`{e}me \ref{TConstrTildeANabla}
 (qui est possible gr\^{a}ce \`{a} la condition (\ref{EqCondcmsansk})).
 Puisque l'accouplement naturel entre $\tilde{A}[\nabla]$ et
 $\tilde{\mathsf{D}}^k\tilde{\varphi}$ d\'{e}finit la prescription d'ordre
 standard $\rho_s[\tilde{\nabla}](\tilde{A}[\nabla])(\tilde{\varphi})$
 il vient que ce dernier est $(-\frac{b+c}{a})$-\'{e}quivariant.
 Il s'ensuit que la formule pour $\rho_L[\nabla](A)(\varphi)$ donne
 une $b+c$-densit\'{e} bien d\'{e}finie. De plus, la diff\'{e}rentielle sym\'{e}trique
 $\tilde{\mathsf{D}}^k\tilde{\varphi}$ contient une seule fois le terme
 $(\mathsf{D}^k\varphi)^h$ d'apr\`{e}s le corollaire \ref{CorTildeDketDk} et
 des diff\'{e}rentielles sym\'{e}triques $(\mathsf{D}^l\varphi)^h$ de degr\'{e} $l<k$.
 La section $\tilde{A}[\nabla]$ commence par le terme $A^h$ d'apr\`{e}s la formule
 (\ref{EqExpliRelNatANabla}). Il s'ensuit que $\rho_L[\nabla](A)$ est
 un op\'{e}rateur diff\'{e}rentiel de $\Ginf(|\Lambda^mT^*M|^b)$ dans
 $\Ginf(|\Lambda^mT^*M|^{b+c})$ dont le symbole principal est \'{e}gal $A$, donc
 $\rho_L$ d\'{e}finit une prescription d'ordre. De plus, le nombre r\'{e}el
 $\bar{a}=(m+1)a$ n'y appara\^{\i}t que dans les combinaisons
 $(-\omega^a/\bar{a})$
 dans $\tilde{\mathsf{D}}^k\tilde{\varphi}$ et $-\bar{a}\mathsf{E}$
 dans $\tilde{A}[\nabla]$. Puisque seul l'accouplement naturel entre
 $-\omega^a/\bar{a}$ et $-\bar{a}\mathsf{E}$ figure dans la formule
 (\ref{EqPreOrdStandard}) pour la prescription d'ordre standard, on
 voit que $\rho_L$ ne d\'{e}pend pas du nombre r\'{e}el $a$. \\
 Par construction, la connexion $\tilde{\nabla}$ et le rel\`{e}vement
 $\tilde{A}[\nabla]$ ne d\'{e}pendent que de la classe projective de la
 connexion $\nabla$. Il vient que $\rho_L$ est projectivement invariante.\\
 Soit $N$ une autre vari\'{e}t\'{e} de dimension $m$ et $\Phi:N\ra M$ une
 immersion. Soit $y'\in\tilde{N}^a$ et $\tilde{q}_{b+c}^{\prime a}:
 \tilde{N}^a\times\real
 \ra |\Lambda^mT^*N|^{b+c}$ la projection du fibr\'{e} associ\'{e}. Alors,
 gr\^{a}ce \`{a} la naturalit\'{e} des rel\`{e}vements $\nabla\mapsto\tilde{\nabla}$
 et $(A,\nabla)\mapsto \tilde{A}[\nabla]$ et \`{a} la naturalit\'{e} de
 la prescription d'ordre standard on a
 \beas
  \left( \rho_L[\Phi^*\nabla](\Phi^*A)(\Phi^*\varphi)\right)_{{\tau'}^a(y')}
     & = &
    \tilde{q}_{b+c}^{\prime a}\left(y',
    \Big(\rho_s[\widetilde{\Phi^*\nabla}]
       \big(\widetilde{\Phi^*A}[\Phi^*\nabla]\big)
         (\widetilde{\Phi^*\varphi})\Big)_{y'}
                         \right) \\
     & = &
    \tilde{q}_{b+c}^{\prime a}\left(y',
    \Big(\rho_s[\tilde{\Phi}^{a*}\widetilde{\nabla}]
    \big(\tilde{\Phi}^{a*}(\widetilde{A}[\nabla])\big)
         (\tilde{\Phi}^{a*}\widetilde{\varphi})\Big)_{y'}
                         \right) \\
     & = &
    \tilde{q}_{b+c}^{\prime a}\left(y',
    \big(\rho_s[\widetilde{\nabla}](\widetilde{A}[\nabla])
         (\widetilde{\varphi})\big)_{\tilde{\Phi}^a(y')}
                         \right) \\
     & = &
       \left(\Phi^*\big(\rho_L[\nabla](A)(\varphi)\big)
         \right)_{{\tau'}^a(y')}
 \eeas
 ce qui montre que $\rho_L$ est naturelle.

 (2) Ceci est le cas particulier $c=0$ de la partie (1).
\ebew

Dans le cas particulier o\`{u} $M$ est munie d'une connexion sans torsion
$\nabla$ dont la partie sym\'{e}trique du tenseur de Ricci s'annule (en particulier
$M=\real^m$) on retrouve la formule explicite dans \cite{LO99},
p.182, eqs (4.13), (4.14):
\bcor
 Soit $M$ une vari\'{e}t\'{e} diff\'{e}rentiable de dimension $m\geq 2$, soient
 $b,c\in\real$ tels que la condition (\ref{EqCondcmsansk})
 soit satisfaite et soit
 $\nabla$ une
 connexion sans torsion dans le fibr\'{e} tangent de $M$ telle que la partie
 sym\'{e}trique du tenseur de Ricci s'annule.\\
 Alors il y a la formule explicite suivante pour la prescription d'ordre
 $\rho_L[\nabla]$:
 \beq
   \rho_L[\nabla](A)\big(\varphi\big)=\sum_{l=0}^k{k\choose l}
                  {\textstyle
                  \frac{\big(k-1+(m+1)b\big)\cdots\big(k-l+(m+1)b\big) }
                       {\big(2k-1+m-(m+1)c\big)\cdots\big(2k-l+m-(m+1)c\big)} }
                       \rho_s\big(\mathsf{Div}^l(A)\big)\big(\varphi\big)
 \end{equation}
 quelles que soient les sections $A\in\Ginf(S^kTM\otimes
 |\Lambda^mT^*M|^c)$ et $\varphi$ appartenant \`{a} $\Ginf(|\Lambda^mT^*M|^b)$.
\ecor
\bbew
 D'apr\`{e}s (\ref{EqExpliRelNatANabla}) on trouve --gr\^{a}ce \`{a} $\mathsf{r}_k=0$:
 \[
   A_{k-l}={\textstyle
                  \frac{(-\bar{a})^l }
                       {l!\big(2k-1+m-(m+1)c\big)\cdots\big(2k-l+m-(m+1)c\big)} }
                       \mathsf{Div}^l(A)
 \]
 De plus, quand on fait l'{\em Ansatz} suivant pour
 $\tilde{\mathsf{D}}^k\tilde{\varphi}$:
 \[
   \tilde{\mathsf{D}}^k\tilde{\varphi}
       =\sum_{l=0}^k u^{(k)}_l
            \big({\textstyle \frac{-\omega^a}{\bar{a}}}\big)^{\vee l}
            \vee (\mathsf{D}^{k-l}\varphi)^h
 \]
 avec des nombres r\'{e}els $u^{(k)}_l$, alors l'\'{e}quation (\ref{EqRecurrenceTildeD})
 nous donne la relation suivante entre les coefficients $u^{(k)}_l$:
 \[
    u^{(k+1)}_l= u^{(k)}_l+(2k-l+1+(m+1)b)u^{(k)}_{l-1}
 \]
 dont la solution unique se trouve par une r\'{e}currence \'{e}l\'{e}mentaire:
 \[
    u^{(k)}_l={k\choose l}\big(k-1+(m+1)b\big)\cdots\big(k-l+(m+1)b\big),
 \]
 d'o\`{u} l'on obtient la formule \'{e}nonc\'{e}e en utilisant la formule
 (\ref{EqPreOrdStandard})
 pour la prescription d'ordre standard.
\ebew

\section{Prescriptions d'ordre \'{e}quivariantes}

Soit $M$ une vari\'{e}t\'{e} diff\'{e}rentiable de dimension $m\geq 2$, soit
$G$ un groupe de Lie et soit $\Phi:G\times M\ra M$ une action de $G$
\`{a} gauche sur $M$. On \'{e}crit $\Phi_g$ pour $\Phi(g,~)$ quel que soit
$g\in G$. Soit $a$ un nombre r\'{e}el non nul. Gr\^{a}ce \`{a} la fonctorialit\'{e}
de $\tilde{F}^a$, les diff\'{e}omorphismes $\tilde{\Phi}^a_g:=
\widetilde{\Phi_g}^a$ d\'{e}finissent une action \`{a} gauche `relev\'{e}e' de
$G$ sur le fibr\'{e} principal de $a$-densit\'{e}s $\tilde{M}^a$. De cette
mani\`{e}re, on obtient une action sur tous les fibr\'{e}s de densit\'{e}s.

On a le corollaire suivant du th\'{e}or\`{e}me \ref{TExPreNatProjInv}:

\bcor \lbl{CorEquivar}
Soit $b,c\in\real$ tels que la condition (\ref{EqCondcmsansk}) soit satisfaite.
Soit $\nabla$ une connexion sans torsion dans le fibr\'{e} tangent de $M$.
On suppose que l'action relev\'{e}e $\tilde{\Phi}^a$ de $G$ sur
$\tilde{M}^a$ laisse invariante la connexion relev\'{e}e $\tilde{\nabla}$,
voir le th\'{e}or\`{e}me \ref{TExUnConnProjEq}, c.-\`{a}-d.
\beq
    \tilde{\Phi}^{a*}_g\tilde{\nabla}=\tilde{\nabla}~~~\forall g\in G.
\end{equation}
Alors la prescription d'ordre $\rho_L[\nabla]$ d\'{e}finie dans le th\'{e}or\`{e}me
\ref{TExPreNatProjInv} est $G$-\'{e}qui\-va\-ri\-ante.
\ecor
\bbew
 La $G$-\'{e}quivariance \'{e}nonc\'{e}e est une cons\'{e}quence imm\'{e}diate de la
 d\'{e}finition (\ref{EqDefRhoL}) et de la
 naturalit\'{e} (\ref{EqDefNaturalitePreOrd}) de $\rho_L$ et du fait que
 \[
  \widetilde{\Phi_g^*\nabla}=\tilde{\Phi}^{a*}_g\tilde{\nabla}
                            = \tilde{\nabla}.
 \]
\ebew

\noindent On voit qu'il n'est pas n\'{e}cessaire que $\nabla$ soit invariante
par l'action de $G$ sur $M$.

Le cas particulier $M=S^m$ avec $m\geq 2$ (voir l'exemple \ref{ExemSmI})
est int\'{e}ressant:
\bexem
 Soit $M=S^m$ et $\tilde{M}^a=\real^{m+1}\setminus\{0\}$ (o\`{u}
 $a=\frac{1}{m+1}$). Soit $\Phi:S^m\ra S^m$ un diff\'{e}omorphisme. En
 utilisant la section $\sigma=|\mathsf{g}|^{a/2}$
 (voir (\ref{EqDefDetmathsfgademi})) on arrive \`{a} la formule
 \beq
   \tilde{\Phi}^a(y)=\frac{|y|}{f_\Phi(y/|y|)}\Phi(y/|y|)~~~
   \forall y\in \real^{m+1}\setminus\{0\}
 \end{equation}
 avec la valeur absolue du jacobien de $\Phi$ \`{a} la puissance $a$,
 \beq
     f_{\Phi}(x):= |\det(\Phi(x),T_x\Phi e_1,\ldots, T_x\Phi e_m)|^a~~~
   \forall x\in S^m
 \end{equation}
 o\`{u} $e_1,\ldots,e_m$ est une base orthonormale de l'espace tangent
 $T_xS^m$. En particulier, soit $g\in GL(m+1,\real)$ et soit
 $\Phi_g:S^m\ra S^m$ l'homographie $\Phi_g(x):=\frac{1}{|gx|}gx$.
 Evidemment, l'application $g\mapsto \Phi_g$ d\'{e}finit une action \`{a} gauche
 de $GL(m+1,\real)$ sur $S^m$ dont le noyau est \'{e}gal \`{a} l'ensemble de tous
 les multiples strictement positifs de l'application identique. Un calcul
 \'{e}l\'{e}mentaire montre que
 \beq
    \tilde{\Phi}^a_g(y) = |\det g|^{-\frac{1}{m+1}}gy.
 \end{equation}
\eexem
Ceci nous donne le
\bcor
 Soit $b,c\in\real$ tels que la condition (\ref{EqCondcmsansk}) soit satisfaite.
 \ben
  \item
   Il existe une prescription d'ordre $\rho$ sur $S^m$ qui soit \'{e}quivariante
 par l'action du groupe $GL(m+1,\real)/\real^+$ de toutes les homographies.
  \item
   Il existe une prescription d'ordre $\rho$ sur $\real^m$ qui soit
   \'{e}quivariante par l'action de l'alg\`{e}bre de Lie
   $\mathfrak{sl}(m+1,\real)$ de toutes les homographies infinit\'{e}simales.
 \een
\ecor
\bbew
 (1) Puisque les diff\'{e}omorphismes $\tilde{\Phi}^a_g$ sont \'{e}videmment des
 applications lin\'{e}aires de $\real^{m+1}$, ils pr\'{e}servent donc la connexion
 canonique $\tilde{\nabla}$ dans $\real^{m+1}$. Alors on peut utiliser
 le corollaire \ref{CorEquivar} pour voir que la prescription d'ordre
 $\rho_L$ est $GL(m+1,\real)/\real^+$-\`{e}quivariante.\\
 (2) De l'\'{e}nonc\'{e} pr\'{e}c\'{e}dent on d\'{e}duit l'\'{e}quivariance infinit\'{e}simale de $\rho_L$
 par rapport \`{a} l'alg\`{e}bre de Lie $\mathfrak{sl}(m+1,\real)$. Puisque
 l'\'{e}quivariance
 infinit\'{e}simale peut toujours \^{e}tre localis\'{e}e \`{a} n'importe quel ouvert de $S^m$,
 alors on obtient le r\'{e}sultat en choisissant l'ouvert
 $U:=S^m\setminus \{y_0\}\cong \real^m$ o\`{u} $y_0$ est un vecteur arbitraire
 de $S^m$.
\ebew

\section*{Probl\`{e}mes ouverts}
 \addcontentsline{toc}{section}{Probl\`{e}mes ouverts}

Les probl\`{e}mes suivants me semblent int\'{e}ressants:

\ben
 \item Puisque la prescription d'ordre $\rho_L[\nabla]$ sur $\real^m$
 est \'{e}quivariante par $\mathfrak{sl}(m+1,\real)$ et une telle prescription
 d'ordre est unique (voir par exemple \cite{LO99}), on peut imaginer
 que $\rho_L$ est l'unique
 prescription d'ordre naturelle projectivement \'{e}quivariante entre deux
 fibr\'{e}s de densit\'{e}s. Pour le
 montrer il faudra \'{e}tudier les op\'{e}rateurs naturels entre les
 foncteurs de fibr\'{e}s $Q_\tau P^1M\times S^kTM$ et $S^lTM$ pour $l\leq k-1$.
 Vraisemblablement ces op\'{e}rateurs se r\'{e}duisent aux contractions des d\'{e}riv\'{e}es
 covariantes avec des polyn\^{o}mes des d\'{e}riv\'{e}es covariantes du tenseur de
 courbure; et on peut esp\'{e}rer qu'ils s'annulent quand on impose
 qu'ils soient les m\^{e}mes pour une connexion projectivement \'{e}quivalente.

 \item Que se passe-t-il avec les valeurs `r\'{e}sonnantes'
  \[
   c\in \left\{ {\textstyle \frac{j+m}{m+1}}~|~j\in\nat~\right\}?
  \]
  pour lesquelles le rel\`{e}vement $A\mapsto \tilde{A}[\nabla]$ ne fonctionne
  plus?

 \item Les amateurs des formules explicites peuvent s'amuser \`{a} trouver
 une formule explicite pour $\rho_L[\nabla]$ au cas o\`{u} la partie sym\'{e}trique
 du tenseur de Ricci ne s'annule pas, mais est constante: $\mathsf{Dr}=0$.
 Ceci est le cas pour tout espace affine localement sym\'{e}trique
 (\cite{KN69}, Ch. XI, p.222),
 ou pour une vari\'{e}t\'{e} einsteinienne (une vari\'{e}t\'{e} semi-riemannienne dont
 le tenseur de Ricci est un multiple constant de la m\'{e}trique, voir par
 exemple \cite{KN63}, p.293). Dans ce cas, les op\'{e}rations $\mathsf{Div}$
 et $i(\mathsf{r})$ commutent.

 \item La m\'{e}thode de relever les structures de $M$ \`{a} $\tilde{M}^a$
  se g\'{e}n\'{e}ralise-t-elle \`{a} d'autres fibr\'{e}s naturels? Par exemple,
  aux foncteurs de
  fibr\'{e}s $F=\Lambda T^*$, $F'=\Lambda^0 T^*$ ce qui nous donne des
  prescriptions d'ordre projectivement invariantes pour les op\'{e}rateurs
  diff\'{e}rentiels des formes diff\'{e}rentielles dans les fonctions.
  L'\'{e}quivariance dans le cas $\real^m$ a \'{e}t\'{e} \'{e}tudi\'{e} dans
  \cite{BHMP02} et \cite{Pon02}.

 \item Peut-on r\'{e}soudre le probl\`{e}me analogue des structures conformes
  (voir \cite{DLO99}) \`{a} l'aide d'un rel\`{e}vement \`{a} une vari\'{e}t\'{e} plus grande?
\een

\begin{appendix}

\section{Fibr\'{e}s principaux} \lbl{SecApp1}
{\small
Soit $m$ un entier positif et ${\mathcal{M}f}_m$ la cat\'{e}gorie des vari\'{e}t\'{e}s
diff\'{e}rentiables de
dimension $m$ dont les morphismes sont des immersions (voir \cite{KMS93},
p. 56). Soit $\mathcal{FM}$ la cat\'{e}gorie des vari\'{e}t\'{e}s fibr\'{e}es
(voir \cite{KMS93}, p. 15) dont les objets sont des triplets $(M,\pi,N)$
o\`{u} $\pi:M\ra N$ est une submersion surjective entre deux vari\'{e}t\'{e}s
diff\'{e}rentiables $M$ et $N$, et dont les morphismes
$(\Phi,\underline{\Phi}):(M,\pi,N)\ra(\bar{M},\bar{\pi},\bar{N})$ sont des
applications $\Phi:M\ra \bar{M}$
de classe $\Cinf$ qui pr\'{e}servent les fibres, c.-\`{a}-d. qui induisent
des applications $\underline{\Phi}:N\ra \bar{N}$ telles que
$\underline{\Phi}\circ \pi = \bar{\pi} \circ \Phi$. Le foncteur de base
$B:\mathcal{FM}\ra{\mathcal{M}f}$ associe \`{a} $(M,\pi,N)$ la vari\'{e}t\'{e} de base
$N$ et \`{a} une application de classe $\Cinf$ pr\'{e}servant les fibres
l'application induite $\underline{\Phi}$. Pour tout $x\in N$ on appelle
{\em la fibre $M_x$} la sous-vari\'{e}t\'{e} $\pi^{-1}(x)$ de $M$. Dans toute vari\'{e}t\'{e}
fibr\'{e}e $(M,\pi,N)$
il existe un sous-fibr\'{e} int\'{e}grable naturel de $TM$, \`{a} savoir le {\em sous-fibr\'{e}
vertical $VM:=\mathrm{Ker}\hspace{1mm}T\pi$} qui est la r\'{e}union de tous
les {\em sous-espaces verticaux} $V_xM:=\mathrm{Ker}\hspace{1mm}T_x\pi$,
$x\in M$. Les champs de vecteurs qui
ont leurs valeurs dans $VM$ sont dits {\em champs de vecteurs verticaux}.
Soit $U$ un ouvert de la base $N$, soit $M|_U:=\pi^{-1}(U)$ et
$p_U$ la restriction de la projection $\pi$ \`{a} $M|_U$. Alors
$(M|_U,\pi|_U,U)$ est une vari\'{e}t\'{e} fibr\'{e}e.
Une {\em section locale de classe
$\Cinf$ de $\pi_U:M|_U \ra U$} est une application
$\varphi:U\ra M$ de classe $\Cinf$ telle que
$\pi\big(\varphi(u)\big)=u$ quel que soit $u\in U$. Si $U=N$ on enl\`{e}ve
l'adjectif `local'. On \'{e}crit $\Ginf(U,M|_U)$ pour l'ensemble de toutes les
sections locales de classe $\Cinf$ de $\pi|_U:M|_U\ra U$.

Un {\em fibr\'{e} principal} $(P,\pi,M,G)$ sur la {\em vari\'{e}t\'{e} de base} $M$
\`{a} {\em groupe structural} $G$ est une vari\'{e}t\'{e} fibr\'{e}e $(P,\pi,M)$ telle
que $\pi:P\ra M$ soit
localement trivial sur $M$ et {\em l'espace total $P$} est muni d'une action
\`{a} droite $r:P\times G\ra P$ d'un groupe de Lie $G$ (\`{a} l'alg\`{e}bre de Lie
$\mathfrak{g}$) libre et propre telle que $M$ soit l'espace quotient de
cette action. On \'{e}crit $pg$ ou $r_g(p)$ pour $r(p,g)$ o\`{u} $p\in P,g\in G$. Le
{\em champ fondamental} $\xi^*$ associ\'{e} \`{a} $\xi\in\mathfrak{g}$
est le g\'{e}n\'{e}rateur infinit\'{e}simal de l'action \`{a} droite, c.-\`{a}-d.
$\xi^*(p):=\frac{d}{dt}\big(p \exp(t\xi)\big)|_{t=0}$ pour tout $p\in P$.
On a la formule $[\xi^*,\eta^*]=[\xi,\eta]^*$ $\forall
\xi,\eta\in\mathfrak{g}$.
Le sous-espace vectoriel $V_p:=\{\xi^*_p~|~\xi\in\mathfrak{g}\}$ de l'espace
tangent co\"{\i}ncide avec le sous-espace {\em vertical} en $p$ quel que
soit $p\in P$.

Un morphisme de fibr\'{e}s principaux $(P,\pi,M,G)\ra (P',\pi',M',G')$ est un
triplet $(\Phi,\underline{\Phi},\phi)$ d'applications de classe $\Cinf$
o\`{u} $\Phi:P\ra P'$, $\underline{\Phi}:M\ra M'$ et $\phi:G\ra G'$ telles
que $(\Phi,\underline{\Phi})$ soit un morphisme de vari\'{e}t\'{e}s fibr\'{e}es,
$\Phi(pg)=\Phi(p)\big((\phi(g)\big)~\forall p\in P,g\in G$ et $\phi$ est
un homomorphisme de groupes de Lie.

Si $S$ est une autre vari\'{e}t\'{e} et $\ell:G\times S\ra S$
est une action gauche de $G$ sur $S$ (o\`{u} l'on \'{e}crit parfois $gs$ ou $\ell_g(s)$
pour $\ell(g,s)$, $\forall g\in G,s\in S$), alors le {\em fibr\'{e} associ\'{e} $P\times_G
S$ \`{a} fibre type $S$}
est d\'{e}fini par l'espace quotient de l'action \`{a} droite
$P\times S\times G\ra P\times S$ donn\'{e}e par
$\big(p,s\big)g:=\big(pg,\ell(g^{-1})s\big)$. On
note la projection $P\times S\ra P\times_G S$ par $q$ (et parfois par $q_M$).
La projection
$\tau:P\times_G S\ra M$ est donn\'{e}e par $\tau\big(q(p,s)\big):=\pi(p)$.
Au cas o\`{u} $S$ est un espace vectoriel $V$ de dimension finie sur $\real$ ou
$\complex$ et $\ell$ une repr\'{e}sentation de $G$, le fibr\'{e} associ\'{e} $P\times_G V$
est un fibr\'{e} vectoriel sur $M$. Soit $U\subset M$ un ouvert.  A chaque
section locale $\varphi$ dans $\Ginf(U,(P\times_G S)|_U)$ on peut associer une unique
application de classe
$\Cinf$, $f_\varphi:P|_U\ra S$ qui est $G$ \'{e}quivariante,
c.-\`{a}-d. $f_\varphi(pg)=\ell(g^{-1})f_\varphi(p)$ $\forall p\in P|_U,g\in G$,
et telle que $\varphi\big(\pi(p)\big)=q\big(p,f_\varphi(p)\big)$.

Une {\em $1$-forme de connexion} $\omega$ dans un fibr\'{e} principal est une
$1$-forme \`{a}
valeurs dans $\mathfrak{g}$ sur $P$ telle que ($i$) $\omega(\xi^*)=\xi$
quel que soit $\xi\in\mathfrak{g}$, ($ii$) $\omega_{pg}(T_pr_g~v)=
\mathrm{Ad}(g^{-1})\big(\omega_p(v)\big)$ quels que soient $p\in P,v\in T_pP$
et $g\in G$. Le sous-espace $\mathrm{Ker}\hspace{1mm}\omega_p=:H_pP$ de l'espace
tangent $T_pP$ est isomorphe --via $T_p\pi$-- \`{a} l'espace tangent $T_{\pi(p)}M$ et
s'appelle le {\em sous-espace horizontal de $T_pP$}. La r\'{e}union de tous
les sous-espaces horizontaux est dite le {\em sous-fibr\'{e} horizontal $HP$ de $TP$}.
Le fibr\'{e} tangent de $P$ est toujours la somme directe du fibr\'{e} vertical
et du fibr\'{e} horizontal.
Un vecteur tangent $w\in T_{\pi(p)}M$ correspond \`{a} un unique vecteur
tangent horizontal $w^h(p)\in H_pP$: $w^h(p)$ est dit le {\em rel\`{e}vement
horizontal de $w$ \`{a} $p$}. Le {\em rel\`{e}vement horizontal d'un champ de vecteurs
$X\in\Ginf(TM)$}, $X^h$, est d\'{e}fini par $X^h_p:=(X_{\pi(p)})^h(p)$.
On a $T\pi~X^h=X\circ\pi$. De plus, $X^h$ est $G$-invariant,
c.-\`{a}-d. $X^H(pg)=T_pr_g~X^h(p)$. Pour deux champs de vecteurs $X,Y$
on a la formule
\beq \lbl{EqCrLieHorPrincipal}
     [X^h,Y^h]=[X,Y]^h-\Omega(X^h,Y^h)^*
\end{equation}
o\`{u} $\Omega(Z_1,Z_2):=d\omega(Z_1,Z_2)+[\omega(Z_1),\omega(Z_2)]$ est la
{\em $2$-forme de courbure de la connexion $\omega$} quels que soient les
champs de vecteurs $Z_1,Z_2$ sur $P$.

Une $1$-forme de connexion peut \^{e}tre d\'{e}termin\'{e}e par une section de classe
$\Cinf$ d'un fibr\'{e} $QP$ localement trivial sur $M$, dit le {\em fibr\'{e} des
connexions principales}, dont la fibre type est un
espace affine model\'{e} sur l'espace vectoriel $\mathfrak{g}\otimes {\real^m}^*$:
dans le fibr\'{e} vectoriel $Hom(TP,\mathfrak{g})$ sur $P$
on consid\`{e}re le sous-fibr\'{e} affine $\check{\tau}:\check{Q}P\ra P$ d\'{e}fini par tous
les \'{e}l\'{e}ments $\omega_p$ de
$Hom(T_pP,\mathfrak{g})$ pour lesquels $\omega_p(\xi^*_p)=\xi$ quels que
soient $p\in P$ et $\xi\in\mathfrak{g}$. On voit que $\check{Q}P$ est
invariant par l'action \`{a} droite de $G$ sur $Hom(TP,\mathfrak{g})$ d\'{e}finie
par $\big(\check{r}_g(\omega_p)\big)(v):=\mathrm{Ad}(g^{-1})
\big(\omega_p(T_{pg}r_{g^{-1}}\hspace{1mm}v\big)$ quels que soient $v\in
T_{ph}P$ et $g\in G$. Puisqu'on a $\check{\tau}\circ \check{r}_g =
r_g\circ \check{\tau}$ l'espace quotient $\check{Q}P/G:=QP$ est un fibr\'{e}
localement trivial sur $M=P/G$. Puisque toute $1$-forme de connexion est
une section $G$-invariante de $\check{Q}P$ elle induit et est d\'{e}termin\'{e}e
par une section de $\tau:QP\ra M$, voir \'{e}galement \cite{KMS93}, p.159,
17.4.

Soit $(\Phi,\underline{\Phi},\phi):
(P,\pi,M,G)\ra (P',\pi',M',G')$ un morphisme de fibr\'{e}s principaux. Il y
maintenant
deux situations importantes pour `promouvoir' et `retirer' des $1$-formes de
connexion:\\
1. Au cas o\`{u} $\underline{\Phi}:M\ra M'$ est
un diff\'{e}omorphisme et $\omega$ est une $1$-forme de connexion sur $P$ il
existe une unique $1$-forme de connexion $\omega'$ sur $P'$ telle que
$\mathrm{Ker}\hspace{1mm}\omega'_{\Phi(p)}=
T_p\Phi\big(\mathrm{Ker}\hspace{1mm}\omega_{p}\big)$ $\forall p\in P$.
Dans ce cas, $\omega'_{\Phi(p)}\big(T_p\Phi~v\big)
=T_e\phi\big(\omega_p(v)\big)$ $\forall p\in P,v\in T_pP$, voir \cite{KN63},
p. 79, Prop. 6.1.\\
2. Au cas o\`{u} $T_e\phi:\mathfrak{g}\ra\mathfrak{g}'$ est un isomorphisme
d'alg\`{e}bres de Lie et $\omega'$ est une $1$-forme de connexion sur $P'$
il existe une unique $1$-forme de connexion $\omega=:\Phi^*\omega'$ sur $P$;
dite la {\em $1$-forme de connexion retir\'{e}e}, telle
que $\mathrm{Ker}\hspace{1mm}\omega'_{\Phi(p)}=
T_p\Phi\big(\mathrm{Ker}\hspace{1mm}\omega_{p}\big)$ $\forall p\in P$
d\'{e}finie par $\omega_p(v)
:=(T_e\phi)^{-1}\big(\omega'_{\Phi(p)}\big(T_p\Phi~v\big)\big)$,
voir \cite{KN63}, p. 81, Prop. 6.2.

Soit $\ell:G\times S\ra S$ une action \`{a} gauche de $G$ sur une vari\'{e}t\'{e}
diff\'{e}rentiable $S$ et $\omega$ une $1$-forme de connexion sur le fibr\'{e}
principal $(P,\pi,M,G)$. Pour $p\in P$ soit $H_pP$ le sous-espace horizontal
de $T_pP$. Pour tout $s\in S$ on d\'{e}finit le {\em sous-espace horizontal
$H_{q(p,s)}E$ de $T_{q(p,s)}E$} par
$T_{(p,s)}q\hspace{1mm}(H_pP\times\{0\}_s)$. La r\'{e}union de tous ces
sous-espaces horizontaux d\'{e}finit le {\em sous-fibr\'{e} horizontal $HE$ du fibr\'{e}
tangent de $E$}. De fa\c{c}on analogue on montre que $TE=VE\oplus HE$ et
qu'il existe un {\em rel\`{e}vement horizontal} $(~)^{\mathbf{h}}_e: T_xM\ra HE_e$ o\`{u}
$\tau(e)=x$ et le rel\`{e}vement horizontal $X^{\mathbf{h}}\in \Ginf(E,HE)$ d'un
champ de vecteurs $X\in\Ginf(TM)$.

Soit $S=V$ un espace vectoriel de
dimension finie sur $\real$ ou $\complex$, $\ell:G\times V\ra V$ une
r\'{e}pr\'{e}sentation de classe $\Cinf$ de $G$ et $\omega$ une $1$-forme de
connexion sur le fibr\'{e} principal $(P,\pi,M,G)$. Pour une section $\varphi$
du fibr\'{e} vectoriel
$E:=P\times_G V$ et un champ de vecteurs $X$ sur $M$ on d\'{e}finit la {\em d\'{e}riv\'{e}e
covariante $\nabla_X\varphi$} par $f_{\nabla_X\varphi}:=X^h(f_\varphi)$.
On a la formule pour le {\em tenseur de courbure}
$R(X,Y)\phi:=\nabla_X\nabla_Y\phi-\nabla_Y\nabla_X\phi-\nabla_{[X,Y]}\phi$:
\beq
    f_{R(X,Y)\phi}
      =\dot{\ell}\big(\Omega(X^h,Y^h)\big)(f_\phi)
\end{equation}
o\`{u} $\dot{\ell}:\mathfrak{g}\ra Hom(V,V)$ est la repr\'{e}sentation de l'alg\`{e}bre
de Lie $\mathfrak{g}$ d\'{e}finie par $\dot{\ell}(\xi)v=
\frac{d}{dt}\big(l(exp(t\xi))v\big)|_{t=0}$. Le tenseur de courbure $R$ est
une section dans $\Ginf\big(M,\Lambda^2T^*M\otimes Hom(E,E)\big)$ \\
Soient $e,e'\in E$ tels que $\tau(e)=x=\tau(e')$. On d\'{e}finit le {\em rel\`{e}vement
vertical de $e'$ \`{a} $e$}, not\'{e} ${e'}^v_e$ comme le vecteur tangent vertical
en $e$ d\'{e}fini par ${e'}^v_e:=\frac{d}{dt}(e+te')|_{t=0}$. Il existe
toujours un champ de vecteurs canonique sur $E$, \`{a} savoir {\em le champ
d'Euler}
$\mathsf{E}(e):=e^v_e:=\frac{d}{dt}(e\exp t)|_{t=0}$.
De la m\^{e}me mani\`{e}re on d\'{e}finit le {\em rel\`{e}vement vertical d'une section} de
classe $\Cinf$ $\varphi$ de $E$ par
$\varphi^v(e):=\varphi\big(\tau(e)\big)^v(e)$.
Le rel\`{e}vement horizontal $X^h$ de $X$ \`{a} $P$, vu comme un champ
de vecteurs sur $P\times G$ via $X^h_{(p,g)}=(X^h_p,0_g)$, est $q$-li\'{e} avec
le rel\`{e}vement horizontal $X^{\mathbf{h}}$ de $X$ \`{a} $E$, c.-\`{a}-d.
$T_{(p,g)}q\hspace{1mm}X^H(p,g)=X^{\mathbf{h}}\big(q(p,g)\big)$ quels que soient
$p\in P,g\in G$. On en
d\'{e}duit
la formule suivante pour le crochet de Lie de deux rel\`{e}vements
horizontaux, qui est analogue \`{a} (\ref{EqCrLieHorPrincipal}):
\beq \lbl{EqCrLieHorHor}
     [X^{\mathbf{h}},Y^{\mathbf{h}}]_e
       =[X,Y]^{\mathbf{h}}_e-\big(R(X,Y)e\big)^v_e ~~\forall e\in E.
\end{equation}
Soient $\phi,\phi'$ deux sections de classe $\Cinf$ du fibr\'{e} vectoriel
$E$. On d\'{e}duit les formules suivantes pour les crochets de Lie:
\bea
 ~   [X^{\mathbf{h}},\varphi^v] & = & (\nabla_X\varphi)^v. \lbl{EqCrLieHorVer} \\
  ~  [\varphi^v,{\varphi'}^v]   & = & 0. \lbl{EqCrLieVerVer}
\eea

Soit $M$ une vari\'{e}t\'{e} de dimension $m$ et soit $r$ un entier positif.
Soit $P^r M$ le fibr\'{e} de tous les jets
d'ordre $r$ en $0$ des applications $f:\real^m\ra M$ telles que $T_0f$ est
inversible, voir \cite{KMS93}, Chapitre IV, p. 116. $P^rM$ est un fibr\'{e}
localement trivial sur $M$ o\`{u} la projection
$\pi^r_0:P^rM\ra M$ est donn\'{e}e par $\pi^r_0\big(j^r_0(f)\big):=f(0)$. Soit $G^r_m$
d\'{e}finie par la fibre $({\pi^r_0})^{-1}(0)$ de $P^r\real^m$. Alors $G^r_m$ est un
groupe de Lie dont la multiplication se d\'{e}crit par la composition des applications
polynomiales $\real^m\ra\real^m$ pr\'{e}servant l'origine \`{a} degr\'{e} $r+1$ pr\`{e}s.
La composition $P^rM\times G^r_m\ra P^rM:\big(j^r_0(f),j^r_0(g)\big)\ra
j^r_0(f\circ g)$ est une action \`{a} droite de $G^r_m$ qui donne \`{a} $P^rM$ la structure
d'un fibr\'{e}
principal sur $M$ \`{a} groupe structural $G^r_m$ qui est dit le {\em fibr\'{e}
des rep\`{e}res d'ordre $r$ de $M$}. Au cas o\`{u} $r\geq r'$ la projection canonique
$j^r_0(f)\mapsto j^{r'}_0(f)$ des jets induit un morphisme de fibr\'{e}s
principaux $\pi^r_{r'}:P^rM\ra P^{r'}M$ avec
$\underline{\pi^r_{r'}}:=\id_M$ et $\phi:G^r_m\ra G^{r'}_m$ \'{e}tant la
projection des jets. Il est clair que $P^0M\cong M$ o\`{u} $G^0_m=\{\id\}$.

Pour $r=1$, le jet
$p:=j^1_0(f)$ s'identifie de fa\c{c}on canonique \`{a} la base
$P:=(e_1,\ldots,e_m)=(T_0f\hspace{1mm}a_1,\ldots,T_0f\hspace{1mm}a_m)$
de l'espace tangent $T_xM$ o\`{u} $x=\pi^1_0(p)=f(0)$ et
$(a_1,\ldots,a_m)$ est la base canonique de $\real^m$. Le groupe de Lie
$G^1_m$ est isomorphe \`{a} $GL(m,\real)$ o\`{u}
la multiplication \`{a} droite par $g=(g^i_j)_{1\leq i,j\leq m}\in GL(m,\real)$ est
d\'{e}finie par
$(e_1,\ldots,e_m)g:= (\sum_{i=1}^m g^i_1 e_i,\ldots,\sum_{i=1}^m g^i_m e_i)$.
Parfois, le fibr\'{e} principal $\big(P^1M,$ $\pi^1_0,M,GL(m,\real)\big)$ est dit le
{\em fibr\'{e} des rep\`{e}res lin\'{e}aires de $M$}. On note la repr\'{e}sentation canonique
de $GL(m,\real)\times \real^m\ra \real^m$ par
$(g,v)\mapsto gv$. Soit $(g,\gamma)\mapsto g\gamma:=\gamma\circ g^{-1}$
la repr\'{e}sentation contragr\'{e}diente de $GL(m,\real^m)$ sur l'espace dual
$\real^{m*}$ de $\real^m$. Pour deux entiers positifs $k,l$ soit
$V^{(k,l)}:={\real^m}^{\otimes k}\otimes {\real^{m*}}^{\otimes l}$, et
on a la repr\'{e}sentation de $GL(m,\real)$ sur $V^{(k,l)}$ par
$g(v_1\otimes\cdots \otimes v_k\otimes\gamma_1\otimes\cdots\otimes
\gamma_l):=(gv_1\otimes\cdots \otimes gv_k\otimes g\gamma_1\otimes\cdots\otimes
g\gamma_l)$ o\`{u} $v_1,\ldots,v_k\in\real^m$ et $\gamma_1,\ldots,\gamma_l\in
\real^{m*}$.
Evidemment, ${\real^m}^{\otimes 0}:=\real={\real^{*m}}^{\otimes 0}$.
Alors le fibr\'{e} associ\'{e} $P^1M\times_{GL(m,\real)}V^{(k,l)}$
est isomorphe au fibr\'{e} vectoriel $TM^{\otimes k}\otimes
T^*M^{\otimes l}$. En particulier, $P^1M\times_{GL(m,\real)}\real^m\cong
TM$. De fa\c{c}on analogue on voit que les fibr\'{e}s $S^kTM$, $S^lT^*M$,
$\Lambda^kTM$ et $\Lambda^lT^*M$ s'obtiennent en tant que fibr\'{e}s associ\'{e}s
de $P^1M$. On mentionne la {\em $1$-forme canonique $\theta_M$ \`{a} valeurs dans
$\real^m$} sur $P^1M$ d\'{e}finie par $\theta_M(\xi^*)=0$ $\forall
\xi\in\mathfrak{gl}(m,\real)$ et
$q\big(p,\theta_{Mp}(v)\big):=T_p\pi^1_0v$ $\forall p\in P^1M,v\in
T_pP^1M$.

Soit $\nabla$ une connexion affine dans le fibr\'{e} tangent. On peut lui associer
une unique $1$-forme de connexion $\omega$ d\'{e}finie sur le fibr\'{e} des
rep\`{e}res lin\'{e}aires $P^1M$. R\'{e}ciproquement, \`{a} chaque $1$-forme de connexion
il correspond une connexion affine $\nabla$ comme d\'{e}crit ci-dessus pour le
cas d'un fibr\'{e} principal g\'{e}n\'{e}ral. On \'{e}crit $Q_\tau P^1M$ pour le
sous-fibr\'{e} affine du fibr\'{e} $QP^1M$ des connexions principales de $P^1M$
dont les sections sont des $1$-formes de connexions correspondant aux
connexion sans torsion dans $TM$. \\
Soit $\omega^{(r)}$ une $1$-forme de connexion dans le fibr\'{e} $P^rM$. Le
morphisme de fibr\'{e}s principaux $\pi^r_1:P^rM\ra P^1M$ induit l'identit\'{e}
sur la base $M$, alors on a la $1$-forme de connexion $\omega$ sur $P^1M$
qui correspond \`{a} $\omega^{(r)}$ via $\pi^r_1$. On \'{e}crit $Q_\tau P^rM$
pour le sous-fibr\'{e} du fibr\'{e} des connexions principales sur $P^rM$,
$QP^rM$, qui correspondent aux connexions sans torsion sur $P^1M$.

Pour chaque
immersion $\Phi:M\ra M'$ o\`{u} $M'$ est une vari\'{e}t\'{e} diff\'{e}rentiable de
dimension $m$, on associe l'application $P^r\Phi:P^rM\ra P^rM'$ d\'{e}finie
par $P^r\Phi\big(j^r_0(f)\big):= j^r_0\big(\Phi(f)\big)$.
Il s'ensuit que $(P^r\Phi,\underline{P^r\Phi}:=\Phi,\phi:=\id_{G^r_m})$
est un morphisme de fibr\'{e}s principaux. Dans le cas $r=1$, la base
$(e_1,\ldots,e_m)$ de l'espace tangent $T_xM$ est envoy\'{e}e \`{a} la base
$(T_x\Phi e_1,\ldots,T_x\Phi e_m)$ de $T_{\Phi(x)}M$.
Soit $\omega'$ une $1$-forme de
connexion sur $P^1M'$. Puisque $P^rM$ et $P^rM'$ ont le m\^{e}me groupe
structurel $G^r_m$ on peut
construire la $1$-forme de connexion retir\'{e}e $\omega:=
(P^1\Phi)^*\omega'$
\`{a} laquelle une connexion
$\nabla=:\Phi^*\nabla'$ dans le fibr\'{e} tangent $TM$ est associ\'{e}e. Il y a
une formule plus directe pour $\Phi^*\nabla'$: pour un champ de vecteurs
$X'$ sur $M'$ on rappelle le {\em champ de vecteurs retir\'{e} $\Phi^*X'$ sur $M$}
d\'{e}fini par $(\Phi^*X')_x:=(T_x\Phi)^{-1}\hspace{1mm}X'_{\Phi(x)}$ $\forall x\in M$;
et pour
un autre champ de vecteurs $Y'$ sur $M'$ il vient
\beq \lbl{EqConnRetTan}
   (\Phi^*\nabla')_{\Phi^*X'}\Phi^*Y':=\Phi^*\big(\nabla'_{X'}Y'\big).
\end{equation}
Au cas o\`{u} $\Phi$ est un
diff\'{e}omorphisme on a
$(\Phi^*\nabla')_XY=\Phi^*\big(\nabla'_{\Phi_*X}\Phi_*Y\big)$ quels que
soient les champs de vecteurs $X,Y\in\Ginf(TM)$ avec
$\Phi_*:=(\Phi^{-1})^*$.}

\section{Fibr\'{e}s et op\'{e}rateurs naturels} \lbl{SecApp2}

{\small
Un {\em foncteur de fibr\'{e}s} o\`{u} un {\em fibr\'{e} naturel}
(voir \cite{KMS93}, p. 138) est un foncteur covariant $F$ de ${\mathcal{M}f}_m$
dans $\mathcal{FM}$ tel que les conditions suivantes soient satisfaites:\\
($i$) $B\circ F = \mathrm{Id}_{\mathcal{M}f_m}$
(c.-\`{a}-d. $FM$ est l'espace total d'une vari\'{e}t\'{e} fibr\'{e}e $\pi_M:FM\ra M$ sur la
m\^{e}me vari\'{e}t\'{e} $M$).\\
 ($ii$) pour toute inclusion $i:U\ra M$ d'une sous-vari\'{e}t\'{e} ouverte
 on a $FU=FM|_U=p_M^{-1}(U)$ et $Fi$ est \'{e}gale \`{a} l'inclusion $p_M^{-1}(U)\subset
 FM$ ({\em localit\'{e}}).\\
($iii$) Si $f:P\times M\ra N$ est de classe $\Cinf$ telle que $f_p:=f(p,~):M\ra
N$ est un diffeomorphisme local quel que soit $p\in P$, alors
$\tilde{F}f:P\times FM\ra FN$ d\'{e}finie par $\tilde{F}f(p,~):= Ff_p$ est de
classe $\Cinf$ ({\em r\'{e}gularit\'{e}}).\\
Dans \cite{KMS93}, p.187 les auteurs montrent
que la condition de
r\'{e}gularit\'{e} (iii) est superflue.

On \'{e}crit $F_xM$ pour la fibre $p_M^{-1}(x)$ sur $x\in M$ et $F_x\Phi$
pour la restriction $(F\Phi)|_{F_xM}$.
La vari\'{e}t\'{e} fibr\'{e}e $FM$ est toujours un fibr\'{e} localement trivial sur $M$
\`{a} fibre type $S:=F_0\real^m$, voir \cite{KMS93}, p.139.
La condition de localit\'{e} ($ii$) implique que pour toute immersion
$\Phi:M\ra N$ (o\`{u} $\dim M =m = \dim N$) et pour tout ouvert $i:U\ra M$
on a $F(\Phi|_U)=(F\Phi)_{FM|_U}$ car $\Phi|_U=\Phi\circ i$.
Soit $x\in M$. Puisqu'il existent un voisinage ouvert $U$ de $x$ et un
voisinage $V$ de $y:=\Phi(x)$ tel que $\Phi|_U$ est un diff\'{e}omorphisme
de $U$ sur $V$, il vient que $F(\Phi|_U)$ est un diff\'{e}omorphisme de
$FU$ sur $FV$. En particulier, $F$ d\'{e}finit un foncteur de la cat\'{e}gorie
$\mathcal{M}f_m$ dans la cat\'{e}gorie $\mathcal{M}f_{m+\dim S}$. En outre,
il s'ensuit que $F_x\Phi$ est toujours un diff\'{e}omorphisme
de la fibre $F_xM$ sur la fibre $F_yN$.

Un foncteur de fibr\'{e}s est dit {\em d'ordre $r$}, $r\in\nat$, lorsque
la condition suivante est satisfaite pour toutes les vari\'{e}t\'{e}s $M,N$
de dimension $m$ et pour tous les morphismes $\Phi,\Phi':M\ra N$:
Pour tout $x\in M$, si les $r$-jets en $x$ de $\Phi$ et de $\Phi'$
co\"{\i}ncident, alors $F\Phi|_{F_xM}=F\Phi'|_{F_xM}$. L'objet $FM$ d'un
foncteur de fibr\'{e}s
d'ordre $r$ s'obtient toujours comme un fibr\'{e} associ\'{e} de $P^rM$, c.-\`{a} d.
$FM\cong P^rM\times_{G^r_m} S$ avec une action \`{a} gauche $\ell:G^r_m\times S\ra S$,
voir \cite{KMS93}, p.140. Les morphismes $F\Phi$ sont induits par $P^r\Phi$
c.-\`{a}.-d. $F\Phi\big(q_M(p,s)\big)=q_N\big(P^r\Phi(p),s\big)$ quels que soient
$p\in P^rM,s\in S$.

Il est
clair que tout $P^r$, vu comme foncteur avec les d\'{e}finitions ci-dessus, est un
foncteur de fibr\'{e}s d'ordre $r$. Ensuite, le fibr\'{e} des connexions
principales de $P^rM$, $QP^rM$, peut \^{e}tre vu en tant qu'objet du foncteur de
fibr\'{e}s $QP^r$ d'ordre $r+1$, voir \cite{KMS93}, p.162, ainsi que
$Q_\tau P^1M$ (pour les connexion sans torsion).

Soient $F$ et $F'$ deux foncteurs de fibr\'{e}s (d'ordre $r$ et $r'$).
Alors l'association $F\times
F'$ d\'{e}finie par $(F\times F')M:=\bigcup_{x\in M}F_xM\times F'_xM$ ($\subset
FM\times F'M$) et
$(F\times F')\Phi:=F\Phi\times F'\Phi|_{(F\times F')M}$ est un foncteur de
fibr\'{e}s (d'ordre $\max(r,r')$).

Un {\em foncteur de fibr\'{e}s vectoriels} $F$ est
un foncteur de fibr\'{e}s dont les objets $FM$ sont des fibr\'{e}s vectoriels et
les morphismes $F\Phi$ sont les morphismes de fibr\'{e}s vectoriels. Pour un
foncteur de fibr\'{e}s vectoriels d'ordre $r$ l'objet $FM$ est toujours un
fibr\'{e} associ\'{e} $P^rM\times_{G^r_m} V$ de $P^rM$ o\`{u} le groupe $G^r_m$ est
repr\'{e}sent\'{e} sur un espace vectoriel $V$, voir \cite{KMS93}, p.141.
On voit que le foncteur tangent $T$ et cotangent $T^*$ sont des exemples de
foncteurs de fibr\'{e}s vectoriels d'ordre $1$.

Soient $F$ et $F'$ deux
foncteurs de fibr\'{e}s vectoriels (d'ordre $r$ et $r'$). Alors
l'association $F\otimes F'$ d\'{e}fini par $(F\otimes F')M:=FM\otimes GM$
et $(F\otimes F')\Phi:=F\Phi\otimes F'\Phi$ est un foncteur de fibr\'{e}s
vectoriels (d'ordre $\max(r,r')$). En outre, l'association $Hom(F,F')$ d\'{e}finie par
le fibr\'{e} vectoriel $Hom(F,F')M$ dont la fibre en $x\in M$ est donn\'{e}e
par $Hom(F,F')_xM:=Hom(F_xM,F'_xM)$ et par l'application
$Hom(F,F')\Phi$ qui envoie $a_x\in Hom(F_xM,F'_xM)$ \`{a} $F'_x\Phi\circ a_x
\circ (F_x\Phi)^{-1}$ est \'{e}galement un foncteur de fibr\'{e}s vectoriels
(d'ordre $\max(r,r')$).

Soit $F$ un foncteur de fibr\'{e}s, $\Phi:M\ra M'$ une immersion
entre deux vari\'{e}t\'{e}s de dimension $m$ et soit $\varphi'$ une section de classe
$\Cinf$ du fibr\'{e} $FM'$. Alors la {\em section retir\'{e}e} (ou la
section pull-back) $\Phi^*\varphi'$ du fibr\'{e} $FM$ est d\'{e}finie par
\beq
   \big(\Phi^*\varphi'\big)(x)
    :=\big(F_x\Phi\big)^{-1}\hspace{1mm}\varphi'\big(\Phi(x)\big)~~~
                                       \forall x\in M.
\end{equation}

On va maintenant donner une description des op\'{e}rateurs diff\'{e}rentiels en
termes de fibr\'{e}s de jets: Soient $\tau:E\ra M$ et $\bar{\tau}:\bar{E}\ra N$
deux fibr\'{e}s sur deux vari\'{e}t\'{e}s $M$ et $N$ de dimension $m$. On note par
$J^kE$ la $k^{\mathrm{\grave{e}me}}$ prolongation jet du fibr\'{e} $E$,
c.-\`{a}-d. l'ensemble de tous les $k$-jets des sections locales de $E$,
voir \cite{KMS93}, p. 124. Au cas o\`{u} $E$ est un fibr\'{e} vectoriel alors
$J^kE$ a la structure d'un fibr\'{e} vectoriel sur $M$.
Soit $\Phi:E\ra \bar{E}$ un morphisme de fibr\'{e}s
(vectoriels) tel que l'application induite $\underline{\Phi}$ est une
immersion. En consid\'{e}rant les jets d'ordre $k$ on obtient un morphisme
de vari\'{e}t\'{e}s fibr\'{e}es $J^k\Phi:J^kE\ra J^k\bar{E}$, voir \cite{KMS93},
p.124. Alors $J^k$ est un foncteur dans la cat\'{e}gorie $\mathcal{FM}_m$
des vari\'{e}t\'{e}s fibr\'{e}es (m\^{e}me des fibr\'{e}s vectoriels) dont les bases sont de
dimension $m$ et les
morphismes induits sur la base, $\underline{\Phi}$, sont des immersions.
Il s'ensuit que
pour tout foncteur de fibr\'{e}s $F$ d'ordre $r$ la composition
$J^k\circ F$ est un foncteur de fibr\'{e}s d'ordre $k+r$, voir
\cite{KMS93}, p.144.

Soient $F$ et $F'$ deux foncteurs de fibr\'{e}s. Un {\em op\'{e}rateur naturel}
$A:F\leadsto F'$ est d\'{e}fini par la suivante: Pour chaque
vari\'{e}t\'{e} diff\'{e}rentiable $M$ de dimension $m$ il y a une application
$A_M:\Ginf(FM)\ra \Ginf(F'M)$ telle que pour  toute section
$\varphi\in\Ginf(FM)$ on a les condition suivantes:\\
($i$) Pour tout diff\'{e}omorphisme $\Phi:M\ra N$ il vient
\[
    A_N(F\Phi\circ \varphi\circ\Phi^{-1})
              = F'\Phi\circ A_M\varphi\circ \Phi^{-1}.
\]
($ii$) Pour tout ouvert $U\subset M$: $A_U(\varphi|_U)=(A_M\varphi)|_U$
(localit\'{e}).\\
($iii$) Pour toute vari\'{e}t\'{e} $M_1$ et pour toute application
$\psi:M_1\times M \ra FM$ de classe $\Cinf$ telle que pour tout $z\in M_1$
l'application $\psi_z:=\psi(z,~)$ est une section de $FM$ il vient que
l'application $M_1\times M \ra F'M:(z,x)\mapsto \big(A_M(\psi_z)\big)(x)$
soit de classe $\Cinf$ (r\'{e}gularit\'{e}).\\
Un op\'{e}rateur naturel $A:F\leadsto F'$ est dit {\em d'ordre $k$} lorsque
pour toute vari\'{e}t\'{e} $M$ de dimension $m$ il existe une application
$\mathcal{A}_M:J^kFM\ra F'M$ (dite {\em l'application associ\'{e}e \`{a} $A_M$})
telle que $A_M(\varphi)=\mathcal{A}_M\big(j^r(\varphi)\big)$.
Les $\mathcal{A}_M$ sont en bijection avec les transformations naturelles
du foncteur $J^k\circ F$ \`{a} $F'$, voir \cite{KMS93}, p.144. Pour les
op\'{e}rateurs naturels d'ordre $k$ on peut reformuler les conditions ($i$) et
($ii$) par l'\'{e}quation suivante: pour toutes vari\'{e}t\'{e}s $M,N$ de dimension
$m$, pour toute immersion $\Phi:M\ra N$ et pour toute section $\varphi'$
dans $\Ginf(FN)$
\beq
       A_M(\Phi^*\varphi')=\Phi^*(A_N\varphi').
\end{equation}

Dor\'{e}navant, soit $M=N$. Tout op\'{e}rateur diff\'{e}rentiel
$D\in \Dop(E,E')$ d'ordre $r$ correspond uniquement \`{a} une section
$\check{D}$ du fibr\'{e} $Hom(J^rE,E')$ en d\'{e}finissant
$D\varphi:=\check{D}(j^r\varphi)$ quel que soit la section $\varphi\in
\Ginf(E)$. Il s'ensuit que si $E$ et $E'$ sont les objets
$FM$ et $F'M$ des foncteurs de fibr\'{e}s vectoriels $F$ et $F'$, alors le fibr\'{e}
vectoriel $Hom(J^rE,E')$ est l'objet $Hom(J^r\circ F,F')M$ du
foncteur de fibr\'{e}s vectoriels $Hom(J^r\circ F,F')$.}

\section{D\'{e}monstration de la proposition \ref{PRelNatConn}} \lbl{SecApp3}

{\small
 1. La vari\'{e}t\'{e}
 $\tilde{M}^a$ est un ouvert de l'espace total $E^a$ du fibr\'{e}
 $|\Lambda^m T^*M|^a$, et le sous-fibr\'{e} horizontal $H\tilde{M}^a$ de
 $T\tilde{M}^a$ d\'{e}finie par la $1$-forme de connexion $\omega^a$ co\"{\i}ncide
 avec la restriction du sous-fibr\'{e} horizontal $HE^a$ de $TE^a$ induit
 par celui d\'{e}fini par la $1$-forme de connexion $\omega$ sur $P^1M$.
 Donc les restrictions \`{a} $\tilde{M}^a$ des rel\`{e}vements horizontaux
 $X^{\mathbf{h}}$ et $Y^{\mathbf{h}}$ de $X$ et $Y$ \`{a} $E^a$
 co\"{\i}ncident avec les rel\`{e}vements horizontaux $X^h$ et $Y^h$.
 Soient $\varphi,\varphi'\in\Ginf(|\Lambda^m T^*M|^a)$ et
 $\varphi^v,{\varphi'}^v$ leurs rel\`{e}vements verticaux en tant que champs
 de vecteurs verticaux sur $E^a$. Il s'ensuit que les formules suivantes
 d\'{e}finissent une connexion $\check{\nabla}$ dans le fibr\'{e} tangent de $E^a$,
 voir par exemple \cite{KMS93}, p.410, Proposition:
 \[
   \begin{array}{rcl}
    \check{\nabla}_{X^{\mathbf{h}}} Y^{\mathbf{h}}=(\nabla_XY)^{\mathbf{h}},
            & ~~~ &
            \check{\nabla}_{X^{\mathbf{h}}}{\varphi'}^v =(\nabla_X\varphi')^v,
                         \\
    \check{\nabla}_{\varphi^v}{Y^{\mathbf{h}}}=0,
            & ~~~ &
            \check{\nabla}_{\varphi^v}{\varphi'}^v=0.
   \end{array}
 \]
 Puisque la torsion de $\nabla$ s'annule la torsion de $\check{\nabla}$ restreinte
 \`{a} $\tilde{M}^a$ est \'{e}gale \`{a}
 \beas
     Tor_{\check{\nabla}}(X^{\mathbf{h}},Y^{\mathbf{h}})
       & = & [X,Y]^{\mathbf{h}} - [X^{\mathbf{h}},Y^{\mathbf{h}}] \\
       & = & [X,Y]^{h} - [X^{h},Y^{h}]
                           \stackrel{(\ref{EqCrLieHorHorMtilde})}{=}
                           -a(\mathrm{tr}R)(X,Y)\mathsf{E}, \\
     Tor_{\check{\nabla}}(X^{\mathbf{h}},{\varphi'}^v) &
               \stackrel{(\ref{EqCrLieHorVer})}{=} & 0, \\
     Tor_{\check{\nabla}}({\varphi}^v,{\varphi'}^v)    &
               \stackrel{(\ref{EqCrLieVerVer})}{=} & 0.
 \eeas
 Par cons\'{e}quent, la connexion $\hat{\nabla}$ d\'{e}finie par
 $\hat{\nabla}_ZZ':=\check{\nabla}_ZZ'-\frac{1}{2}Tor_{\check{\nabla}}(Z,Z')$
 est sans torsion. On voit ais\'{e}ment que
 \[
    \hat{\nabla}_{X^{\mathbf{h}}}Y^{\mathbf{h}}=
       (\nabla_XY)^{\mathbf{h}} +\frac{a}{2}(\mathrm{tr}R)(X,Y)\mathsf{E}
 \]
 tandis que les trois autres termes qui d\'{e}finissent $\check{\nabla}$ ne
 changent pas. L'\'{e}quation
 $\hat{\nabla}_{\varphi^v}X^{\mathbf{h}}=\check{\nabla}_{\varphi^v}X^{\mathbf{h}}=0$
 implique que $\hat{\nabla}_{V}X^{\mathbf{h}}=0$ pour tout champ de vecteurs vertical,
 alors, en particulier, pour $V=\mathsf{E}$ on a
 $0=\hat{\nabla}_{\mathsf{E}}X^{\mathbf{h}}=\hat{\nabla}_{\mathsf{E}}X^{h}.$
 Puisque les rel\`{e}vements horizontaux sont invariants on a
 \beas
  0 & = & [X^h, \mathsf{E}]=[X^{\mathbf{h}},\mathsf{E}]
              = \hat{\nabla}_{X^{\mathbf{h}}}\mathsf{E}
               -\hat{\nabla}_{\mathsf{E}}X^{\mathbf{h}}
     =   \hat{\nabla}_{X^{\mathbf{h}}}\mathsf{E} - 0
              = \hat{\nabla}_{X^{h}}\mathsf{E}.
 \eeas
 Finalement, l'\'{e}quation $\hat{\nabla}_{\varphi^v}{\varphi'}^v
 =\check{\nabla}_{\varphi^v}{\varphi'}^v=0$
 implique que $\hat{\nabla}_{V}{\varphi'}^v=0$ pour tout champ de vecteurs vertical,
 alors, en particulier, pour $V=\mathsf{E}$ on a
 $0=\hat{\nabla}_{\mathsf{E}}{\varphi'}^v.$
 Le flot $F_t$ de $\varphi^v$ est donn\'{e} par la translation
 $F_s(y)=y+s\varphi\big(\tau^a(y)\big)$, tandis que le flot du champ
 d'Euler s'\'{e}crit $G_t(y)=e^ty$ quel que soit $y\in E^a$. Par cons\'{e}quent,
 $(G_{-t}\circ F_{s}\circ G_t)(y)=y+se^{-t}\varphi\big(\tau^a(y)\big)$,
 alors le champ retir\'{e} est de la forme $G_t^*(\varphi^v)=e^{-t}\varphi^v$,
 donc
 \[
   \hat{\nabla}_{\varphi^v}\mathsf{E}=
    \hat{\nabla}_{\varphi^v}\mathsf{E}-\hat{\nabla}_{\mathsf{E}}{\varphi}^v
    =[\varphi^v,\mathsf{E}]=
        -\frac{d}{dt}\big( G_t^*(\varphi^v)\big)|_{t=0}
        =-\frac{d}{dt}\big( e^{-t}(\varphi^v)\big)|_{t=0}=\varphi^v.
 \]
 On voit que l'\'{e}quation en r\'{e}sultant
 $\hat{\nabla}_{\varphi^v}\mathsf{E}=\varphi^v$ reste vraie pour tout
 champ de vecteurs vertical $V$, alors $\hat{\nabla}_V\mathsf{E}=V$, en
 particulier, pour $V=\mathsf{E}$ il vient
 $\hat{\nabla}_{\mathsf{E}}\mathsf{E}=\mathsf{E}$.
 Alors la connexion $\hat{\nabla}$, restreinte \`{a} $\tilde{M}^a$, se r\'{e}crit
 de la fa\c{c}on suivante:
 \[
    \begin{array}{rcl}
    \hat{\nabla}_{X^{h}} Y^{h}=(\nabla_XY)^{h}
       +\frac{a}{2}\big({\tau^a}^*(\mathrm{tr}R)(X,Y)\big)\mathsf{E},
            & ~~~ &
            \hat{\nabla}_{X^{h}}\mathsf{E} =0,
                         \\
    \hat{\nabla}_{\mathsf{E}}{X^{h}}=0,
            & ~~~ &
            \hat{\nabla}_{\mathsf{E}}\mathsf{E}=\mathsf{E}.
   \end{array}
 \]
 Toute autre connexion sans torsion dans le fibr\'{e} tangent de est \'{e}gale
 \`{a} la somme de $\hat{\nabla}$ et d'un champ de tenseur sym\'{e}trique
 $W$ appartenant \`{a}
 $\Ginf(T\tilde{M}^a\otimes S^2T^*\tilde{M}^a)$.
 Puisque les valeurs de $X^h$ et et $\mathsf{E}$ en tout $y\in\tilde{M}^a$
 engendrent l'espace tangent $T_y\tilde{M}^a$ ce champ $W$ peut \^{e}tre
 d\'{e}finie sur les champs de vecteurs $X^h$ et et $\mathsf{E}$. En
 d\'{e}finissant
 $W(X^h,Y^h):=\mu{\tau^a}^*\big(Ric(X,Y)+Ric(Y,X)\big)\hspace{1mm}\mathsf{E}$,
 $W(X^h,\mathsf{E}):=\nu X^h$ et
 $W(\mathsf{E},\mathsf{E}):=(\rho-1)\mathsf{E}$, on voit que
 $\bar{\nabla}$ est une connexion sans torsion bien d\'{e}finie dans le fibr\'{e}
 tangent de $\tilde{M}^a$.

 Pour l'invariance de $\bar{\nabla}$
 il suffit de v\'{e}rifier
 que le champ de tenseurs $L_{\mathsf{E}}\bar{\nabla}$ s'annule sur les
 rel\`{e}vements horizontaux et sur le champ d'Euler: masi ceci se voit
 facilement avec la formule (\ref{EqDefLENabla})
 car $[\mathsf{E},X^h]=0=[\mathsf{E},\mathsf{E}]$ et les fonctions
 ${\tau^a}^*(\mathrm{tr}R)(X,Y)$ et $\mu{\tau^a}^*\big(Ric(X,Y)+Ric(Y,X)\big)$
 sont invariantes par le flot de $\mathsf{E}$
 parce qu'elles ne d\'{e}pendent pas de la fibre sur tout point $x\in M$.

 Soit $\Phi:M\ra M'$ une immersion dans une autre vari\'{e}t\'{e} diff\'{e}rentiable
 $M$ de dimension $m$. Soit $\nabla'$ une connexion sans torsion dans le
 fibr\'{e} tangent de $M'$ et soit $\nabla:=\Phi^*\nabla'$ la connexion retir\'{e}e.
 Pour montrer la naturalit\'{e} de
 $\nabla\mapsto\hat{\nabla}$ on observe que la diff\'{e}rence de connexions
 $\overline{\Phi^*\nabla'}-\tilde{\Phi}^{a*}\hat{\nabla}$ d\'{e}finit un
 champ de tenseurs sym\'{e}trique $B$ dans
 $\Ginf(T\tilde{M}^a\otimes S^2T^*\tilde{M}^a)$. Donc,
 pour v\'{e}rifier qu'il s'annule, il suffit de le v\'{e}rifier sur les champs de
 vecteurs retir\'{e}s $\tilde{\Phi}^{a*}{X'}^{h'}$ (o\`{u} $X'$ est un champ de vecteurs
 sur $N$ et $(~)^{h'}E$ d\'{e}signe le rel\`{e}vement horizontal par rapport \`{a} $\nabla'$)
 et le champ d'Euler $\mathsf{E}'$ de $\tilde{M}^a$ retir\'{e},
 $\tilde{\Phi}^{a*}\mathsf{E}'$. Puisque $\tilde{\Phi}^a$ est un morphisme
 de fibr\'{e}s principaux il vient $\tilde{\Phi}^{a*}\mathsf{E}'=\mathsf{E}$
 et $\tilde{\Phi}^{a*}{X'}^{h'}=(\Phi^*X')^h$ o\`{u} $(~)^h$ d\'{e}signe le
 rel\`{e}vement horizontal par rapport \`{a} $\nabla$. De plus, le tenseur de
 courbure $R$ de $\nabla$ est \'{e}gal au tenseur de courbure retir\'{e}
 $R:=\Phi^*R'$ du tenseur de courbure $R'$ de $\nabla'$. Il en est de m\^{e}me
 pour le tenseur $\mathrm{tr}R$ et le tenseur de Ricci $Ric$. En notant
 $X:=\Phi^*X'$, $Y:=\Phi^*Y'$ pour un autre champ de vecteurs $Y'$ sur
 $N$ et ${\tau'}^a$ pour la projection $\tilde{N}^a\ra N$ on a d'une part:
 \beas
   \lefteqn{\overline{\Phi^*\nabla'}_{\tilde{\Phi}^{a*}{X'}^{h'}}
             \tilde{\Phi}^{a*}{Y'}^{h'}
        =  \bar{\nabla}_{X^h}Y^h }\\
       & = & (\nabla_XY)^h
              +\frac{a}{2}\big({\tau^a}^*(\mathrm{tr}R)(X,Y)\big)\mathsf{E}
       +\mu{\tau^a}^*\big(Ric(X,Y)+Ric(Y,X)\big)\hspace{1mm}
                            \mathsf{E}
 \eeas
 et d'autre part
 \beas
   \lefteqn{\big(\tilde{\Phi}^{a*}\hat{\nabla'}\big)_{\tilde{\Phi}^{a*}{X'}^{h'}}
         \tilde{\Phi}^{a*}{Y'}^{h'}
        =
         \tilde{\Phi}^{a*}\big(
          \hat{\nabla'}_{{X'}^{h'}}{Y'}^{h'}\big)} \\
       & = &
         \tilde{\Phi}^{a*}
            \big((\nabla'_{X'}Y')^{h'}
              +\frac{a}{2}\big({{\tau'}^a}^*(\mathrm{tr}R')(X',Y')\big)
               \mathsf{E}'\big) \\
       &   &~~~~~~~~~+\mu \tilde{\Phi}^{a*}\big(
            {{\tau'}^a}^*\big(Ric'(X',Y')+Ric'(Y',X')\big)\mathsf{E}'\big)
            \\
       & = & \big(\Phi^*(\nabla'_{X'}Y')\big)^h
             + \frac{a}{2}{\tau^a}^*
                 \big(\Phi^*\big((\mathrm{tr}R')(X',Y')\big)\big)
                 \tilde{\Phi}^{a*}\mathsf{E}' \\
       &   &~~~~~~~~~+\mu
            {{\tau}^a}^*\Phi^*\big(Ric'(X',Y')+Ric'(Y',X')\big)
                  \tilde{\Phi}^{a*} \mathsf{E}' \\
       & = &(\nabla_XY)^h
              +\frac{a}{2}\big({\tau^a}^*(\mathrm{tr}R)(X,Y)\big)\mathsf{E}
       +\mu{\tau^a}^*\big(Ric(X,Y)+Ric(Y,X)\big)\hspace{1mm}
                            \mathsf{E}.
 \eeas
 Pour les autres combinaisons $(\tilde{\Phi}^{a*}{X'}^{h'},\mathsf{E})$ et
 $(\mathsf{E},\mathsf{E})$ on fait un calcul analogue.\\
 Alors l'association $\nabla\mapsto \bar{\nabla}$ est un rel\`{e}vement
 naturel de connexions quels que soient les nombres r\'{e}els $\mu,\nu,\rho$.
 En particulier, $\nabla\mapsto\hat{\nabla}$ l'est aussi.

 2. Toute autre connexion sans torsion $\nabla^{\circ}$ sur $\tilde{M}^a$
 qui d\'{e}finit un
  rel\`{e}vement naturel de connexions est \'{e}gal \`{a} la somme de $\bar{\nabla}$
 et d'un champ de tenseur sym\'{e}trique $C$ appartenant \`{a}
 $\Ginf(T\tilde{M}^a\otimes S^2T^*\tilde{M}^a)$. $C$ peut \^{e}tre d\'{e}fini sur
 les rel\`{e}vements horizontaux et le champ d'Euler et doit \^{e}tre invariant.
 Soit $G_t$ le flot du champ d'Euler $\mathsf{E}$, et soient $Z,Z'$ des
 champs de vecteurs invariants sur $\tilde{M}^a$ (par exemple $X^h$ ou
 $\mathsf{E}$). L'invariance de $C$ veut dire
 \[
     T_yG_t\hspace{1mm}C_y(Z_y,Z'_y)\stackrel{!}{=}
      C_{ye^t}(T_yG_t\hspace{1mm}Z_y,T_yG_t\hspace{1mm}Z'_y)
      =C_{ye^t}(Z_{ye^t},Z'_{ye^t})
 \]
 quels que soient $y\in \tilde{M}^a,t\in\real$. Puisque $\tau^a\circ
 G_t=\tau^a$, l'application de $T_{ye^t}\tau^a$ \`{a} cette \'{e}quation montre
 que $T_y\tau^a C_y(Z_y,Z'_y)\in T_xM$ ne d\'{e}pend pas de la fibre $y$ sur
 $x=\tau^a(y)$. Pour le choix $Z=X^h,Z'=Y^h$ on voit que cette expression
 d\'{e}finit un champ de tenseur $C_1[\nabla]$ appartenant \`{a} $\Ginf(TM\otimes S^2T^*M)$
 par $C_1[\nabla]_{\tau^a(y)}(X_{\tau^a(y)},Y_{\tau^a(y)}):=T_y\tau^a
 C_y(X^h_y,Y^h_y)$. Pour le choix $Z=X^h, Z'=\mathsf{E}$ on obtient de la
 m\^{e}me fa\c{c}on un champ de tenseur $C_2[\nabla]$ dans
 $\Ginf\big(Hom(TM,TM)\big)$, et pour le choix $Z=Z'=\mathsf{E}$ on
 obtient un champ de vecteurs $C_3[\nabla]$ sur $M$. On voit que les expressions
 $C(X^h,Y^h)-\big(C_1[\nabla](X,Y)\big)^h$,
 $C(X^h,\mathsf{E})-\big(C_2[\nabla](X)\big)^h$
 et $C(\mathsf{E},\mathsf{E})-C_3[\nabla]^h$ sont des champs de vecteurs verticaux
 invariants, alors ils sont des multiples invariants du champ d'Euler.
 A l'aide d'un raisonnement analogue, on trouve une $2$-forme sym\'{e}trique
 $C_4[\nabla]$ appartenant \`{a} $\Ginf(S^2T^*M)$, une $1$-forme $C_5[\nabla]$
 appartenant \`{a}
 $\Ginf(T^*M)$ et une fonction $C_6[\nabla]\in\Cinf(M,\real)$ telles que la
 connexion $\nabla^\circ$ est de la forme suivante
 \beas
     {\nabla^\circ}_{X^h}Y^h & = & \hat{\nabla}_{X^h}Y^h
                                    + \big(C_1[\nabla](X,Y)\big)^h
                + {\tau^a}^*C_4[\nabla](X,Y)\hspace{1mm}\mathsf{E}
                                       \\
     {\nabla^\circ}_{X^h}\mathsf{E}={\nabla^\circ}_{\mathsf{E}}X^h
                             & = & \hat{\nabla}_{X^h}\mathsf{E}
                                   + \big(C_2[\nabla](X)\big)^h
                                   + {\tau^a}^*C_5[\nabla](X)\hspace{1mm}\mathsf{E}
                                   \\
     {\nabla^\circ}_{\mathsf{E}}\mathsf{E}
                             & = & \hat{\nabla}_{\mathsf{E}}\mathsf{E}
                                    + C_3[\nabla]^h
                                    + {\tau^a}^*C_6[\nabla]\hspace{1mm}\mathsf{E}
 \eeas
 Quand on examine la naturalit\'{e} de l'association
 $\nabla\mapsto \nabla^\circ$ comme dans la premi\`{e}re partie de cette
 d\'{e}monstration on voit ais\'{e}ment que --gr\^{a}ce \`{a} la naturalit\'{e} de
 $\nabla\mapsto \hat{\nabla}$--que les associations $\nabla\mapsto
 C_k[\nabla]$, $1\leq l\leq 6$, d\'{e}finissent des op\'{e}rateurs naturels
 \[
  \begin{array}{rcl}
    C_1:Q_\tau P^1M\leadsto T\otimes S^2T^*,
              & C_2:Q_\tau P^1M\leadsto Hom(T,T), &
                     C_3:Q_\tau P^1M\leadsto T, \\
    C_4:Q_\tau P^1M\leadsto S^2T^*, & C_5:Q_\tau P^1M\leadsto T^*, &
                C_6: Q_\tau P^1M\leadsto S^0T^*.
  \end{array}
 \]
 Dans \cite{KMS93}, p.239, Example 28.7, on trouve la d\'{e}monstration pour
 le fait que $C_4$ est proportionnel \`{a} la partie sym\'{e}trique du tenseur de
 Ricci. Les autres cas se font d'une mani\`{e}re similaire: on note d'abord
 que tous les op\'{e}rateurs $C_l$, $1\leq l\leq 6$ sont d'ordre fini $r$, voir
 \cite{KMS93}, p. 207, Proposition 23.5 et p. 208, 23.6 Examples. Un tel
 op\'{e}rateurs est toujours d\'{e}termin\'{e}s par une application $f$ de classe $\Cinf$ de
 l'espace affine $J^rQ_\tau P^1_0\real^m$ dans un $GL(m,\real)$-module $V_l$
 $1\leq l\leq 6$, qui soit $G^{r+2}_m$-\'{e}quivariante (o\`{u} $g\in G^{r+2}_m$ agit
 sur $V_l$ via la projection de jets $G^{r+2}_m\ra G^1_m=GL(m,\real)$),
 voir \cite{KMS93}, p.145, Theorem. Dans notre cas, on a
 $V_1=\real^m\otimes S^2{\real^m}^*$, $V_2=Hom(\real^m,\real^m)$,
 $V_3=\real^m$, $V_4=S^2{\real^m}^*$, $V_5={\real^m}^*$ et $V_6=\real$
 (module trivial). D'abord on v\'{e}rifie l'homog\'{e}n\'{e}it\'{e} des $f_l$ \`{a} l'aide des
 homoth\'{e}ties $s\mathbf{1}\in GL(m,\real)\in G^{r+2}_m$ quel que soit
 $s\in \real$:
 \[
     f_l(s\Gamma_0,s^2\Gamma_1,\ldots,s^{r+1}\Gamma_r)
        =s^{a_l}f_l(\Gamma_0,\Gamma_1,\ldots,\Gamma_r)
 \]
 o\`{u} $(\Gamma_0,\Gamma_1,\ldots,\Gamma_r)$ d\'{e}signe les valeurs d'un jet
 dans $J^rQ_\tau P^1_0\real^m$.
 On a $a_1=1$, $a_2=0$, $a_3=-1$, $a_4=2$, $a_5=1$ et $a_6=0$. Dans les
 cas $l=2$ et $l=6$ il vient d'apr\`{e}s le th\'{e}or\`{e}me de la fonction homog\`{e}ne
 (voir \cite{KMS93}, p.213) que $C_2$ et $C_6$ sont d'ordre $0$ et que
 $f_2$ et $f_6$ sont constantes \`{a} valeurs
 $GL(m,\real)$-invariantes, ce qui nous ne donne que
 $C_2[\nabla]=\nu\mathbf{1}$ (le champs de l'homomorphisme identique)
 et $C_6[\nabla]=\rho 1$ (la fonction constante \'{e}gale \`{a} $1$). Pour $l=3$
 la seule solution est \'{e}gale $f_4=0$, d'o\`{u} $C_3[\nabla]=0$. Pour $l=1$ et
 $l=5$ il s'ensuit que $r=0$ et $f_1$, $f_5$ sont homog\`{e}nes d'ordre $1$,
 donc lin\'{e}aires. Gr\^{a}ce \`{a} l'invariance de $f_1$ et $f_5$ par le nilradical
 du groupe $G^2_m$ (qui agit par translation) il s'ensuit pour le jet d'une
 connexion sans torsion $\Gamma_0\in\real^m\otimes S^2{\real^m}^*$
 que $f_l(\Gamma_0)=f_l(0)=0$ pour $l=1,5$. Alors $C_1[\nabla]=0$ et
 $C_5[\nabla]=0$. L'\'{e}nonc\'{e} sur $C_4$ est d\'{e}duit du th\'{e}or\`{e}me 28.6, p.238,
 \cite{KMS93} de la r\'{e}duction aux sous-espaces de courbure. Ainsi
 l'unicit\'{e} de la famille \`{a} trois param\`{e}tres $\mu,\nu,\rho$ de rel\`{e}vements
 naturels de connexions sans torsion pour le foncteur $\tilde{F}^a$ est
 d\'{e}montr\'{e}e.
}

\end{appendix}

\end{document}